%% file: fi23_R2_rev2.tex
\documentclass[12pt]{article}
\usepackage{amsmath,amssymb,amsthm,amsfonts,hhline,color}
\usepackage{bbm}
\usepackage[titletoc,title]{appendix}

\usepackage[dvipdfmx]{hyperref}
\usepackage{hyperref}

\usepackage{comment}
\usepackage{color}
\usepackage[dvipdfmx]{graphicx}

\newcommand{\ds}{\displaystyle}
\newcommand{\fr}[2]{\frac{#1}{#2}}
\newcommand{\dfr}[2]{\dfrac{#1}{#2}}
\newcommand{\cd}{\cdot}
\newcommand{\cds}{\cdots}
\newcommand{\dsum}{\displaystyle \sum}

\renewcommand{\l}{\left}
\renewcommand{\r}{\right}

\newcommand{\q}{\quad}
\newcommand{\qq}{\qquad}

\newcommand{\la}{\langle}
\newcommand{\ra}{\rangle}

\newcommand{\abs}[1]{\lvert{#1}\rvert}

\newcommand{\Z}{\mathbb{Z}}
\newcommand{\C}{\mathbb{C}}
\newcommand{\R}{\mathbb{R}}

\newcommand{\M}{\mathbb{M}}

\newcommand{\aut}{\mathrm{Aut}}
\newcommand{\wt}{\mathrm{wt}}
\newcommand{\tr}{\mathrm{tr}}
\renewcommand{\hom}{\mathrm{Hom}}

\newcommand{\al}{\alpha}

\newcommand{\Span}{\mathrm{Span}}

\newcommand{\w}{\omega}
\newcommand{\vac}{\mathbbm{1}}
\newcommand{\irr}{\mathrm{Irr}}
\newcommand{\Fi}{\mathrm{Fi}}
\newcommand{\mymid}{\,|\,}

\makeatletter
\@addtoreset{equation}{section}

\makeatother

\theoremstyle{plain}

\newtheorem{thm}{Theorem}[section]
\newtheorem{prop}[thm]{Proposition}
\newtheorem{lem}[thm]{Lemma}
\newtheorem{cor}[thm]{Corollary}

\theoremstyle{definition}
\newtheorem{df}[thm]{Definition}
\newtheorem{conj}[thm]{Conjecture}

\newtheorem{nota}[thm]{Notation}
\newtheorem{rem}[thm]{Remark}

\newcommand{\sfr}[2]{\leavevmode\kern-.05em
  \raise.5ex\hbox{\the\scriptfont0 #1}\kern-.1em
  /\kern-.15em\lower.25ex\hbox{\the\scriptfont0 #2}\kern.02em}

\newcommand{\shf}{\sfr{1}{2}}
\DeclareMathOperator*{\tensor}{\otimes}
\DeclareMathOperator*{\fusion}{\boxtimes}

\newcommand{\pf}{\noindent {\bf Proof:}\q }

\newcommand{\com}{\mathrm{Com}}

\newcommand{\dih}[2]{\mathrm{DIH}_{#1}(#2)}
\renewcommand{\ker}{\mathrm{Ker}}
\renewcommand{\o}{\mathrm{o}}
\newcommand{\Co}{\mathrm{Co}}
\newcommand{\VF}{V\!\! F^\natural}

\title{The Conway-Miyamoto correspondences \\
for the Fischer 3-transposition groups}

\author{
   Ching Hung Lam\footnote{Partially supported by MoST grant 104-2115-M-001-004-MY3 of Taiwan.}%
  \medskip\\
  {\small \it Institute of Mathematics, Academia Sinica, Taipei, Taiwan 10617}\\
 {\small \it  and National Center for Theoretical Sciences, Taipei,  Taiwan}\\
  {\small e-mail: {\tt chlam@math.sinica.edu.tw}}
  \bigskip\\
  Hiroshi Yamauchi\footnote{Partially supported by JSPS Grant-in-Aid for Young Scientists 
  (B) No 24740027. }
  \medskip\\
  {\small \it Department of Mathematics,
  Tokyo Woman's Christian University}\\
  {\small \it 2-6-1 Zempukuji, Suginami-ku, Tokyo 167-8585, Japan}\\
  {\small e-mail: \texttt{yamauchi@lab.twcu.ac.jp}}
}
\date{}

\begin{document}

\maketitle

\begin{abstract}
In this paper, we present a general construction of 3-transposition groups as automorphism 
groups of vertex operator algebras.
Applying to the moonshine vertex operator algebra, 
we establish the Conway-Miyamoto correspondences between Fischer 3-transposition groups 
$\Fi_{23}$ and $\Fi_{22}$ and $c=25/28$ and $c=11/12$ Virasoro vectors of subalgebras of 
the moonshine vertex operator algebra.
\end{abstract}

\pagestyle{plain}

\baselineskip 6mm

\tableofcontents

\section{Introduction}

A \emph{$3$-transposition group} is a pair $(G,D)$ of a group $G$ and a set of 
involutions $D$ of $G$, called 3-transpositions, such that $G=\la D\ra$, $D=D^G$ and  
$\abs{ab}\leq 3$ for any $a$, $b\in D$. 
A subgroup $H$ of $G$ is called a \emph{$D$-subgroup} if $H=\la H\cap D\ra$ 
which is again a 3-transposition group.
Let $a\in D$ and set $D_a=\{ x\in D \mid ax=xa\}$.
The $D$-subgroup $G_a=\la D_a\ra$ is uniquely determined up to conjugation 
if $D$ is a conjugacy class in $G$. 
We call the 3-transposition group $G_a/\la a\ra$ an \emph{inductive structure} of $G$.
Given a 3-transposition group $(G,D)$, the intersection $D\cap P$ of $D$ with a Sylow 
2-subgroup $P$ of $G$ is called a \emph{basic} set which is a maximal subset of $D$ 
consisting of mutually commutative elements.
The size of a basic set is called the \emph{width} of $G$.
Since $D^G=D$, the width is an \emph{invariant} of a 3-transposition group.
3-transposition groups were first studied by Fischer \cite{Fi} (see also \cite{As,CH}).   
Fischer also discovered the Fischer groups $\Fi_{22}$, $\Fi_{23}$ and $\Fi_{24}$ in the 
classification of almost simple 3-transposition groups.
It turns out that Fischer groups $\Fi_{22}$ and $\Fi_{23}$ are obtained by taking 
inductive structures of $\Fi_{23}$ and $\Fi_{24}$, respectively, 
whereas the largest Fischer group $\Fi_{24}$ is best understood via the inclusion 
of its cover $3.\Fi_{24}$ into the Monster (cf.~Section 11.34 of \cite{As}).
The subscripts of $\Fi_{22}$, $\Fi_{23}$ and $\Fi_{24}$ indicate their widths,  
and the normalizers of the elementary abelian subgroups generated by basic sets 
of $\Fi_{22}$, $\Fi_{23}$ and $\Fi_{24}$ are $2^{10}{:}M_{22}$, $2^{11}{\cd}M_{23}$ and 
$2^{12}{\cd}M_{24}$, respectively, where $M_{22}, M_{23},M_{24}$ are sporadic Mathieu groups 
(cf.~\cite{ATLAS}).

In this article, we will describe a general construction of 3-transposition groups 
as automorphism groups of vertex operator algebras (VOAs). 
We also try  to establish some non-trivial interconnections between the Fischer 
3-transposition groups and certain vertex operator subalgebras of the moonshine vertex 
operator algebra and to ensure that vertex operator algebras are canonical and natural 
objects for some of sporadic groups in the Monster family. 
Our methods are similar to that in \cite{HLY1,HLY2}. 
The main idea is to study certain vertex subalgebras generated by the so called 
Ising vectors (cf.~Definition \ref{df:2.1}) and try to analyze  the corresponding 
commutant subalgebras.
In particular, we construct two vertex operator algebras $\VF_{23}$ and $\VF_{22}$ 
as commutant vertex subalgebras of the famous moonshine VOA $V^\natural$. 
The Fischer $3$-transposition group $\Fi_{23}$ (resp.~$\Fi_{22}$) acts naturally on 
$\VF_{23}$ (resp.~$\VF_{22}$).  
We also prove the 3-transposition property using the theory of VOA and establish 
a natural correspondence between the 3-transpositions of $\Fi_{23}$ (resp.~$\Fi_{22}$) 
and some simple $c=25/28$ (resp.~$c=11/12$) Virasoro vectors in $\VF_{23}$ (resp.~$\VF_{22}$). 
We call such a correspondence the Conway-Miyamoto correspondence 
(see Definition \ref{df:5.1} for the precise meaning).   

The most prominent example of Conway-Miyamoto correspondences is given by 
the 2A-involutions of the Monster and the Ising vectors of the moonshine VOA.
Conway \cite{C} constructed a slightly modified version of the original Griess algebra 
\cite{GrO} of dimension $1+196883$ and described a 2A-involution of the Monster 
in terms of the adjoint action of an idempotent called the transposition axis. 
Miyamoto \cite{M1} showed that Conway's axis corresponds to a $c=1/2$ Virasoro vector, 
called an \emph{Ising vector}, of $V^\natural$ and he associated each Ising vector with  
an involution of $V^\natural$ based on the fusion rules of the Virasoro VOA, 
which may be viewed as a converse of Conway's axis; it turns out that there exists 
a one-to-one correspondence between the 2A-elements of the Monster and the Ising vectors 
of the moonshine VOA (cf.~\cite{C,M1,Ma1,Ho}). 
In \cite{HLY1,HLY2}, some similar correspondences between the 2A-elements of 
the Baby Monster and $c=7/10$ Virasoro vectors of $\sigma$-type of the Baby Monster VOA 
and that between the 2C-elements of the largest Fischer 3-transposition group and 
the $c=6/7$ Virasoro vectors of $\sigma$-type of the Fischer VOA were established. 
Both Baby Monster VOA and Fischer VOA are commutant subalgebras of the moonshine VOA 
and these correspondences are derived from the Monster.
Once a correspondence has been established, we can internally describe the corresponding 
automorphism group based on the VOA structure and study the group as an intrinsic 
symmetry of the corresponding VOA.

Now let us explain our setting. 
For a positive integer $n$, let $c_n=1- 6/(n+2)(n+3)$  
be the central charge of the unitary series of the Virasoro algebra.
If a VOA $V$ contains a simple $c=c_n$ Virasoro vector $e$,  then the linear map 
$\tau_e=(-1)^{4(n+2)(n+3)\o(e)}$ is well-defined, where $\o(e)=e_{(1)}$ is the zero-mode of $e$, 
and gives rise to an automorphism of $V$ called the \emph{Miyamoto involution} associated to $e$ 
(cf.~\cite{M1}).
The case $c_1=1/2$ in $V^\natural$ produces a 2A-element of the Monster.
Similarly, one can define another Miyamoto involution $\sigma_e$ on a certain subVOA of $V$ 
using a similar method (see Theorem \ref{thm:2.4} for detail). 
When $\tau_e$ is trivial and $\sigma_e$ is well-defined on the whole space $V$, 
such a Virasoro vector $e$ is said to be of $\sigma$-type on $V$.

A VOA $V$ is said to be \emph{of OZ-type} (which stands for one-dimensional, zero-dimensional) 
if $V_{n<0}=0$,  $V_0=\C \vac$ and $V_1=0$. 
We also assume that there exists a real form $V_\R$ of $V$ such that the invariant bilinear 
form on $V_\R$ is positive definite.
Our method is to study vertex subalgebras generated by Ising vectors and try to analyze 
the action of various types of Miyamoto involutions on some commutant subalgebras. 
In particular, some vertex subalgebras generated by two Ising vectors, 
which we call \emph{dihedral algebras}, studied in \cite{LM,LYY1,LYY2,S} play some 
important roles for our study. 
In \cite{S}, Griess algebras of dihedral subalgebras of a VOA of OZ-type are classified. 
It turns out that all of them are realizable inside the moonshine VOA (cf.~\cite{LM,LYY1,LYY2}). 
We will mainly consider the dihedral subalgebras of 2A, 3A and 6A-types, 
which are simply called the 2A, 3A and 6A-algebras, respectively (cf.~Theorem \ref{thm:2.6}). 

Let $E_V$ be the set of Ising vectors of $V_\R$.  
We assume that there exist $a$, $b\in E_V$ such that the dihedral subalgebra $\la a,b\ra$ 
generated by them is isomorphic to the 3A-algebra. 
In this case, the subgroup  $\la \tau_a,\tau_b\ra$ is isomorphic to $\mathrm{S}_3$ \cite{M3,S}. 
We set $I_{a,b}=\{ x\in E_V \mid (a\mymid x)=(b\mymid x)=2^{-5}\}$. 
Then  $x\in I_{a,b}$ if and only if both $\la a,x\ra$ and $\la b,x\ra$ are isomorphic to 
the 2A-algebra and each $\tau_x$, $x\in I_{a,b}$, centralizes $\la \tau_a,\tau_b\ra$ 
(cf.~\cite{M1}). 
We will prove in Corollary \ref{cor:3.7} that the group $G=\la \tau_x \mymid x\in I_{a,b}\ra$ 
is a 3-transposition subgroup of $C_{\aut(V)}(\tau_a,\tau_b)$, that is, if $x$, $y\in I_{a,b}$ 
then the order of $\tau_x\tau_y$ is bounded by three.
Since $G$  acts trivially on $\la a,b\ra$, it is natural to regard $G$ as a subgroup of 
the automorphism group of the commutant $\com_V\la a,b\ra$ of the 3A-subalgebra 
$\la a,b\ra$ in $V$.
As an application, we will obtain a realization of $\Fi_{23}$ as an automorphism group of 
the moonshine VOA.
Namely, if we apply Corollary \ref{cor:3.7} to $V^\natural$, we obtain the second largest 
Fischer group $\Fi_{23}\cong C_\M(\tau_a,\tau_b)$ acting on the commutant of the 3A-algebra 
$\la a,b\ra$ in $V^\natural$.

Our construction can be viewed as a sequel of the works in \cite{HLY1,HLY2}.
In \cite{HLY2}, the $W_3$-algebra $W_3(\sfr{4}{5})$ at $c=4/5$ is used to define 
the Fischer VOA $V\!\!F^\natural_{24}$ acted by the largest Fischer group $\Fi_{24}$ 
as the commutant of $W_3(\sfr{4}{5})$ in the moonshine VOA.
The $W_3$-algebra $W_3(\sfr{4}{5})$ is a subalgebra of the 3A-algebra so that 
the commutant subalgebra $V\!\!F_{23}^\natural= \com_{V^\natural}\la a,b\ra$ affording 
the action of $\Fi_{23}$ is a subVOA of the Fischer VOA $V\!\!F^\natural_{24}$ 
considered in \cite{HLY2}.
It is known that $\Fi_{24}$ is not a subgroup of the Monster but its extension $3.\Fi_{24}$ 
is a subgroup.  
The group $3.\Fi_{24}$ is not a 3-transposition group so that we have to use 
$\sigma$-involutions, which are well-defined only on the subalgebra $V\!\!F^\natural_{24}$ 
of $V^\natural$, to realize the 3-transposition group $\Fi_{24}$ in \cite{HLY2}.
On the other hand, the 3-transposition group obtained by Corollary \ref{cor:3.7} is indeed 
a subgroup of the automorphism group of the whole VOA. 
Hence, our construction naturally reflects the inclusion $\Fi_{23}< \M=\aut(V^\natural)$. 

In order to study inductive structures and basic sets of the 3-transposition group 
obtained by Corollary \ref{cor:3.7}, we consider the following subalgebras.
Let $a$, $b\in E_V$ be a pair such that $\la a,b\ra$ is isomorphic to the 3A-algebra and 
define $I_{a,b}$ as before.
We define $X^{[n]}=\la a,b,x^1,\dots,x^n\ra$ inductively as follows.
We set $X^{[0]}=\la a,b\ra$.
Suppose we have chosen $x^1,\dots,x^i\in I_{a,b}$ and $X^{[i]}=\la a,b,x^1,\dots,x^i\ra$ 
has been defined.
We choose $x^{i+1}\in I_{a,b}$ such that $x^{i+1}\not\in X^{[i]}$ and 
$(x^{i+1}\mymid x^j)=2^{-5}$ for all $1\leq j\leq i$ and define $X^{i+1}=\la X^{[i]},x^{i+1}\ra$
as long as possible.
Suppose we have obtained $X^{[n]}=\la a,b,x^1,\dots,x^n\ra$ in $V_\R$.
Set $D^{[0]}=\{ \tau_y \mid y\in I_{a,b}\}$ and define 
$D^{[i]}=\la \tau_y \in D^{[i-1]} \mid \tau_y \tau_{x^i}=\tau_{x^i}\tau_y\ra$ for $1\leq i\leq n$.
Then $D^{[i]}$ preserves the commutant $\com_V X^{[i]}$ and we obtain the following series 
of commutant subalgebras acted on by inductive substructures 
of the 3-transposition group $G^{[0]}=\la D^{[0]}\ra$.
\begin{equation}\label{eq:1.1}
\begin{array}{ccccccccc}
  ~~~X^{[0]} &\subset& ~~~X^{[1]} &\subset& ~~~X^{[2]} &\subset
  & \cds & \subset & ~~~X^{[n]}
  \vspace{-10pt}\\
  \mbox{\large \rotatebox{270}{$\leadsto$}} 
  && \mbox{\large \rotatebox{270}{$\leadsto$}} 
  && \mbox{\large \rotatebox{270}{$\leadsto$}} 
  && 
  && \mbox{\large \rotatebox{270}{$\leadsto$}} 
  \\
  \com_V X^{[0]} &\supset& \com_V X^{[1]} &\supset& \com_V X^{[2]} &\supset 
  & \cds & \supset & \com_V X^{[n]}
  \\ 
  \circlearrowleft && 
  \circlearrowleft && 
  \circlearrowleft && 
  &&   \circlearrowleft 
  \smallskip\\
  G^{[0]} &&  G^{[1]} && G^{[2]} && \cds && G^{[n]}
  \\
  \mbox{\large \rotatebox{90}{$\twoheadrightarrow$}} 
  && \mbox{\large \rotatebox{90}{$\twoheadrightarrow$}} 
  && \mbox{\large \rotatebox{90}{$\twoheadrightarrow$}}
  && && \mbox{\large \rotatebox{90}{$\twoheadrightarrow$}}
  \\ 
  \la D^{[0]}\ra & \supset & \la D^{[1]}\ra & \supset & \la D^{[2]}\ra & 
  \supset & \cds & \supset & \la D^{[n]}\ra
\end{array}
\end{equation}
where $G^{[i]}$ is the image of $\la D^{[i]}\ra$ in $\aut(\com_V X^{[i]})$.
In the process above, $\{ \tau_{x^1},\dots,\tau_{x^n}\}$ gives a set of mutually commutative 
involutions in $D^{[0]}$ and $X^{[n]}$ is related to the basic sets of a 3-transposition group 
$G^{[0]}$.
We will prove in Theorem \ref{thm:4.17} that the Griess algebra generated by 
$a$, $b$, $x^1,\dots,x^n$ is uniquely determined and $X^{[n]}$ has a full subVOA isomorphic to 
a tensor product of Virasoro VOAs
\begin{equation}\label{eq:1.2}
  L(c_3,0)\tensor L(c_4,0)\tensor \cds \tensor L(c_{n+4},0).
\end{equation}

Our construction of $X^{[n]}$ is compatible with inductive structures of 
3-transposition groups in \eqref{eq:1.1}; thereby applying to the moonshine VOA, 
we will obtain a series of subVOAs affording the actions of the Fischer 3-transposition groups.
As for the moonshine VOA, we obtain $G^{[0]}\cong \Fi_{23}$, $G^{[1]}\cong \Fi_{22}$ and 
$G^{[2]}\cong \Fi_{21}=\mathrm{PSU}_6(2)$ (cf.~Proposition \ref{prop:5.4}).
We know that there exists a bijection between $I_{a,b}\subset V^\natural$ and $D^{[0]}\subset \M$ 
by the Conway-Miyamoto correspondence of Ising vectors and 2A elements of the Monster. 
However, if we regard $G^{[i]}$ as a subgroup of $\aut(\com_{V^\natural} X^{[i]})$, 
the Ising vector $x\in I_{a,b}$ corresponding to $\tau_x\in D^{[i]}$ is not contained in 
the commutant $\com_{V^\natural} X^{[i]}$.
Therefore, we need to replace Ising vectors in $V^\natural$ by suitable Virasoro vectors 
in $\com_{V^\natural} X^{[i]}$ for $i=0,1$, which requires some detailed analysis of the subVOAs 
$X^{[i]}$ and their corresponding Griess algebras.
To establish the Conway-Miyamoto correspondences for $G^{[0]}\cong \Fi_{23}$ and 
$G^{[1]}\cong \Fi_{22}$, we also need to determine the invariants of the centralizers 
$C_{\Fi_{23}}(\mathrm{2A})$ and $C_{\Fi_{22}}(\mathrm{2A})$ based on explicit calculations on 
the Griess algebra of $X^{[2]}$. 
Showing that they are given by $c=c_5$ and $c=c_6$ Virasoro vectors, we will establish  
the Conway-Miyamoto correspondences for $\Fi_{23}$ and $\Fi_{22}$ in Theorems \ref{thm:5.8} 
and \ref{thm:5.15}. 
In principle, we can continue our argument to the next case $\Fi_{21} \cong \mathrm{PSU}_6(2)$ 
where 3-transpositions of $\Fi_{21}$ are described by $c=c_7$ Virasoro vectors.
However, by a technical reason, we cannot determine the Griess algebra and 
its invariant subalgebra of the centralizer (cf.~Remark \ref{rem:5.16}).
The Griess algebra corresponding to $\Fi_{21}$ is rather small and the Conway-Miyamoto 
correspondences seem to terminate at $\Fi_{22}$ in this series.

The study of $X^{[n]}$ has its own interest since it is related to basic sets of $G^{[0]}$ 
in \eqref{eq:1.1}.
Since the size of a basic set of $\Fi_{23}$ is 23, there exist 23 Ising vectors 
$x^1,\dots,x^{23}\in I_{a,b}\subset V^\natural$ satisfying our conditions and  
there is a subalgebra $X^{[23]}=\la a,b,x^1,\dots,x^{23}\ra$ in $V^\natural$.
Surprisingly, the central charge of the Virasoro frame in \eqref{eq:1.2} is equal to 24 
when $n=23$ and hence $X^{[23]}$ is a full subVOA of $V^\natural$.
Therefore, it is possible to analyze $V^\natural$ as a finite extension of $X^{[23]}$. 
We will discuss a construction of  $V^\natural$ starting from $X^{[23]}$ in \cite{CLY}.

The organization of this article is as follows. 
In Section 2, we recall some basic properties about Virasoro VOAs \cite{FZ,W} and 
the dihedral algebras constructed in \cite{LYY1,LYY2}. 
In Section 3, we give a construction of 3-transposition groups using the 3A-algebra.
We first fix a pair of Ising vectors $a, b\in V$ such that $a$ and $b$ generate 
a 3A-algebra in $V$ and then consider the set of Ising vectors $I_{a,b}$ of $V$ such that 
$\la a,x\ra$ and $\la b,x\ra$ are isomorphic to the 2A-algebra for all $x\in I_{a,b}$. 
We will show that the Miyamoto involutions associated to Ising vectors of $I_{a,b}$ generate a 
$3$-transposition group in $\aut(V)$. 
In Section 4, we will determine all possible structures of the subVOA generated by $a$, $b$, 
$x$, $y$ for any $x$, $y\in I_{a,b}$ and then introduce inductive subalgebras 
$X^{[i]}=\la a,b,x^1,\dots,x^i\ra$ for $i\geq 0$. 
The commutants of subalgebras $X^{[i]}$ are used to study inductive structures  
of the 3-transposition group obtained in Section 3.
The results in Section 4.2 will not be used in the later section but we will need 
these results in our future work so that we include them.
In Section 5, we will apply our results to the moonshine VOA and the Fischer 
3-transposition groups $\Fi_{23}$ and $\Fi_{22}$.
The Conway-Miyamoto correspondence between the 2A-involutions of $\Fi_{23}$ and $c=c_5$ 
Virasoro vectors of $\VF_{23}=\com_{V^\natural} X^{[0]}$ and that between the 2A-involutions 
of $\Fi_{22}$ and $c=c_6$ Virasoro vectors of $\sigma$-type of $\VF_{22}=\com_{V^\natural} X^{[1]}$ 
will be established. 
In the appendix, we will give explicit constructions of the VOAs studied in Section 4.

The authors used computer algebra systems Risa/Asir for the calculations in Griess algebras 
and GAP 4.7.4 for Linux for the character calculations of finite groups.

\paragraph{Acknowledgement.}
Part of this work has been done while the author(s) were staying 
at Academia Sinica in September 2011, March 2013 and March 2015, 
at National Taitung University 
in March 2013 and 
March 2015,  
at Imperial College London during ``Majorana Theory, the Monster and Beyond'' 
in September 2013,  
at National Dong Hwa University 
in March 2014, 
at Mathematisches Forschungsinstitut Oberwolfach during ``Subfactor and Conformal Field Theory'' 
in March 2015, 
and at Sichuan University 
in September 2015. 
They gratefully acknowledge the hospitality there. 
H.Y. thanks Masaaki Kitazume, Atsushi Matsuo, Masahiko Miyamoto, Hiroki Shimakura 
and Katsushi Waki for valuable comments and helpful discussions.

\paragraph{Notation and terminology.}
In this paper, vertex operator algebras (VOAs) are defined over the complex number field $\C$.
A VOA $V$ is called {\it of OZ-type} if it has the $L(0)$-grading 
$V=\oplus_{n\geq 0}V_n$ such that $V_0=\C \vac$ and $V_1=0$.
In this case, $V$ is equipped with a unique invariant bilinear form such that $(\vac|\vac)=1$.
A real form $V_\R$ of $V$ is called \emph{compact} if the associated bilinear form 
is positive definite.
We will mainly consider a VOA of OZ-type having a compact real form.
For a subset $A$ of $V$, the vertex subalgebra generated by $A$ is denoted by $\la A\ra$. 
(In most cases, $\la A\ra$ indeed becomes a subVOA in this paper.)
For $a\in V_n$ we define $\wt(a)=n$.
We write $Y(a,z)=\sum_{n\in \Z}a_{(n)}z^{-n-1}$ for $a\in V$ and define 
its {\it  zero-mode} by $\o(a):=a_{(\wt(a)-1)}$ if $a$ is homogeneous and extend linearly.
The weight two subspace $V_2$ carries a structure of a commutative algebra defined by 
the product $\o(a)b=a_{(1)}b$ for $a$, $b\in V_2$.
This algebra is called the \emph{Griess algebra} of $V$.
A \emph{Virasoro vector} is $a\in V_2$ such that $a_{(1)}a=2a$.
The subalgebra $\la a\ra$ is isomorphic to a Virasoro VOA with 
central charge $2(a|a)$.
If two Virasoro vectors $a$ and $b$ are orthogonal, we will denote their sum by $a \dotplus b$.
We denote by $L(c,h)$ the irreducible highest weight module over the Virasoro algebra 
with central charge $c$ and highest weight $h$.
A \emph{simple} $c=c_e$ Virasoro vector $e\in V$ is is a Virasoro vector such that 
$\la e\ra \cong L(c_e,0)$.
If $e$ is taken from a compact real form of $V$ then $e$ is always simple.
A Virasoro vector $\w$ is called the \emph{conformal vector} of $V$ if each graded 
subspace $V_n$ agrees with $\ker_V\, (\o(\w)-n)$ and satisfies 
$\w_{(0)}a =a_{(-2)}\vac$ for all $a\in V$.
The half of the conformal vector gives the unit of the Griess algebra and 
hence uniquely determined.
We write $\w_{(n+1)}=L(n)$ for $n\in \Z$.
A Virasoro vector $e$ of $V$ is called \emph{characteristic} if it is fixed by 
$\aut(V)$.
Clearly the conformal vector is characteristic in $V$.
A \emph{subVOA} $(W,e)$ of $V$ is a pair of a subalgebra $W$ of $V$ 
together with a Virasoro vector $e$ in $W$ such that $e$ is the conformal vector of $W$.
We often omit $e$ and simply denote by $W$.
The commutant subalgebra of $(W,e)$ in $V$ is defined by $\com_V W:=\ker_V\,e_{(0)}$ 
(cf.~\cite{FZ}).
A subVOA $W$ of $V$ is called {\it full} if $V$ and $W$ shares the same conformal vector.
For a subgroup $G$ of $\aut(V)$, we denote the set of $G$-invariants by $V^G$.

\section{The dihedral subalgebras}
First we recall some basic properties about Virasoro VOAs \cite{FZ,W} and 
the dihedral algebras constructed in \cite{LYY1,LYY2}.

\subsection{Virasoro vertex operator algebras}\label{sec:2.1}

Let
\begin{equation}\label{eq:2.1}
\begin{array}{rl}
  c_n &:= 1-\dfr{6}{(n+2)(n+3)},\qq n=1,2,3,\dots ,
  \medskip\\
  h_{r,s}^{(n)} &:= \dfr{\{ r(n+3)-s(n+2)\}^2-1}{4(n+2)(n+3)},~~
  1\leq r\leq n+1,~~1\leq s\leq n+2.
\end{array}
\end{equation}
It is shown in \cite{W} that $L(c_n,0)$ is rational and $L(c_n,h_{r,s}^{(n)})$,
$1\leq s\leq r\leq n+1$, are all the irreducible $L(c_n,0)$-modules 
(see also \cite{DMZ}). 
This is the so-called unitary series of the Virasoro VOAs. 
Note that $h_{r,s}^{(n)}=h_{n+2-r,n+3-s}^{(n)}$.
The fusion rules among $L(c_n,0)$-modules are computed in \cite{W} and given by
\begin{equation}\label{eq:2.2}
  L\l( c_n,h^{(n)}_{r_1,s_1}\r)\fusion L\l( c_n,h^{(n)}_{r_2,s_2}\r)
  = \dsum_{1\leq i\leq M\atop 1\leq j\leq N}
    L(c_n,h^{(n)}_{\abs{r_1-r_2}+2i-1,\abs{s_1-s_2}+2j-1}),
\end{equation}
where $M=\min \{ r_1,\,r_2,\,n+2-r_1,\,n+2-r_2\}$ and $N=\min \{ s_1,\,s_2,\,n+3-s_1,\,n+3-s_2\}$.

\begin{df}\label{df:2.1}
  A Virasoro vector $e$ with central charge $c$ is called {\it simple\/}
  if $\la e\ra\cong L(c,0)$.
  A simple $c=1/2$ Virasoro vector is called an {\it Ising vector}.
\end{df}

\medskip

The fusion rules among $L(c_m,0)$-modules have a canonical
$\Z_2$-symmetry. 

\begin{thm}[\cite{M1}]\label{thm:2.2}
  Let $V$ be a VOA and $e\in V$ a simple Virasoro vector
  with central charge $c_n$.
  Denote by $V_e[h^{(n)}_{r,s}]$ the sum of irreducible
  $\la e\ra\cong L(c_n,0)$-submodules isomorphic to $L(c_n,h^{(n)}_{r,s})$,
  $1\leq s\leq r\leq n+1$.
  Then the linear map
  \[
    \tau_e
    := (-1)^{4(n+2)(n+3)\o(e)}
    =
    \begin{cases}
      (-1)^{r+1}
      & \text{on }\ V_e[h^{(n)}_{r,s}]\q \text{if } n \ \text{is even},
      \medskip\\
      (-1)^{s+1}
      & \text{on }\ V_e[h^{(n)}_{r,s}]\q \text{if } n \ \text{is odd},
    \end{cases}
  \]
  defines an automorphism of $V$. 
\end{thm}

The automorphism $\tau_e$ is called the \emph{Miyamoto involution} associated to $e$.

\begin{df}\label{df:2.3}
  Let $e$ be a simple $c=c_n$ Virasoro vector of $V$.
  Set
  \[
    P_n:=\begin{cases}
    \{ h_{1,s}^{(n)} \mid 1\leq s\leq n+2\} & \mbox{if $n$ is even,}
    \medskip\\
    \{ h_{r,1}^{(n)} \mid 1\leq r\leq n+1\} & \mbox{if $n$ is odd.}
    \end{cases}
  \]
  It follows from the fusion rules in \eqref{eq:2.3} that the subspace 
  \[
    V_e[P_n]=\bigoplus_{h\in P_n} V_e[h]
  \]
  forms a subalgebra of $V$.
  We say that $e$ is \emph{of $\sigma$-type on $V$} if $V=V_e[P_n]$.
\end{df}

\begin{thm}[\cite{M1}]\label{thm:2.4}
  The linear map 
  \[
    \sigma_e :=\begin{cases}
      (-1)^{s+1} ~~\mbox{on~ $V_e[h_{1,s}^{(n)}]$~ if~ $n$~ is even,}
      \medskip\\
      (-1)^{r+1} ~~\mbox{on~ $V_e[h_{r,1}^{(n)}]$~ if~ $n$~ is odd,} 
    \end{cases}
  \]
  defines an element of $\aut(V_e[P_n])$. 
\end{thm}

The automorphism $\sigma_e$ is called a Miyamoto involution of $\sigma$-type or  simply a $\sigma$-involution.
Suppose $e$ is an Ising vector, i.e., simple $c=1/2$ Virasoro vector. 
Then we have 
\begin{equation}\label{eq:2.3}
  \tau_e =\begin{cases}
    ~1 & \text{on }~ V_{e}[0]\oplus V_{e}[\sfr{1}{2}],\\
    -1 & \text{on }~ V_{e}[\sfr{1}{16}].
  \end{cases}
\end{equation}
If $e$ is of $\sigma$-type on $V$, the linear map $\sigma_{e}$ is defined by
\begin{equation}\label{eq:2.4}
\sigma_{e}:=\begin{cases}
~1 & \text{on }~ V_{e}[0],\\
-1 & \text{on }~ V_{e}[\sfr{1}{2}].
\end{cases}
\end{equation}
Alternatively, one can define $\sigma_e=(-1)^{2\o(e)}$.
By the definition of Miyamoto involutions,  we have the following conjugation.

\begin{prop}\label{prop:2.5}
  Let $e$ be a simple $c=c_n$ Virasoro vector of $V$.
 Then $\tau_{ge}=g\tau_e g^{-1}$ for any $g\in \aut(V)$. 
 If $e$ is of $\sigma$-type on $V$, we have  $\sigma_{ge}=g\sigma_e g^{-1}$.
\end{prop}

For an Ising vector $e$ and $v\in V_2$, one has $v+\tau_e v \in V^{\la \tau_e\ra}$ 
and hence $\sigma_e(v+\tau_e v)$ is a well-defined element in $V_2$.
Moreover, we have (cf.~Eq.~(2.2) of \cite{S}):
\begin{equation}\label{eq:2.5}
  e_{(1)}v=8(e \mymid v)e+\dfr{5}{32}v+\dfr{3}{32}\tau_e v-\dfr{1}{8}\sigma_e(v+\tau_e v).
\end{equation}
In particular, if $\tau_e v= v$, we have 
\begin{equation}\label{eq:2.6}
  \sigma_e v=v+32(e \mymid v)e-4e_{(1)}v.
\end{equation}
These relations will be used in Section \ref{sec:3}.

\medskip

In \cite{LYY1,LYY2},  subalgebras generated by two Ising vectors are constructed using 
the $E_8$-lattice and such subalgebras are classified in \cite{S} at the level of 
Griess algebras.

\begin{thm}[\cite{S}]\label{thm:2.6}
  Let $V$ be a VOA of OZ-type with compact real form $V_\R$ and let $e$ and $f$ be 
  distinct Ising vectors of $V_\R$.
  Then the Griess algebra of the subalgebra $\la e,f\ra$ is isomorphic to one of the eight 
  algebras called 2A, 3A, 4A, 5A, 6A, 4B, 2B and 3C-type constructed in \cite{LYY1,LYY2}.
  The inner product and the dimension of the Griess algebra are as follows\footnote{%
  Here we include the case $e=f$ which is called the 1A-type.}.
  \[
  \renewcommand{\arraystretch}{1.2}
  \begin{array}{c|ccccccccc}
    \mbox{\rm{Type}} & \mbox{\rm{1A}} & \mbox{\rm{2A}} & \mbox{\rm{3A}} 
    & \mbox{\rm{4A}} &\mbox{\rm{5A}} & \mbox{\rm{6A}} & \mbox{\rm{4B}} 
    & \mbox{\rm{2B}} & \mbox{\rm{3C}} 
    \\ \hline
    2^{10} (e\mymid f) & 2^8 & 2^5 & 13 & 8 & 6 & 5 & 4 & 0 & 4
    \\ \hline
    \dim \la e,f\ra_2 & 1 & 3 & 4 & 5 & 6 & 8 & 5 & 2 & 3
    \\ \hline
    \# ~\mbox{{\rm of Ising vectors}} & 1& 3& 3& 4& 5& 7& 5 & 2 & 3
  \end{array}
  \renewcommand{\arraystretch}{1}
  \]
\end{thm}

We call $\la e,f\ra$ the \emph{dihedral subalgebra} and we denote the dihedral 
subalgebra of type $n$X by $U_{n\mathrm{X}}$.
As a by-product of the classification of dihedral subalgebras, 
the following 6-transposition property was established in \cite{S}.

\begin{thm}[\cite{S}]\label{thm:2.7}
  Let $V$ be a VOA of OZ-type with compact real form $V_\R$ and let $e$ and $f$ be 
  distinct Ising vectors of $V_\R$.
  Then
  \\
  (1)~ If $\la e,f\ra$ is 2A or 2B-type then $\abs{\tau_e\tau_f}$ divides 2.
  \\
  (2)~ If $\la e,f\ra$ is 3A or 3C-type then $\abs{\tau_e\tau_f}=3$.
  \\
  (3)~ If $\la e,f\ra$ is 4A or 4B-type then $\abs{\tau_e\tau_f}=2$ or $4$.
  \\
  (4)~ If $\la e,f\ra$ is 5A-type then $\abs{\tau_e\tau_f}=5$.
  \\
  (5)~ If $\la e,f\ra$ is 6A-type then $\abs{\tau_e\tau_f}=3$ or $6$.
  \\
  (6)~ $\tau_e f=f$ if and only if $\la e,f\ra$ is either 2A or 2B-type.
  \\
  (7)~ $\tau_e f=\tau_f e$ if and only if $\la e,f\ra$ is either 3A or 3C-type.
  \\
  In particular, the order of $\tau_e\tau_f$ is bounded by 6.
\end{thm}

The next lemma also follows from the classification of the dihedral algebras.

\begin{lem}\label{lem:2.8}
  Let $a^1$, $a^2$, $b^1$, $b^2$ be Ising vectors and set $U^1=\la a^1,b^1\ra$ and 
  $U^2=\la a^2,b^2\ra$.
  If $U^1$ is a proper subalgebra of $U^2$ then the types of $(U^1,U^2)$ are either 
  $\mathrm{(2A,4B)}$, $\mathrm{(2A,6A)}$, $\mathrm{(2B,4A)}$ or $\mathrm{(3A,6A)}$.
\end{lem}

In the rest of this section, we will summarize some properties of the dihedral 
subalgebras of 2A, 3A and 6A-types which will be used in the later sections.
See \cite{LYY2} for details.

\subsection{2A-algebra}\label{sec:2.2}

It follows from Theorem \ref{thm:2.6} that $\la a,b\ra$ is of 2A-type if and only if 
$(a\mymid b)=2^{-5}$.
Let $U_{\mathrm{2A}}=\la a,b\ra$ be the 2A-algebra. 
Then $a$ and $b$ are of $\sigma$-type on $\la a,b\ra$ and satisfies $\sigma_a b=\sigma_b a$.
Set $a\circ b:=\sigma_a b=\sigma_b a$.
It follows from \eqref{eq:2.6} that
\begin{equation}\label{eq:2.7}
  a\circ b = a+b-4a_{(1)}b.
\end{equation}

\begin{thm}\label{thm:2.9}
  Let $U_{\mathrm{2A}}=\la a,b\ra$ be the 2A-algebra.
  \\
  (1)~There are exactly three Ising vectors in of $U_{\mathrm{2A}}$, namely, $a$, $b$ and $a\circ b$.
  \\
  (2)~The Griess algebra of $U_{\mathrm{2A}}$ is 3-dimensional spanned by these three Ising vectors. 
  \\
  (3)~$(a \mymid b)=(a \mymid a\circ b)=(b \mymid a\circ b)=2^{-5}$.
  \\
  (4)~$\aut(U_{\mathrm{2A}})=\la \sigma_a,\sigma_b\ra \cong \mathrm{S}_3$.
  \\
  (5)~Suppose $U_{\mathrm{2A}}=\la a,b\ra$ is a subalgebra of a VOA $V$.
  Then $\tau_{a\circ b}=\tau_a\tau_b=\tau_b\tau_a$ on $V$.
\end{thm}

We will call the set of Ising vectors $\{ a,b,a\circ b\}$ of $\la a,b\ra$ a  
\emph{2A-triple}.
The subalgebra generated by three Ising vectors of $\sigma$-type was classified 
in \cite{Ma2} and the following ``No 2A-tetrahedron lemma'' holds 
(cf.~Proposition 1 of \cite{Ma2}\footnote{Proposition 3.3.8 of \texttt{arXiv:math/0311400}.}).

\begin{lem}\label{lem:2.10}
Let $a$, $b$, $c$ be Ising vectors such that $(a \mymid b)=(a\mymid c)=(b\mymid c)=2^{-5}$ and $c\notin \langle a, b\rangle$.
Then we have $(a\circ b \mymid c)=(a \mymid  b \circ c)=(b \mymid a\circ c)=0$.
\end{lem}

\subsection{3A-algebra}\label{sec:2.3}

It follows from Theorem \ref{thm:2.6} that $\la a,b\ra$ is of 3A-type if and only if 
$(a \mymid b)=13\cd 2^{-10}$.
Let $U_{\mathrm{3A}}=\la a,b\ra$ be the 3A-algebra.
Then $\tau_a b= \tau_b a$ holds.
Set $c=\tau_a b$ and define 
\begin{equation}\label{eq:2.8}
  u=u_{a,b}=\dfr{2^6}{135}\l( 2a+2b+c-16a_{(1)}b\r) .
\end{equation}
Then $u$ is a $c=4/5$ Virasoro vector in $U_{\mathrm{3A}}$.

\begin{thm}\label{thm:2.11}
  Let $U_{\mathrm{3A}}=\la a,b\ra$ be the 3A-algebra.
  \\
  (1)~There are exactly three Ising vectors in $U_{\mathrm{3A}}$, namely,  $a,b$ and $c$.
  \\
  (2)~$\aut(U_{\mathrm{3A}})=\la \tau_a,\tau_b\ra\cong \mathrm{S}_3$.
  \\
  (3)~The Griess algebra of $U_{\mathrm{3A}}$ is 4-dimensional and  spanned by $u$, $a$, $b$ and $c$.
  \\
  (4)~The $c=4/5$ Virasoro vector $u$ is characteristic in $U_{\mathrm{3A}}$.
  \\
  (5)~The multiplications and inner products of the Griess algebra are as follows.
  \[
  \begin{array}{l}
    a_{(1)}b=-\dfr{135}{2^{10}}u+\dfr{1}{2^3}a+\dfr{1}{2^3}b+\dfr{1}{2^4}c, ~~~
    a_{(1)}u=\dfr{5}{2^4}u+\dfr{4}{3^2}a-\dfr{2}{3^2}b-\dfr{2}{3^2}c,
    \medskip\\
    (a \mymid b)=(a \mymid c)=(b \mymid c)=\dfr{13}{2^{10}},~~~
    (a \mymid u)=(b \mymid u)=(c \mymid u)=\dfr{1}{2^4}.
  \end{array}
  \]
  (6)~There is a unique characteristic $c=6/7$ Virasoro vector in $\la a,b\ra$, which is given by 
  \[
    v=-\dfr{5}{14}u+\dfr{16}{21}(a+b+c).
  \]
  Moreover, $u$ and $v$ are mutually orthogonal and $u+v$ is the conformal vector of 
  $\la a,b\ra$.
\end{thm}

We call the set of Ising vectors in $U_{\mathrm{3A}}$ a \emph{3A-triple}.
We will need the following relations later.

\begin{lem}\label{lem:2.12}
  In the 3A-algebra,  the following holds.
  \[
  \begin{array}{l}
    \sigma_{a}u=-\dfr{1}{2^2}u+\dfr{2}{3^2}a+\dfr{8}{3^2}(b+c),~~~
    \sigma_{a}(b+c)=\dfr{135}{2^{7}}u-\dfr{3}{2^{4}}a+\dfr{1}{2^2}(b+c).
  \end{array}
  \]
\end{lem}

\subsection{6A-algebra}\label{sec:2.4}

It follows from Theorem \ref{thm:2.6} that $\la a,b\ra$ is of 6A-type if and only if 
$(a \mymid b)=5\cd 2^{-10}$.
Let $U_{\mathrm{6A}}=\la a,b\ra$ be the 6A-algebra.
There are 7 Ising vectors in the 6A-algebra.
Namely, $a$, $b$, $\tau_a b$, $\tau_b a$, $\tau_a\tau_b a$, $\tau_b\tau_a b$ and 
\begin{equation}\label{eq:2.9}
  x=a\circ \tau_b\tau_a b= b\circ \tau_a\tau_b a=\tau_b a \circ \tau_a b.
\end{equation}
The Ising vector $x$ is of $\sigma$-type  and defines a $\sigma$-involution
on $\la a,b\ra$.
Inside the 6A-algebra, there are two 3A-subalgebras and the 3-sets $\{ a,\tau_b a,\tau_a\tau_b a\}$
and $\{ b, \tau_a b, \tau_b\tau_a b\}$ form the 3A-triples with the same characteristic 
$c=4/5$ Virasoro vector $u_{a,\tau_b a}=u_{b,\tau_a b}$.
The automorphism $\tau_a\tau_b$ acts as a 3-cycle on each of the 3A-triples and 
has order 3 on $\la a,b\ra$.

\begin{thm}\label{thm:2.13}
  Let $U_{\mathrm{6A}}=\la a,b\ra$ be the 6A-algebra.
  \\
  (1)~The set of Ising vectors of $U_{\mathrm{6A}}$ is given by 
  $\{ x,a,b,\tau_a b,\tau_b a,\tau_a\tau_b a, \tau_b \tau_a b\}$.
  \\
  (2)~$x$ is the unique characteristic $\sigma$-type Ising vector in $U_{\mathrm{6A}}$ 
  and $\sigma_x$ is central in $\aut(U_{\mathrm{6A}})$.
  \\
  (3)~$\{ x,a,\tau_b\tau_ab\}$, $\{ x,b,\tau_a\tau_ba\}$ and $\{ x,\tau_b a,\tau_ab\}$ 
  form 2A-triples in $\la a,b\ra$.
  \\
  (4)~$\la a,\tau_b a,\tau_a\tau_b a\ra$ and $\la b,\tau_ab,\tau_b\tau_a b\ra$ 
  are isomorphic to the 3A-algebra with the same characteristic $c=4/5$ Virasoro vector 
  $u=u_{a,\tau_b a}=u_{b,\tau_a b}$.
  \\
  (5)~The Griess algebra of $U_{\mathrm{6A}}$ is 8-dimensional with a basis 
  $\{ u,x,a,b,\tau_a b,\tau_b a,\tau_a\tau_b a,\tau_b \tau_a b\}$.
  \\
  (6)~$\aut(U_{\mathrm{6A}})=\la \tau_a\sigma_x,\tau_b\ra$ is isomorphic to the dihedral group 
  of order 12.
  \\
  (7)~Set $\rho=\tau_a\tau_b\sigma_x$ and $e^i=\rho^i a$.
  Then $\{ a,b,\tau_a b, \tau_b a, \tau_a \tau_b a, \tau_b \tau_a b\}=\{ e^i \mid 0\leq i\leq 5\}$ 
  and the multiplications and inner products of the Griess algebra are as follows.
  \[
  \begin{array}{l}
    u_{(1)}x=0,~~~ 
    u_{(1)}e^i=-\dfr{135}{2^{10}}u+\dfr{4}{3^2}e^i-\dfr{2}{3^2}(e^{i+2}+e^{i+4}),~~~
    (u \mymid x)=0,~~~ 
    (u \mymid e^i)=\dfr{1}{2^4},~~~
    \medskip\\
    x_{(1)}e^i=\dfr{1}{4}x+e^i-e^{i+3},~~~
    e^0_{(1)}e^1= \dfr{45}{2^{10}}u+\dfr{1}{2^5}\l( x+e^0+e^1-e^2-e^3-e^4-e^5\r),~~~
    \medskip\\
    (x \mymid e^i)=\dfr{1}{2^5},~~~
    (e^i \mymid e^{i+1})=\dfr{5}{2^{10}},~~~
    (e^i \mymid e^{i+2})=\dfr{13}{2^{10}},~~~
    (e^i \mymid e^{i+3})=\dfr{1}{2^5}.
  \end{array}
  \]
  (8)~The commutant of $\la a,\tau_b a\ra$ in $\la a,b\ra$ is generated by a 
  $c=25/28$ Virasoro vector
  \[
    f=-\dfr{15}{56} u + \dfr{1}{2}x -\dfr{2}{21}(a+\tau_b a +\tau_a\tau_b a)
      +\dfr{2}{3}(b+ \tau_a b + \tau_b \tau_a b).     
  \]
  Therefore, $\la a,b\ra$ has a unitary Virasoro frame isomorphic to 
  $L(\sfr{4}{5},0)\tensor L(\sfr{6}{7},0)\tensor L(\sfr{25}{28},0)$.
\end{thm}

The characteristic Ising vector $x$ is called \emph{central} since 
$\sigma_x$ is central in the automorphism group of the 6A-algebra.
We collect some technical relations needed in the later arguments.

\begin{lem}\label{lem:2.14}
  Let $\la a,b\ra$ be the 6A-algebra and let $u$ be the characteristic
  $c=4/5$ Virasoro vector and $x$ the central Ising vector.
  Then one has
  \[
  \begin{array}{l}
    a_{(1)}(\tau_b a)_{(1)}x=-\dfr{45}{2^{10}}u+\dfr{1}{2^7}
    \l( 7x+11a+b+5\tau_b a-7\tau_b\tau_a b +3\tau_a\tau_b a-\tau_a b\r) ,
    \medskip\\
    \sigma_{a} (b+\tau_a b ) = -\dfr{45}{2^{7}}u_{a,b} -\dfr{1}{2^{2}}x+\dfr{1}{2^{4}}a
    +\dfr{1}{2^{2}}\tau_b \tau_a b +b+\tau_a b +\dfr{1}{2^{2}}(\tau_b a +\tau_a \tau_b a).
  \end{array}
  \]
\end{lem}

\begin{rem}\label{rem:2.15}
  In some situation, it is  useful  to display the relations above using another labeling.
  We rename the Ising vectors so that $a$ and $b$ generate a $3A$-algebra and $x$ is 
  the central Ising vector in a 6A-algebra. Let $c=\tau_a b=\tau_b a$. Then $\{a, b,c\}$ forms a 3A-triple and 
  $\{a\circ x, b\circ x, c \circ x\}$ gives another 3A-triple. 
  Lemma \ref{lem:2.14} can be rewritten as follows.
  \begin{equation}\label{eq:2.10}
  \begin{array}{l}
    a_{(1)}b_{(1)}x=-\dfr{45}{2^{10}} u_{a,b} +\dfr{1}{2^7}\l( 7x+11a+5b+3c-7a\circ x
    -b\circ x+c\circ x\r) , 
    \medskip\\
    \sigma_{a} (b\circ x+c\circ x ) = -\dfr{45}{2^{7}}u -\dfr{1}{2^{2}}x+\dfr{1}{2^{4}}a
    +\dfr{1}{2^{2}}a\circ x +\dfr{1}{2^{2}}(b  +c)+b\circ x+c\circ x.
  \end{array}
  \end{equation}
\end{rem}

\section{3-transposition property}\label{sec:3}

Suppose that $V$ has a grading $V=\bigoplus_{n\geq 0}V_n$
with $V_0=\C\vac$ and $V_1=0$.
We also assume that $V$ has a compact real form $V_\R$ and
all Ising vectors are taken from it.

\begin{lem}\label{lem:3.1}
  Let $e$ be an Ising vector and $v\in V_2$.
  Then $\tau_e v\in \la e,v\ra$.
  If $v\in V^{\la \tau_e\ra}$, then $\sigma_e v\in \la e,v\ra$.
\end{lem}

\pf
We have the eigenspace decomposition 
\[
  V_2=\C e \oplus (V_2)_e[0] \oplus (V_2)_e[\shf]\oplus (V_2)_e[\sfr{1}{16}].
\]
Decompose $v=\lambda e +v_0+v_{1/2}+v_{1/16}$ with $\lambda \in \C$ and $v_h\in (V_2)_e[h]$.
Then we have 
\[
  16e_{(1)}v=32\lambda e+ 8 v_{1/2}+v_{1/16} \mbox{~~and~~}
  256e_{(1)}e_{(1)}v=1024 \lambda e+ 64 v_{1/2}+v_{1/16}.
\]
It follows that all of $v_0$, $v_{1/2}$, $v_{1/16}$ are in 
$\Span\{ v,e,e_{(1)}v,e_{(1)}e_{(1)}v\}$.
Therefore, $\tau_e v=v-2v_{1/16}\in \la e,v\ra$.
If $v\in V^{\la \tau_e\ra}$ then it follows from \eqref{eq:2.6} that $\sigma_e v \in \la e,v\ra$. 
\qed

\begin{lem}\label{lem:3.2}
  Let $a,b,x$ be Ising vectors of $V_\R$ such that  $\tau_a b=\tau_b a$ and
  $\la a,x\ra \cong \la b,x\ra \cong U_{\mathrm{2A}}$.
  Then $\la a,b\ra \cong U_{\mathrm{3A}}$, $\la a,b,x\ra \cong U_{\mathrm{6A}}$ and $x$ is
  the central Ising vector of $\la a,b,x\ra$.
\end{lem}

\pf
Since $a,b\in V^{\la \tau_x\ra}$,
$x$ is of $\sigma$-type on $\la a,b,x\ra$ and one has $\tau_{a\circ x}=\tau_a\tau_x$.
We will show that $\la a\circ x,b\ra = \la a,b,x\ra$.
By Lemma \ref{lem:3.1}, $\la a\circ x,b\ra\subset \la a,b,x\ra$.
Since $\tau_x b=b$, we have $\tau_{a\circ x}b=\tau_a\tau_xb=\tau_a b$.
By our assumption, we have $\tau_a b =\tau_b a$.
Then $\tau_{a\circ x}b=\tau_a\tau_xb=\tau_a b=\tau_b a$ and 
$a=\tau_{b}\tau_{a\circ x}b \in \la a\circ x,b\ra$ by Lemma \ref{lem:3.1}.
Since $a\circ (a\circ x)=x$, we also have $x=a+a\circ x-4a_{(1)}(a\circ x)
\in \la a\circ x,b\ra$.
Therefore, $a,b,x\in \la a \circ x,b\ra$ and $\la a,b,x\ra =\la a \circ x,b\ra$.
By (7) of Theorem \ref{thm:2.7}, the assumption $\tau_a b=\tau_b a$ implies 
$\la a,b\ra$ is isomorphic to either $U_{\mathrm{3A}}$ or $U_{\mathrm{3C}}$.
Since $\la a,b\ra \subset \la a \circ x,b\ra$, the possible types of
$\la a,b\ra$ and $\la a \circ x,b\ra$ are only $U_{\mathrm{3A}}$ and $U_{\mathrm{6A}}$, respectively,
by Lemma \ref{lem:2.8}.
As is of $\sigma$-type on $\la a,b,x\ra\cong U_{\mathrm{6A}}$,
$x$ is the unique central Ising vector in $\la a,b,x\ra$.
This completes the proof.
\qed

\begin{lem}\label{lem:3.3}
  Let $a$, $b$, $x$ be Ising vectors of $V_\R$ and such that 
  $\la a,x\ra\cong \la b,x\ra \cong U_{\mathrm{2A}}$ and $\la a,b\ra \cong U_{\mathrm{6A}}$.
  Let $z$ be the central Ising vector of $\la a,b\ra$.
  Then $\tau_x=\tau_z$ on $V_\R$.
\end{lem}

\pf
Suppose $\tau_x\ne \tau_z$ on $V_\R$.
Let $a'=\tau_b\tau_a b$.
Then $a$, $a'$ and $z=a\circ a'$ form a 2A-triple by (3) of Theorem \ref{thm:2.13}.
By Proposition \ref{prop:2.5} and (5) of Theorem \ref{thm:2.9}, we have 
\[
  \tau_z=\tau_a \tau_{a'}=\tau_a \tau_{\tau_b\tau_a b}=(\tau_a\tau_b)^3.
\]
Since $\la a,x\ra$ and $\la b,x\ra$ are of 2A-type, we have $[\tau_x,\tau_a]=[\tau_x,\tau_b]=1$. 
Then $\tau_x$ commutes with $\tau_z=(\tau_a\tau_b)^3$ and $\tau_{a'}=\tau_a\tau_z$.
Therefore, $G=\la \tau_a,\tau_{a'},\tau_x\ra \subset \aut(V_\R)$ is an abelian subgroup.
The order of $\tau_a\tau_b$ is either 3 or 6 on $V_\R$. 
Since 
\[
  (\tau_{a\circ x}\tau_b)^3
  = (\tau_x\tau_a\tau_b)^3=\tau_x^3(\tau_a\tau_b)^3
  = \tau_x  \tau_z\neq 1,
\]
the order of $\tau_{a\circ x}\tau_b$ is 6 by the 6-transposition property and 
we have $\la a\circ x,b\ra \cong U_{\mathrm{6A}}$ by Theorem \ref{thm:2.7}.
It follows from (3) of Theorem \ref{thm:2.13} that $a\circ x$ and $\tau_b\tau_{a\circ x}b$ 
generates a 2A-subalgebra in $\la a\circ x,b\ra$.
On the other hand, by (5) of Theorem \ref{thm:2.9} we have 
$\tau_b\tau_{a\circ x}b=\tau_b\tau_a\tau_x b =\tau_b\tau_a b= a'$.
Therefore $(a\circ x \mymid a')=(a\circ x \mymid \tau_b\tau_{a\circ x}b)=2^{-5}$.
Since $a$ and $b$ are fixed by $\tau_x$, $z=a\circ a'$ is also fixed by $\tau_x$.
Clearly $a$ and $a'$ are fixed by $\la \tau_a,\tau_{a'}\ra$ and hence $a,a',z\in V_\R^G$.
Similarly, $x$ is fixed by $\la \tau_a,\tau_b\ra$ so that $x$ is fixed by $\tau_a$ and 
$\tau_{a'}=\tau_{\tau_b\tau_a b}=\tau_b\tau_a\tau_b\tau_a\tau_b$.
Therefore, all of $a$, $a'$, $z$ and $x$ belong to $V_\R^G$ and of $\sigma$-type on it.
We have $(a'\mymid x)=(\tau_b\tau_a b \mymid x)=(b \mymid \tau_a\tau_bx)=(b \mymid x)=2^{-5}$ 
and $(z \mymid x)=(a\circ a' \mymid x)=(\sigma_{a'}a \mymid x)=(a \mymid \sigma_{a'}x)
=(a \mymid \sigma_x a')= (\sigma_x a \mymid a')=(a\circ x \mymid a')=2^{-5}$.
Therefore, we have 
$(a \mymid a')=(a \mymid z)=(a \mymid x)=(a' \mymid z)=(a' \mymid x)=(z \mymid x)=2^{-5}$.
Namely, $a$, $a'$, $x$ and $z$ form a 2A-tetrahedron in $V_\R^G$ which contradicts 
Lemma \ref{lem:2.10}.
Thus $\tau_x=\tau_z$.
\qed

\begin{lem}\label{lem:3.4}
  Let $a,b,x,y$ be Ising vectors of $V$ such that
  $\la a,x\ra \cong \la b,x\ra \cong \la a,y\ra \cong \la b,y\ra \cong U_{\mathrm{2A}}$
  and $\la x,y\ra \cong U_{\mathrm{6A}}$.
  Then $\tau_a=\tau_b$.
\end{lem}

\pf
Let $z$ be the central Ising vector of $\la x,y\ra$.
Applying Lemma \ref{lem:3.3} to $x,y,a$, we obtain $\tau_a=\tau_z$, and
applying Lemma \ref{lem:3.3} to $x,y,b$, we obtain $\tau_b=\tau_z$.
Thus $\tau_a=\tau_b$.
\qed
\medskip

Let $E_V$ be the set of Ising vectors of $V_\R$.
Suppose there is a pair $a,b\in E_V$ such that $\la a,b\ra\cong U_{\mathrm{3A}}$. 
Then we set
\begin{equation}\label{eq:3.1}
  I_{a,b}:= \{ x \in E_V \mid \la a,x\ra \cong \la b,x\ra \cong U_{\mathrm{2A}}\} .
\end{equation}
We will show that a 3-transposition group arises from $I_{a,b}$.

\begin{thm}\label{thm:3.5}
  Let $a,b\in V_\R$ be Ising vectors such that $\la a,b\ra \cong U_{\mathrm{3A}}$.
  Define $I_{a,b}$ as in \eqref{eq:3.1}.
  Then for any $x,y \in I_{a,b}$, $\la x,y\ra\cong U_{\mathrm{1A}}$, $U_{\mathrm{2A}}$ 
  or $U_{\mathrm{3A}}$.
\end{thm}

\pf
Let $x,y\in I_{a,b}$.
We eliminate the possibility of the type of $\la x,y\ra$ other than 1A, 2A, and 3A.
By our choice,  $\tau_x$ and $\tau_y$ commute with $\tau_a$ and $\tau_b$.
\\
~(i)~Case 6A:~
Since $\tau_a\tau_b$ is non-trivial on $\la a,b\ra \cong U_{\mathrm{3A}}$,
it is also non-trivial on $V$.
If $\la x,y\ra \cong U_{\mathrm{6A}}$, then $\tau_a=\tau_b$ by Lemma \ref{lem:3.4} and
we have a contradiction.
\medskip\\
~(ii)~Case 5A:~
Suppose $\la x,y\ra \cong U_{\mathrm{5A}}$.
Then $\tau_x\tau_y$ has an order $5$.
Since $\tau_a$ is non-trivial, $\tau_a\tau_x\tau_y=\tau_{a\circ x}\tau_y$ has
order 10 and we have a contradiction by the 6-transposition property in Theorem \ref{thm:2.7}.
\medskip\\
~(iii)~Cases 4A and 4B:~
Suppose $\la x,y\ra \cong U_{\mathrm{4A}}$ or $U_{\mathrm{4B}}$.
Then $\tau_x\tau_y$ has order 2 or 4 by Theorem \ref{thm:2.7}.
By (5) of Theorem \ref{thm:2.9},  we have 
\[
  \tau_{a\circ x}\tau_{b\circ y} =\tau_a\tau_x\tau_b\tau_y=\tau_a\tau_b\tau_x\tau_y.
\] 
If $\tau_x\tau_y$ has order 4, then $\tau_{a\circ x}\tau_{b\circ y}$ has order 12, 
which contradicts the 6-transposition property in Theorem \ref{thm:2.7}.
If $\tau_x\tau_y$ has order 2, then $\tau_{a\circ x}\tau_{b\circ y} $ has order 6 
and the subVOA generated by $a\circ x$ and $b\circ y$ is isomorphic to the 6A-algebra.
Since $\tau_{a\circ x}\tau_{b\circ y}$ has order 3 on $\la a\circ x, b\circ y\ra$, 
$(\tau_{a\circ x}\tau_{b\circ y})^3=\tau_x\tau_y$ acts trivially on $\la a\circ x, b\circ y\ra$.
Thus $\tau_x \tau_y$ fixes both $a\circ x$ and $b\circ y$.
Since $\tau_y$ fixes $b\circ y$, so does $\tau_x$ and hence $\tau_{b\circ y}$ fixes $x$.
Then $\tau_y$ fixes $x$ since $\tau_{b\circ y}=\tau_y\tau_b$ and $\tau_b$ fixes $x$.
This contradicts (6) of Theorem \ref{thm:2.7}.
Hence, $\la x,y\ra$ cannot be the 4A nor the 4B-algebra.
\medskip\\
~(iv)~Case 3C:~
That $\la x,y\ra \cong U_{\mathrm{3C}}$ is impossible by Lemma \ref{lem:3.2} since 
$\la x,y\ra \subset \la a,x,y\ra \cong U_{\mathrm{6A}}$.
\\[5pt]
~(v)~Case 2B:~  Suppose $\la x,y\ra \cong U_{\mathrm{2B}}$.
Then $x_{(1)}y=0$ and $(x \mymid y)=0$. It follows from $a\circ x= a+x-4a_{(1)}x$ that
\[
  (a\circ x \mymid y)
  = (a \mymid y)+(x \mymid y)-4(a_{(1)}x \mymid y)
  = (a \mymid y)-4(a \mymid x_{(1)}y)
  = 2^{-5}-4(a \mymid 0)
  = 2^{-5}.
\]
We consider $( a\circ x\mymid b\circ y) = (a+x-4a_{(1)}x\mymid b\circ y)
= (a\mymid b\circ y) +(x\mymid b\circ y) -4(x\mymid a_{(1)}(b\circ y) )$. 
By (7) of Theorem \ref{thm:2.13}, we have 
\[
  a_{(1)}(b\circ y) = \dfr{45}{2^{10}}u_{a,b}+\dfr{1}{2^5}\l( y+ a+ b\circ y
  -a\circ y- b -\tau_a b- (\tau_a b)\circ y \r),
\]
where $u_{a,b}$ is the characteristic $c=4/5$ Virasoro vector in $\langle a,b\rangle$.
Then we have 
\[
  (x\mymid a_{(1)}(b\circ y)) = -2\cdot2^{-10}
\]
since $(u_{a,b}\mymid x) =(y\mymid x)=0$ and  
$(a\mymid x)= ( b\circ y\mymid x) = (a\circ y\mymid x) = ( b\mymid x) = (\tau_a b\mymid x) 
= ((\tau_a b)\circ y\mymid x) = 2^{-5}$. 
Thus, 
\[
  ( a\circ x\mymid b\circ y) 
  =  (a\mymid b\circ y) +(x\mymid b\circ y) -4(x\mymid a_{(1)}(b\circ y) ) 
  = 5\cdot2^{-10} +2^{-5} +8\cdot 2^{-10} 
  = 45\cdot2^{-10},
\]
which is impossible by Theorem \ref{thm:2.6}.  

Therefore, $\la x,y\ra$ is isomorphic to one of $U_{\mathrm{1A}}$,
$U_{\mathrm{2A}}$ or $U_{\mathrm{3A}}$.
\qed


\medskip
 
As a corollary, we have 

\begin{cor}\label{cor:3.7}
  For $x,y\in I_{a,b}$, one has $\abs{\tau_x\tau_y}\leq 3$.
  Therefore, $G=\la \tau_x \mid x\in I_{a,b}\ra$ forms a 3-transposition group
  in the stabilizer $\{ g\in \aut(V_\R)\mid ga=a,~ gb=b\}$ of $\la a,b\ra$.
\end{cor}

\section{Inductive structures}\label{sec:4}

Let $V$ be an OZ-type VOA with a compact real form $V_\R$, and 
let $E_V$ be the set of Ising vectors in $V_\R$.
Suppose we have a pair $a$, $b \in E_V$ such that $\la a,b\ra$ is isomorphic to the 
3A-algebra, i.e., $(a \mymid b)=13\cd 2^{-10}$,  
and consider the subset $I_{a,b}$ defined as in \eqref{eq:3.1}.
Let $x$ and $y$ be distinct Ising vectors in $I_{a,b}$.
In this section, we first study the subalgebra $\la a,b,x,y\ra$ and then 
introduce a series of subalgebras $X^{[i]}=\la a,b,x^1,\dots,x^i\ra$ to study 
inductive structures of the 3-transposition group obtained 
by Corollary \ref{cor:3.7}. 
Note that $\langle x, y\rangle \cong U_{\mathrm{2A}}$ or $U_{\mathrm{3A}}$ by Theorem \ref{thm:3.5}. 

\begin{lem}\label{lem:4.1}
  Let $a$, $x$, $y$ be Ising vectors such that 
  $\la a , x\ra \cong \la a, y \ra\cong U_{\mathrm{2A}}$ 
  and $a\notin \la x,y\ra$.
  \\
  (1)~If $\la x,y\ra$ is the 2A-algebra then $(a \mymid x_{(1)}y)=2^{-6}$.
  \\
  (2)~If $\la x,y\ra$ is the 3A-algebra then $(a \mymid x_{(1)}y)=5\cd 2^{-9}$.
\end{lem}

\pf
(1):~Suppose $\la x,y\ra$ is the 2A-algebra.
Then $x_{(1)}y=\dfr{1}{4}(x+y-x\circ y)$.
By Lemma \ref{lem:2.10}, one has $(a \mymid x\circ y)=0$.
Therefore
\[
  (a \mymid x_{(1)}y)
  =\dfr{1}{4}(a \mymid x+y-x\circ y)=\dfr{1}{4}\cd 2\cd 2^{-5}=2^{-6}.
\]
(2):~Suppose $\la x,y\ra$ is the 3A-algebra.
Then $\la a,x,y\ra$ is the 6A-algebra with the central Ising vector $a$ by Lemma \ref{lem:3.2}.
By (5) of Theorem \ref{thm:2.11}, we have 
\[
  (a \mymid x_{(1)}y)=-\dfr{135}{2^{10}}(a \mymid u_{x,y})
  +\dfr{1}{2^4}(a \mymid 2x+2y+\tau_x y)
  =\dfr{1}{2^4}(2+2+1)\cd 2^{-5}=5\cd 2^{-9}.
\]
Thus we have the lemma.
\qed



\begin{prop}\label{prop:4.3}
  Let $a$, $b$, $x$, $y$ be Ising vectors such that
  $\la a,b\ra \cong U_{\mathrm{3A}}$ and 
  $\la a, x\ra \cong \la a,y\ra \cong \la b, x\ra \cong \la  b ,y\ra \cong U_{\mathrm{2A}}$.
  \\
  (1)~If $\la x,y\ra$ is the 2A-algebra then $\la a\circ x, b\circ y\ra \cong U_{\mathrm{6A}}$   , 
  i.e., $(a\circ x \mymid b\circ y)=5\cd 2^{-10}$.
  \\
  (2)~If $\la x,y\ra$ is the 3A-algebra then $\la a\circ x, b\circ y\ra \cong U_{\mathrm{3A}}$   , 
  i.e., $(a\circ x \mymid b\circ y)=13\cd 2^{-10}$ .
\end{prop}

\pf 
By definition, we have $a\circ x=a+x-4a_{(1)}x$ and $b\circ y=b+y-4b_{(1)}y$.
So by Lemma \ref{lem:4.1}, one has  
\[
\begin{array}{ll}
  (a\circ x \mymid b \circ y)
  &= (a+x-4a_{(1)}x \mymid b+y-4b_{(1)}y)
  \medskip\\
  &= \underbrace{(a\mymid b)}_{=13\cd 2^{-10}}
  +\underbrace{(a\mymid y)}_{=2^{-5}}
  -4\underbrace{(a_{(1)}b\mymid y)}_{=5\cd 2^{-9}}
  +\underbrace{(x \mymid b)}_{=2^{-5}}
  +(x \mymid y)-4(x_{(1)}y \mymid b)
  \medskip\\
  &~~~-4\underbrace{(x \mymid a_{(1)}b)}_{=5\cd 2^{-9}}
  -4(a\mymid x_{(1)}y)+16(x\mymid a_{(1)}b_{(1)}y)
  \medskip\\
  &= (x \mymid y)-4(a\mymid x_{(1)}y)-4(b\mymid x_{(1)}y)
     +16(x\mymid a_{(1)}b_{(1)}y)-3\cd 2^{-10}.
\end{array}
\]
Notice that   $\la a,b,y\ra$ is the 6A-algebra with 
the central Ising vector $y$ by Lemma \ref{lem:3.2}.  
By Lemma \ref{lem:2.14} (see also Eq.~\eqref{eq:2.10}), we have 
\begin{equation}\label{meq:4.1}
(x \mymid a_{(1)}b_{(1)}y)
    =2^{-7}\l( 7(x\mymid y)-7(x \mymid a\circ y)-(x \mymid b\circ y)
    +(x\mymid \tau_a b\circ y)+19\cd 2^{-5}\r).  
\end{equation} 
(1):~If $\la x,y\ra$ is the 2A-algebra, then $(x \mymid y)=2^{-5}$
and $(a \mymid x_{(1)}y)=(b \mymid x_{(1)}y)=2^{-6}$ by Lemma \ref{lem:4.1}.
Therefore
\[
  (a\circ x \mymid b\circ y)
  =16 (x \mymid a_{(1)}b_{(1)}y) -99\cd 2^{-10}.
\]
In this case,  $(x \mymid a\circ y)=0$ by Lemma \ref{lem:2.10} and $\tau_a$ acts trivially 
on $\la x,y\ra$.
Thus $\tau_a(b\circ y)=\tau_a b\circ \tau_a y =\tau_a b\circ y$ and 
$(x \mymid \tau_a b \circ y)=(x \mymid \tau_a(b\circ y))
=(\tau_a x \mymid b\circ y)=(x \mymid b\circ y)$.
Then by \eqref{meq:4.1}, one has
\[
  (x\mymid a_{(1)}b_{(1)}y)
  = 2^{-7}\l( 7\cd 2^{-5}+19\cd 2^{-5}\r)
  = 13\cd 2^{-11}.
\]
Therefore,
\[
  (a\circ x \mymid b\circ y)
  = 16\cd 13\cd 2^{-11}-99\cd 2^{-10}=5\cd 2^{-10}.
\]
(2):~If $\la x,y\ra$ is the 3A-algebra then $(x \mymid y)=13\cd 2^{-10}$
and $(a \mymid x_{(1)}y)=(b \mymid x_{(1)}y)=5\cd 2^{-9}$ by Lemma \ref{lem:4.1}.
Therefore
\[
  (a\circ x \mymid b\circ y)
  = 16(x \mymid a_{(1)}b_{(1)}y)-35\cd 2^{-9}.
\]
In this case,  $\la x,a\circ y\ra$, $\la x,b\circ y\ra$ and 
$\la x,\tau_a b \circ y\ra$ are all isomorphic to the 6A-algebra by Lemma \ref{lem:3.2}.
Thus $(x \mymid a\circ y)=(x \mymid b\circ y)=(x \mymid \tau_a b\circ y)=5\cd 2^{-10}$.
By \eqref{meq:4.1}, one has
\[
  (x \mymid a_{(1)}b_{(1)}y)
  =2^{-7}\l( 7\cd 13\cd 2^{-10}-7\cd 5\cd 2^{-10} +19\cd 2^{-5}\r)
  = 83\cd 2^{-14}.
\]
Therefore, 
\[
  (a\circ x \mymid b\circ y)
  = 16\cd 83\cd 2^{-14}-35\cd 2^{-9}
  = 13\cd 2^{-10}.
\]
This completes the proof.  
\qed
\medskip

The following lemma is a slight modification of Lemma 3 of \cite{Ma2} and 
is useful in determining the conformal vector of a subVOA.

\begin{lem}\label{lem:4.4}
  Let $V_\R$ be a compact VOA of OZ-type and let $A$ be a set of Virasoro vectors 
  of $V_\R$  such that the linear span of $A$ forms a subalgebra of the Griess algebra.
  Then the subalgebra $\la A\ra$ generated by $A$ has a conformal vector and 
  forms a subVOA of $V_\R$.
  The conformal vector $\eta$ of $\la A\ra$ belongs to the linear span of $A$ and 
  uniquely determined by the condition  $(\eta \mymid a)=(a \mymid a)$ for all $a\in A$.
\end{lem}

\pf
First, we prove that $\la A\ra$ has the conformal vector inside $\R A$, 
the linear span of $A$.
Let $\w$ be the conformal vector of $V_\R$ and let $B$ be the orthogonal complement 
of $\R A$ in $(V_\R)_2$.
Let $\w= \eta +\xi$ be the orthogonal decomposition with respect to 
$(V_\R)_2=\R A \perp B$.
For any $a\in A$, one has 
\[
  (a_{(1)} B \mymid \R A)=(B \mymid a_{(1)}\R A ) \subset (B \mymid \R A )=0.
\]
Therefore $(\R A)_{(1)}B\subset B$.
Then 
\[
  2\eta =\w_{(1)} \eta=(\eta+\xi)_{(1)} \eta=\eta_{(1)}\eta+\eta_{(1)}\xi 
\] 
and we obtain $\eta_{(1)}\xi=2\eta-\eta_{(1)}\eta \in B\cap \R A =0$.
Thus $\eta_{(1)}\eta=2\eta$ and $\eta_{(1)}\xi=0$.
This implies $\xi_{(1)}\xi=(\w-\eta)_{(1)}\xi= 2\xi$ and hence 
$\eta$ and $\xi$ are Virasoro vectors of $V_\R$ by Lemma 5.1 of \cite{M1}.
Since a Virasoro vector $\eta$ satisfies $\eta_{(0)}\eta= \eta_{(-2)}\vac$, we have 
\[
  \xi_{(0)}\eta
  = (\w-\eta)_{(0)}\eta
  = \w_{(0)}\eta-\eta_{(0)}\eta
  = \eta_{(-2)}\vac-\eta_{(-2)}\vac
  = 0
\]
and thus $\eta$ and $\xi$ are mutually orthogonal by Theorem 5.1 of \cite{FZ}.
We will show that $\eta$ is the conformal vector of $\la A\ra$.
For this, it suffices to show that $\eta_{(1)}a=2a$ and $\eta_{(0)}a=a_{(-2)}\vac$ 
for $a \in A$ since $\eta_{(0)}$ is a derivation on $V_\R$.
Let $a\in A$.
Then 
\[
  2a=\w_{(1)}a=(\eta+\xi)_{(1)}a=\eta_{(1)}a+\xi_{(1)}a
\]
and we have $\xi_{(1)}a=2a-\eta_{(1)}a \in B\cap \R A =0$.
Therefore $\eta_{(1)}a=2a$ and $\xi_{(1)}a=0$.
Next we will prove $\eta_{(0)}a=a_{(-2)}\vac$.
Since $\xi_{(2)}a\in (V_\R)_1=0$,  we have 
\[
  (\xi_{(0)}a \mymid \xi_{(0)}a)
  = (a \mymid \xi_{(2)}\xi_{(0)}a)
  = (a \mymid [\xi_{(2)},\xi_{(0)}]a)
  = (a \mymid 2\xi_{(1)}a)=0
\]
and $\xi_{(0)}a=0$ by the positive definiteness of $V_\R$.
Then $\eta_{(0)}a=(\w-\xi)_{(0)}a=\w_{(0)}a-\xi_{(0)}a=\w_{(0)}a=a_{(-2)}\vac$.
Therefore, $\eta$ is the conformal vector of $\la A\ra$.

Finally, we will show that $\eta \in \R A$ is uniquely determined by the condition 
$(\eta \mymid a)=(a \mymid a)$ for all $a\in A$.
Since $\eta_{(1)}a=a_{(1)}a=2a$, one has
\[ 
  2(a \mymid a)=(\eta_{(1)}a \mymid a)=(\eta \mymid a_{(1)}a)=2(\eta \mymid a)
\]
and the conformal vector $\eta \in \R A$ satisfies the condition.
Let $a_1,\dots,a_n\in A$ be a linear basis of $\R A$ and write 
$\eta=c_1a_1+\cds c_n a_n$ with $c_1,\dots,c_n \in \R$.
Then the Gram matrix $\begin{bmatrix} (a_i \mymid a_j)\end{bmatrix}_{1\leq i,j\leq n}$ 
is non-singular and the condition $(\eta \mymid a_i)=(a_i \mymid a_i)$ uniquely determines 
the coefficients $c_i\in \R$.
Therefore, the conformal vector of $\la A\ra$ is determined by the 
condition $\eta\in \R A$ and $(\eta \mymid a)=(a \mymid a)$ for all $a\in A$.
\qed

\begin{rem}\label{rem:4.5}
  From the proof of Lemma \ref{lem:4.4},  one sees that 
  $\la A\ra \subset \com_{V_\R} \la \xi\ra =\ker_{V_\R}\, \xi_{(0)}$ 
  (cf.~Theorem 5.1 of \cite{FZ}).
  It is shown in Proposition 4 of \cite{Ma2} that 
  $\R A= \la A\ra\cap (V_\R)_2$ if $\R A=(V_\R)_2\cap \ker_V\, \xi_{(1)}$.
\end{rem}

\subsection{The case $\la x,y\ra \cong U_{\mathrm{2A}}$}\label{sec:4.1}

In this subsection, we determine the Griess algebra of $\la a,b,x,y\ra$ when $\la x,y\ra$ 
is the 2A-algebra.

\begin{thm}\label{thm:4.6}
  Let $a$, $b$, $x$, $y$ be Ising vectors such that 
  $\la a,b\ra$ is isomorphic to the 3A-algebra,  
  and $\la x,y\ra$, $\la a,x\ra$, $\la a,y\ra$, $\la b,x\ra$ and $\la b,y\ra$ are 
  isomorphic to the 2A-algebra, i.e.,  
  $(x\mymid y)=(a \mymid x)=(a \mymid y)=(b \mymid x)=(b \mymid y)=2^{-5}$.
  Then the Griess subalgebra of $\la a,b,x,y\ra$ generated by $a$, $b$, $x$, $y$ 
  is 13-dimensional with a basis 
  \[
    A=\{ u=u_{a,b}, a, b, c=\tau_a b, a\circ x, b\circ x, c\circ x,
    a\circ y, b\circ y, c\circ y, x,y,z=x\circ y \} .
  \]
  The conformal vector of $\la a,b,x,y\ra$ is 
  \[
    \eta = \dfr{1}{6}u+\dfr{16}{27}(a+b+c+a\circ x +b\circ x +c\circ x
         +a\circ y +b\circ y +c\circ y)+\dfr{4}{9}(x+y+z)
  \]
  and has $52/15$.
  The VOA $\la a,b,x,y\ra$ has a full subVOA isomorphic to
  \[
    L(c_3,0)\tensor L(c_4,0)\tensor L(c_5,0)\tensor L(c_6,0)
    = L(\sfr{4}{5},0)\tensor L(\sfr{6}{7},0)\tensor L(\sfr{25}{28},0)
      \tensor L(\sfr{11}{12},0).
  \]
\end{thm}

\pf
Set $c=\tau_a b=\tau_b a$ and $z=x\circ y$.
We have seen in Lemma \ref{lem:3.2} that both $\la a,b\circ x\ra$ and 
$\la a,b\circ y\ra$ are isomorphic to the 6A-algebra and 
\[
\begin{array}{l}
  \la a,b\circ x\ra_2 
  = \Span \{ u_{a,b},x,a,b,c,a\circ x,b\circ x,c\circ x\} ,~~~
  \medskip\\
  \la a,b\circ y\ra_2 
  = \Span \{ u_{a,b},y,a,b,c,a\circ y,b\circ y,c\circ y\} .
\end{array}
\]
We have $u_{a,b}=u_{a\circ x,b\circ x}$ in $\la a,b\circ x\ra$ and 
$u_{a,b}=u_{a\circ y,b\circ y}$ in $\la a,b\circ y\ra$.
Therefore 
\[
  u=u_{a,b}=u_{a\circ x,b\circ x}=u_{a\circ y,b\circ y} .
\]
It follows from Lemma \ref{lem:2.10} that 
$(a \mymid x\circ y)= (x \mymid a\circ y)=(y\mymid a\circ x)=0$ 
so that $\sigma_a(x\circ y)=x\circ y$ and 
\[
  (a\circ x)\circ (a\circ y)
  = (\sigma_a x)\circ (\sigma_a y)
  = \sigma_a (x\circ y)
  = x\circ y.
\]
Therefore, the 3-set $\{ a\circ x, a\circ y, x\circ y\}$ forms a 2A-triple. 
Similarly, the 3-sets $\{ b\circ x, b\circ y, x\circ y\}$ and 
$\{ c\circ x, c\circ y, x\circ y\}$ also form 2A-triples.
By (1) of Proposition \ref{prop:4.3}, the Griess algebra of the subalgebra 
$\la a\circ x, b\circ y\ra$ is isomorphic to the 6A-algebra with a basis
\[
\begin{array}{l}
  a\circ x,~ 
  b\circ y,~ 
  \tau_{a\circ x}(b\circ y) = c\circ y,~ 
  \tau_{b\circ y}(a\circ x) = c\circ x,~ 
  \tau_{c\circ x}(b\circ y) = a\circ y,~
  \medskip\\
  \tau_{c\circ y}(a\circ x) = b\circ x,~
  (a\circ x)\circ (a\circ y)=x\circ y=z,~
  u_{a\circ x,b\circ x}=u_{a\circ y,b\circ y}=u_{a,b}.
\end{array}
\]
Therefore, the linear span of $A$ forms a subalgebra in the Griess algebra.
The determinant of the Gram matrix of $A$ is $3^{25}\cd 11^5/2^{81}\cd 5$ 
so that this matrix is non-singular. 
Thus $A$ is a basis of the Griess subalgebra generated by $a$, $b$, $x$ and $y$.

One can directly verify that the vector $\eta$ satisfies 
$(\eta \mymid t)=(t \mymid t)$ for all $t\in A$.
Therefore, $\eta$ is the conformal vector of $\la a,b,x,y\ra$ by Lemma \ref{lem:4.4}.
The central charge of $\eta$ is  $2(\eta \mymid \eta)=52/15$.
We know that the subalgebra $\la a,b,x\ra$ is isomorphic to the 6A-algebra and 
has a Virasoro frame 
$L(\sfr{4}{5},0)\tensor L(\sfr{6}{7},0)\tensor L(\sfr{25}{28},0)$ 
(see (8) of Theorem \ref{thm:2.13}).
By a direct computation, one finds that the commutant subalgebra of $\la a,b,x\ra$ 
in $\la a,b,x,y\ra$ is generated by the following $c=11/12$ Virasoro vector:
\begin{equation}\label{eq:4.1}
  f=-\dfr{5}{24}u-\dfr{2}{27}(a+b+c+a\circ x+b\circ x +c\circ x)
  +\dfr{16}{27}(a\circ y+b\circ y +c\circ y)
  -\dfr{1}{18}x +\dfr{4}{9}(y+z).
\end{equation}
Therefore, $\la a,b,x,y\ra$ has a unitary Virasoro frame isomorphic to 
\[
  L(c_3,0)\tensor L(c_4,0)\tensor L(c_5,0)\tensor L(c_6,0)
  = L(\sfr{4}{5},0)\tensor L(\sfr{6}{7},0)\tensor L(\sfr{25}{28},0)
  \tensor L(\sfr{11}{12},0).
\]
This completes the proof.
\qed

\subsection{The case $\la x,y\ra \cong U_{\mathrm{3A}}$}\label{sec:4.2}
Next we will determine the Griess algebra for 
$\langle a,b,x, y\rangle$ when $\la x,y\ra \cong U_{\mathrm{3A}}$. 
The result will be useful for analyzing the structures of the 3-transposition groups 
mentioned in Corollary \ref{cor:3.7}.
However, we are mainly interested in the inductive structures, i.e., 
the case when $\tau_x$ and $\tau_y$ commute. 
Thus, the results in this subsection will not be used in the later section. 
We include this subsection for completeness and for future references.

Suppose $\la a,b\ra$ and $\la x,y\ra$ are isomorphic to the 3A-algebra and  $\la a,x\ra$, $\la a,y\ra$, $\la b,x\ra$ and $\la b,y\ra$ are isomorphic the 2A-algebra, i.e., 
$(a \mymid x)=(a \mymid y)=(b \mymid x)=(b \mymid y)=2^{-5}$. 
Let $c=\tau_a b =\tau_b a$ and $z=\tau_x y=\tau_y x$. 
We also recall from \eqref{eq:2.8} that $u_{a,b}$ and $u_{x,y}$ denote the 
characteristic $c=4/5$ Virasoro vectors in the 3A-algebras $\la a,b \ra$ and 
$\la x,y\ra$, respectively.
Since all of $\la a,b, x\ra$, $\la a,b, y\ra$, $\la a,b, z\ra$, $\la a,x,y\ra$, 
$\la b,x, y\ra$ and $\la c,x, y\ra$ are isomorphic to the 6A-algebra 
by Lemma \ref{lem:3.2} and $U_{\mathrm{6A}}$ has the unique characteristic $c=4/5$ Virasoro vector, 
we have 
\begin{equation}\label{eq:4.2}
  u_{a,b}=u_{a\circ x, b\circ x} =u_{a\circ y, b\circ y}= u_{a\circ z, b\circ z},~~~
  u_{x,y}=u_{a\circ x, a\circ y}= u_{b\circ x, b\circ y}=u_{c\circ x, c\circ y}.
\end{equation}

Consider the subalgebra generated by $a\circ x$, $a\circ y$ and $b\circ x$.
Since $\tau_{a\circ x}=\tau_a\tau_x$, $\tau_{a\circ y}=\tau_a\tau_y$ and 
$\tau_{b\circ x}=\tau_b\tau_x$ by (5) of Theorem \ref{thm:2.11}, we have 
$\tau_{a\circ x}\tau_{a\circ y}=\tau_x\tau_y$ and $\tau_{a\circ x}\tau_{b \circ x}=\tau_a\tau_b$ 
and obtain the following conjugacy.
\begin{equation}\label{eq:4.3}
\begin{array}{ccccc}
  a\circ x &\xrightarrow{~\tau_x\tau_y~} &
  a\circ y &\xrightarrow{~\tau_x\tau_y~} &
  a\circ z 
  \\
  \mbox{\scriptsize $\tau_a\tau_b$}
  \mbox{\raisebox{2ex}[2ex][1ex]{\rotatebox{-90}{$\longrightarrow$}}}~~~~~
  &&
  \mbox{\scriptsize $\tau_a\tau_b$}
  \mbox{\raisebox{2ex}[2ex][1ex]{\rotatebox{-90}{$\longrightarrow$}}}~~~~~
  &&
  \mbox{\scriptsize $\tau_a\tau_b$}
  \mbox{\raisebox{2ex}[2ex][1ex]{\rotatebox{-90}{$\longrightarrow$}}}~~~~~
  \\
  b\circ x &\xrightarrow{~\tau_x\tau_y~} &
  b\circ y &\xrightarrow{~\tau_x\tau_y~} &
  b\circ z 
  \\
  \mbox{\scriptsize $\tau_a\tau_b$}
  \mbox{\raisebox{2ex}[2ex][1ex]{\rotatebox{-90}{$\longrightarrow$}}}~~~~~
  &&
  \mbox{\scriptsize $\tau_a\tau_b$}
  \mbox{\raisebox{2ex}[2ex][1ex]{\rotatebox{-90}{$\longrightarrow$}}}~~~~~
  &&
  \mbox{\scriptsize $\tau_a\tau_b$}
  \mbox{\raisebox{2ex}[2ex][1ex]{\rotatebox{-90}{$\longrightarrow$}}}~~~~~
  \\
  c\circ x &\xrightarrow{~\tau_x\tau_y~} &
  c\circ y &\xrightarrow{~\tau_x\tau_y~} &
  c\circ z 
\end{array}
\end{equation}
Set $H=\la \tau_{a\circ x},\tau_{a\circ y}, \tau_{b\circ x}\ra$.
Then $H\cong 3^2{:}2$ acts transitively on the the set of 9 Ising vectors 
in \eqref{eq:4.3}.
By (2) of Proposition \ref{prop:4.3}, we see that $(a\circ x \mymid b\circ y)=13\cdot 2^{-10}$ 
and hence $\la a\circ x, b\circ y\ra $ is also a 3A-algebra. 
We have $\tau_{a\circ x}(b\circ y)=\tau_a b\circ \tau_x y=c\circ z$ and 
$\{a\circ y, b\circ x, c \circ z\}$ forms the 3A-triple of $\la a\circ y,b\circ x\ra$.
Similarly, by using the conjugacy relations we can verify that each pair of Ising vectors 
in \eqref{eq:4.3} generate a 3A-algebra and we obtain the following 
four $H$-orbits of 3A-triples:
\begin{equation}\label{eq:4.4}
\begin{array}{c}
  \mathcal{L}_1=\{\{ a\circ x, b\circ x, c\circ x\},
  \{ a\circ y, b\circ y, c\circ y\},
  \{ a\circ z, b\circ z, c\circ z\} \} ,
  \\
  \mathcal{L}_2=\{\{ a\circ x, a\circ y, a\circ z\} ,
  \{ b\circ x, b\circ y, b\circ z\},
  \{ c\circ x, c\circ y, c\circ z\} \} ,
  \\
  \mathcal{L}_3=\{\{ a\circ x, b\circ y, c\circ z\},
  \{ a\circ y, b\circ z, c\circ x\},
  \{ a\circ z, b\circ x, c\circ y\} \} ,
  \\
  \mathcal{L}_4=\{
  \{ a\circ x, b\circ z, c\circ y\},
  \{ a\circ y, b\circ x, c\circ z\},
  \{ a\circ z, b\circ y, c\circ x\} \} .
\end{array}
\end{equation}
The configuration of 9 points given by Ising vectors in \eqref{eq:4.3} together 
with the lines given by 3-sets forming 3A-triples is isomorphic to the affine plane 
of order 3 and each $H$-orbit $\mathcal{L}_i$ in \eqref{eq:4.4} is the set of 
parallel lines.
By \eqref{eq:4.2}, all 3A-triples in $\mathcal{L}_1$ define the same 
characteristic $c=4/5$ Virasoro vector.
The subVOA generated by such 9 Ising vectors is determined in \cite{LSo}.

\begin{lem}[Lemma 4.22 and Proposition 4.24 of \cite{LSo}]\label{lem:4.7}
  For $1 \leq i\leq 4$, all 3A-triples in $\mathcal{L}_i$ define the same 
  characteristic $c=4/5$ Virasoro vector $u^i$. Moreover, the $c=4/5$ Virasoro 
  vectors $u^i$, $1\leq i\leq 4$, are mutually orthogonal.
\end{lem}

As a consequence, we have 
\begin{equation}\label{eq:4.5}
  u_{a\circ x,b\circ y}=u_{a\circ y,b\circ z}=u_{a\circ z,b\circ x},\quad \text{ and } \quad 
  u_{a\circ x,b\circ z}=u_{a\circ y,b\circ x}=u_{a\circ z,b\circ y}.
\end{equation}
It follows from the relations above that the linear span of the Ising vectors 
in \eqref{eq:4.3} together with the $c=4/5$ Virasoro vectors $u^1=u_{a,b}$, $u^2=u_{x,y}$, 
$u^3=u_{a\circ x,b\circ y}$ and $u^4=u_{a\circ x,b\circ z}$ forms a subalgebra in 
the Griess algebra.
The products and inner products of this Griess subalgebra is described 
by using the 3A-algebras along with the incidence structure of the affine plane of order 3.
More precisely, the subalgebra generated by Ising vectors in \eqref{eq:4.3} is 
described as follows.

\begin{thm}[Lemma 3.8 of \cite{HLY2} and Theorem 4.25 of \cite{LSo}]\label{thm:4.8}
  The subVOA generated by 9 Ising vectors in \eqref{eq:4.3} is isomorphic the 
  ternary code VOA $M_{\mathcal{C}}$ constructed in \cite{KMY}.
  It is actually generated by 3 Ising vectors $a\circ x$, $a\circ y$ and $b\circ x$ and 
  its Griess algebra is 12-dimensional spanned by 9 Ising vectors in \eqref{eq:4.3} and 
  4 mutually orthogonal $c=4/5$ Virasoro vectors $u_{a,b}$, $u_{x,y}$, 
  $u_{a\circ x,b\circ y}$ and $u_{a\circ x,b\circ z}$ with one linear relation
  \[
  \begin{array}{l}
    u_{a,b}+ u_{x,y}+u_{a\circ x,b\circ y}+u_{a\circ x,b\circ z}
    \medskip\\
    =
    \dfr{32}{45}( a\circ x+b\circ x++c\circ x+ a\circ y+b\circ y+ c\circ y+ 
    a\circ z+b\circ z+ c\circ z). 
  \end{array}
  \]
Moreover, $ \omega=u_{a,b}+ u_{x,y}+u_{a\circ x,b\circ y}+u_{a\circ x,b\circ z}$ gives the conformal vector of the subalgebra.
\end{thm}

Set 
\begin{equation}\label{eq:4.6}
\begin{array}{l}
  A=\{ u_{a,b},~u_{x,y},~ u_{a\circ x,b\circ y},~ u_{a\circ x,b\circ z},~
    a,~b,~c=\tau_a b,~
    x,~y,~z=\tau_x y,~
    \\
~~~~~~~~~~~~~~~~~~~~~
    a\circ x,~ a\circ y,~ a\circ z,~ 
    b\circ x,~ b\circ y,~ c\circ z,~ 
    c\circ x,~ c\circ y,~ c\circ z\} .
\end{array}
\end{equation}
In the rest of this subsection, we will prove that the linear span of $A$ 
forms a subalgebra in the Griess algebra.
By Lemma \ref{lem:3.2},  the following subalgebras are isomorphic to the Griess algebra of the 6A-algebra:
\begin{equation}\label{eq:4.7}
\begin{array}{l}
  \la a,b,x\ra_2=\Span\{ u_{a,b},x,a,b,c,a\circ x, b\circ x,c\circ x\}, 
  \medskip\\
  \la a,b,y\ra_2=\Span\{ u_{a,b},y,a,b,c,a\circ y, b\circ y,c\circ y\}, 
  \medskip\\
  \la a,b,z\ra_2=\Span\{ u_{a,b},z,a,b,c,a\circ z, b\circ z,c\circ z\}, 
  \medskip\\
  \la a,x,y\ra_2=\Span\{ u_{x,y},a,x,y,z,a\circ x, a\circ y,a\circ z\}, 
  \medskip\\
  \la b,x,y\ra_2=\Span\{ u_{x,y},b,x,y,z,b\circ x, b\circ y,b\circ z\}, 
  \medskip\\
  \la c,x,y\ra_2=\Span\{ u_{x,y},c,x,y,z,c\circ x, c\circ y,c\circ z\} .
\end{array}
\end{equation}
Therefore, it is enough to show that the products of $a$, $b$, $c$, $x$, $y$, $z$ 
and $u_{a\circ x,b\circ y}$, $u_{a\circ x,b\circ z}$ lie in the linear span of $A$ and 
determine the inner products among them.
Set $G=\la \tau_a,\tau_b,\tau_x,\tau_y\ra$.
Then the sets of 3A-triples $\mathcal{L}_1$ and $\mathcal{L}_2$ in \eqref{eq:4.4} are 
$G$-stable, whereas $\mathcal{L}_3$ and $\mathcal{L}_4$ are $H$-stable but not $G$-stable.
We have the following conjugacy.

\begin{lem}\label{lem:4.9}
  For $e\in \{ a,b,c,x,y,z\}$, 
  we have $\tau_e u_{a\circ x, b\circ y}= u_{a\circ x, b\circ z}$.
\end{lem}

\pf
Let $e\in \{ a,b,c,x,y,z\}$.
Then we have $\tau_e \mathcal{L}_3=\mathcal{L}_4$ and the conjugacy 
$\tau_e u_{a\circ x,b\circ y}=u_{a\circ x,b\circ z}$ of the characteristic 
$c=4/5$ Virasoro vectors follows from that of the corresponding 3A-triples.
\qed

\begin{lem}\label{lem:4.10}
  Let $e\in \{x,y,z, a,b,c\}$. 
  Then
  \[
    \l( e \mymid u_{a\circ x,b\circ y}\r) 
    = \l( e \mymid u_{a\circ x,b\circ z}\r) =\frac{1}{80}.
  \]
\end{lem}

\pf
Recall that
\[
  u_{a\circ x, b\circ y}
  = \frac{2^6}{135} \l( 2 a\circ x+ 2 b\circ y + c\circ z - 16(a\circ x)_{(1)}(b\circ y)\r).
\]
First we have
\[
\begin{array}{ll}
  \l( a \mymid (a\circ x)_{(1)} (b\circ y)\r) 
  & = \l( a_{(1)}(a\circ x) \mymid b\circ y\r) 
  = \dfr{1}{4} \l( a+a\circ x -x \mymid b\circ y\r) 
  \medskip\\
  &= \dfr{1}{4} \l(\dfr{5}{2^{10}}+ \dfr{13}{2^{10}} -\dfr{5}{2^{10}} \r) 
  = \dfr{13}{2^{12}}.
\end{array}
\]
Hence,
\[
  \l( a \mymid u_{a\circ x, b\circ y}\r)  
  = \dfrac{2^6}{135} \l( 2\cd\dfrac{1}{2^5} +2\cdot \dfrac{5}{2^{10}}+ \dfrac{5}{2^{10}} 
  -16\cdot \dfrac{13}{2^{12}}\r)
  = \dfrac{1}{80}.
\]
The other cases can be proved similarly.
\qed
\medskip

Let us compute $e_{(1)}u$ for $e\in \{ a,b,c,x,y,z\}$ and 
$u\in \{ u_{a\circ x,b\circ y}, u_{a\circ x,b\circ z}\}$. 
Recall that $a\circ x=2^{-2}(a+x-a_{(1)}x)$  in \eqref{eq:2.7}. 
For simplicity, we regard $a\circ x$ as a bilinear product and define 
$a\circ (x+y+z)=a\circ x+a\circ y+a\circ z$.

\begin{lem}\label{lem:4.11}
  Let $u= u_{a\circ x, b\circ y}$ or $u_{a\circ x, b\circ z}$. 
  Then we have
  \[
  \begin{array}{ll}
    \sigma_e(u+\tau_eu)
    &= -\dfr{1}{2}u_{a,b}-u_{x,y} 
    +\dfrac{4}{15}e -\dfr{8}{45}e\circ (x+y+z)
    \medskip\\
    &
    ~~~~-\dfrac{16}{45}(a+b+c) +\dfrac{8}{15}(x+y+z) +\dfrac{32}{45}(a +b +c)\circ (x+y+z) 
  \end{array}
  \]
  if $e\in \{ a,b,c\}$, and 
  \[
  \begin{array}{ll}
    \sigma_e(u+\tau_e u)
    &= -u_{a,b}-\dfr{1}{2}u_{x,y} 
    +\dfrac{4}{15}e -\dfr{8}{15}(a+b+c)\circ e 
    \medskip\\
    &
    ~~~~+\dfrac{8}{15}(a+b+c) -\dfrac{16}{45}(x+y+z) +\dfrac{32}{45}(a +b +c)\circ (x+y+z) 
  \end{array}
  \]
  if $e\in \{ x,y,z\}$.
\end{lem}

\pf
We only compute the case $e=x$.
The case $e=a$ is similar.
By \eqref{eq:2.8}, we have
\[
\begin{array}{lll}
  u_{a\circ x, b\circ y} 
  &=& \dfrac{2^6}{135} \l( 2a\circ x +2b\circ y +c\circ z-16(a\circ x)_{(1)}(b\circ y)\r) ,
  \medskip\\
  u_{a\circ x, b\circ z} 
  &=& \dfrac{2^6}{135} \l( 2a\circ x +2b\circ z +c\circ y-16(a\circ x)_{(1)}(b\circ z)\r) .
\end{array}
\]
Since $\sigma_x(a\circ x)=a$ and 
$u+\tau_x u= u_{a\circ x, b\circ y}+ u_{a\circ x, b\circ z}$, 
we have 
\begin{equation}\label{eq:4.8}
\begin{array}{l}
  \sigma_x(u+\tau_xu)
  = \sigma_x( u_{a\circ x, b\circ y}+ u_{a\circ x, b\circ z})
  \medskip\\
  = \dfrac{2^6}{135}\sigma_x \left( 4a\circ x +2b\circ y + 2b\circ z + c\circ y+c\circ z
    -16(a\circ x)_{(1)}(b\circ y+b\circ z) \r)
  \medskip\\
  = \dfrac{2^6}{135}\left( 4a + 2\sigma_x(b\circ y+b\circ z) 
    + \sigma_x(c\circ y+c\circ z) -16 a_{(1)} \sigma_x(b\circ y+b\circ z)\right).
\end{array}
\end{equation}
By Lemma \ref{lem:2.14} (cf.~Eq.~\eqref{eq:2.10}), we have
\begin{equation}\label{eq:4.9}
\begin{array}{l}
  \sigma_x(b\circ y+b\circ z) 
  = -\dfrac{45}{2^7} u_{x,y} -\dfrac{1}{2^2}b +\dfrac{1}{2^4}x +\dfrac{1}{2^2}b\circ x 
    + \dfrac{1}{2^2}(y+z) +b\circ y+b\circ z,
  \medskip\\
  \sigma_x(c\circ y+c\circ z) 
  = -\dfrac{45}{2^7} u_{x,y} -\dfrac{1}{2^2}c +\dfrac{1}{2^4}x +\dfrac{1}{2^2}c\circ x 
    +\dfrac{1}{2^2}(y+z) +c\circ y+c\circ z.
\end{array}
\end{equation}
By using the Griess algebras of $\la a,b\circ y\ra$ and $\la a,b\circ z\ra$ 
(cf.~(7) of Theorem \ref{thm:2.13}), we have 
\begin{equation}\label{eq:4.10}
\begin{array}{l}
  a_{(1)}\, \sigma_x(b\circ y+b\circ z) 
  \medskip\\ 
  = a_{(1)}\l( -\dfrac{45}{2^7} u_{x,y} -\dfrac{1}{2^2}b +\dfrac{1}{2^4}x 
  +\dfrac{1}{2^2}b\circ x + \dfrac{1}{2^2}(y+z) +b\circ y+b\circ z\r)
  \medskip\\
  = -\dfrac{45}{2^7}a_{(1)}u_{x,y}
  -\dfrac{1}{2^2}\left( -\dfrac{135}{2^{10}}u_{a,b}+\dfrac{1}{2^4}(2a+2b+c) \right) 
  +\dfrac{1}{2^4}\cdot \dfrac{1}4(a+x-a\circ x) 
  \medskip\\
  ~~~~
  +\dfrac{1}{2^2}\left( \dfrac{45}{2^{10}}u_{a,b}
  +\dfrac{1}{2^5}(x+a+b\circ x -b-c-a\circ x -c\circ x) \right)
  \medskip\\
  ~~~~
  +\dfrac{1}{2^2}\left( \dfrac{1}4(a+y-a\circ y)+\dfrac{1}4(a+z-a\circ z)\right)
  \medskip\\
  ~~~~
  +\dfrac{45}{2^{10}} u_{a,b}
  + \dfrac{1}{2^5}\left(y+a+b\circ y-b -c-a\circ y-c\circ y\right) 
  \medskip\\
  ~~~~
  +\dfrac{45}{2^{10}} u_{a,b}
  + \dfrac{1}{2^5}\left(z+a+b\circ z-b -c-a\circ z-c\circ z\right) .
\end{array}
\end{equation}
Note that $a_{(1)}u_{x,y} =0$ in the 6A-algebra $\la a,x,y\ra$ and thus by plugging 
\eqref{eq:4.9} and \eqref{eq:4.10} into \eqref{eq:4.8},
we have
\[
\begin{array}{ll}
  \sigma_x( u+\tau_xu)
  &= \dfrac{8}{15}(a+b+c) -\dfrac{4}{45}x -\dfrac{16}{45}(y+z) 
    +\dfrac{8}{45}(a\circ x + b\circ x+c\circ x)
  \medskip\\
  &~~~~ + \dfrac{32}{45}(a\circ y + b\circ y+c\circ y+ a\circ z + b\circ z+c\circ z) 
  -u_{a,b}-\dfrac{1}2u_{x,y}
\end{array}
\]
as desired.
\qed

\begin{lem}\label{lem:4.12}
  Let $u= u_{a\circ x, b\circ y}$ or $u_{a\circ x, b\circ z}$. 
  Then we have
  \[
  \begin{array}{ll}
    e_{(1)} u 
    &= \dfrac{1}{16}u_{a,b} +\dfrac{1}{8} u_{x,y} +\dfrac{5}{32}u +\dfrac{3}{32}\tau_eu  
    +\dfrac{1}{15}e +\dfrac{1}{15}e\circ (x+y+z)
    \medskip\\
    &~~~~  
    +\dfrac{2}{45}(a+b+c) -\dfrac{1}{15}(x+y+z) -\dfrac{4}{45}(a+b+c)\circ (x+y+z)
  \end{array}
  \]
  if $e\in \{ a,b,c\}$, and 
  \[
  \begin{array}{ll}
    e_{(1)} u 
    &= \dfrac{1}{8}u_{a,b} +\dfrac{1}{16} u_{x,y} +\dfrac{5}{32}u +\dfrac{3}{32}\tau_eu  
    +\dfrac{1}{15}e +\dfrac{1}{15}(a+b+c)\circ e
    \medskip\\
    &~~~~  
    -\dfrac{1}{15}(a+b+c) +\dfrac{2}{45}(x+y+z) -\dfrac{4}{45}(a+b+c)\circ (x+y+z)
  \end{array}
  \]
  if $e\in \{ x,y,z\}$.
\end{lem}

\pf
Again we only compute the case $e=x$.
By \eqref{eq:2.5}, we have
\[
  x_{(1)}u
  =8(x \mymid u)x +\dfr{5}{32}u+\dfr{3}{32}\tau_x u-\dfr{1}{8}\sigma_x(u+\tau_x u) .
\]
Then the lemma immediately follows from Lemmas \ref{lem:4.10} and \ref{lem:4.11}.
\qed

Recall the set $A$ of Virasoro vectors in \eqref{eq:4.6}.
By Lemma \ref{lem:4.12} and $H$-invariance of $u_{a\circ x,b\circ y}$ and 
$u_{a\circ x,b\circ z}$, the products $e_{(1)}u$ for $e\in \{ a,b,c,x,y,z\}$ and 
$u\in \{ u_{a\circ x,b\circ y},u_{a\circ x,b\circ z}\}$ are completely determined 
by linear combinations of the Virasoro vectors in $A$.
Therefore, the linear span of $A$ forms a subalgebra of the Griess algebra of 
$\la a,b,x,y\ra$.
\medskip

Next we determine the conformal vector of $\la a,b,x,y\ra$.

\begin{prop}\label{prop:4.13}
  Set 
  \[
    \eta 
    = \dfr{17}{22}(u_{a,b}+u_{x,y})
    +\dfr{10}{11}(u_{a\circ x,b\circ y}+u_{a\circ x,b\circ z})
    +\dfr{16}{33}(a+b+c+x+y+z).
  \]
  Then $(\eta \mymid t)=(t\mymid t)$ for all $t\in A$ and therefore $\eta$ is 
  the conformal vector of $\la a,b,x,y\ra$ of central charge $228/55$.
\end{prop}

\pf 
Since $u_{a,b}$, $u_{x,y}$, $u_{a\circ x,b\circ y}$, $u_{a\circ x,b\circ z}$ are mutually 
orthogonal by Lemma \ref{lem:4.7}, it is straightforward to verify 
$(\eta\mymid t)=(t\mymid t)$ for all $t\in A$ by using the 6A-algebras given 
in \eqref{eq:4.7} and Lemma \ref{lem:4.10}.
By Lemma \ref{lem:4.4}, $\eta$ is the conformal vector of $\la a,b,x,y\ra$. 
The central charge is given by 
\[
\begin{array}{ll}
  2(\eta|\eta)
  &=\dfr{17}{22}(\eta\mid u_{a,b}+u_{x,y})
    +\dfr{10}{11}(\eta \mid u_{a\circ x,b\circ y}+u_{a\circ x,b\circ z})
    +\dfr{16}{33}(\eta \mid a+b+c+x+y+z)
  \medskip\\
  &=2\cd\l(\dfr{17}{22}\cd \dfr{2}{5}\cd 2+\dfr{10}{11}\cd \dfr{2}{5}\cd 2
  +\dfr{16}{33}\cd \dfr{1}{4}\cd 6\r) 
  = \dfr{228}{55} 
\end{array}
\]
as claimed.
\qed

Summarizing everything, the structure of the Griess algebra of $\la a,b,x,y\ra$ is 
described as follows.

\begin{thm}\label{thm:4.14}
  Let $a$, $b$, $x$, $y$ be Ising vectors of $V_\R$ such that $\la a,b\ra$ and $\la x,y\ra$ 
  are isomorphic to the 3A-algebra and $\la a,x\ra$, $\la a,y\ra$, $\la b,x\ra$ and 
  $\la b,y\ra$ are isomorphic to the 2A-algebra,
  i.e.,  $(a \mymid x)=(a \mymid y)=(b \mymid x)=(b \mymid y)=2^{-5}$.
  Then the Griess subalgebra of $\la a,b,x,y\ra$ generated by $a$, $b$, $x$ and $y$ 
  is 18-dimensional spanned by 
  \[
  \begin{array}{l}
    A=\{ u_{a,b},~u_{x,y},~ u_{a\circ x,b\circ y},~ u_{a\circ x,b\circ z},~
      a,~b,~c=\tau_a b,~
      x,~y,~z=\tau_x y,~
      \\
    ~~~~~~~~~~~~~~~~~~~~~
    a\circ x,~ a\circ y,~ a\circ z,~ 
    b\circ x,~ b\circ y,~ c\circ z,~ 
    c\circ x,~ c\circ y,~ c\circ z\} 
  \end{array}
  \]
  with one linear relation
  \[
    u_{a,b}+u_{x,y}+u_{a\circ x,b\circ y}+u_{a\circ x,b\circ z}
    = \dfr{32}{45}(a\circ x + a\circ y + a\circ z + b\circ x + b\circ y + c\circ z, 
    c\circ x, c\circ y, c\circ z).
  \]
  The conformal vector of $\la a,b,x,y\ra$ is 
  \[
    \eta 
    = \dfr{17}{22}(u_{a,b}+u_{x,y})
    +\dfr{10}{11}(u_{a\circ x,b\circ y}+u_{a\circ x,b\circ z})
    +\dfr{16}{33}(a+b+c+x+y+z) 
  \]
  and has central charge $228/55$.
  The VOA $\la a,b,x,y\ra$ has a full subVOA isomorphic to 
  \[
    L(c_3,0)\tensor L(c_3,0)\tensor L(c_3,0)\tensor L(c_3,0)\tensor L(c_8,0)
    =L(\sfr{4}{5},0)^{\tensor 4}\tensor L(\sfr{52}{55},0).
  \] 
\end{thm}

\pf
We have already shown that the Griess subalgebra generated by $a$, $b$, $x$, $y$ is 
linearly spanned by $A$.
We also know the linear relation in Theorem \ref{thm:4.8}.
It is straightforward to verify that the determinant of the 
Gram matrix of $A\setminus\{ u_{a\circ x,b\circ z}\}$ is equal to 
$3^{52}\cd 11\cd 13^6/2^{138}\cd 5^3$. 
Therefore, the linear span of $A$ is 18-dimensional.
The conformal vector of $\la a,b,x,y\ra$ is given in Proposition \ref{prop:4.13}.

The subalgebra $\la a\circ x,a\circ y,b\circ x\ra$ is isomorphic to the ternary 
code VOA and has the Virasoro frame 
$u_{a,b}+u_{x,y}+u_{a\circ x,b\circ y}+u_{a\circ x,b\circ z}$.
That is, $\la a\circ x,a\circ y,b\circ x\ra$ is a $\Z_3\oplus \Z_3$-graded simple 
current extension of $L(\sfr{4}{5},0)^{\tensor 4}$.
Then the element  
\begin{equation}\label{eq:4.11}
\begin{array}{ll}
  \xi
  &=\eta-(u_{a,b}+u_{x,y}+u_{a\circ x,b\circ y}+u_{a\circ x,b\circ z})
  \medskip\\
  &= \dfr{5}{22}(u_{a,b}+u_{x,y})+\dfr{1}{11}(u_{a\circ x,b\circ y}+u_{a\circ x,b\circ z})
  +\dfr{16}{33}(a+b+c+x+y+z)
\end{array}
\end{equation}
is a Virasoro vector orthogonal to 
$u_{a,b}+u_{x,y}+u_{a\circ x,b\circ y}+u_{a\circ x,b\circ z}$.
The central charge of $\xi$ is $228/55-4\cd 4/5=52/55=c_8$.
Therefore, $\la a,b,x,y\ra$ has the full subVOA 
\[
  \la a\circ x,a\circ y,b\circ x,\xi\ra
  \cong \la a\circ x,a\circ y,b\circ x\ra\tensor \la\xi\ra
\]
which contains a unitary Virasoro frame 
$L(\sfr{4}{5},0)^{\tensor 4}\tensor L(\sfr{52}{55},0)$.
\qed

\subsection{Inductive subalgebras}\label{sec:4.3}

Let $a$, $b$ be Ising vectors of $V_\R$ such that 
$\la a,b\ra$ is isomorphic to the 3A-algebra, i.e.,  
$(a \mymid b)=13\cd 2^{-10}$ and let $I_{a,b}$ be defined as in \eqref{eq:3.1}. 
It is shown in Corollary \ref{cor:3.7} that $\la \tau_x \mid x\in I_{a,b}\ra$ is 
a 3-transposition group.
In order to study the inductive structures and basic sets of this group, 
we will choose $x^i\in I_{a,b}$ and define subalgebras $X^{[i]}=\la X^{[i-1]},x^i\ra$ 
inductively as follows.

\begin{df}\label{df:4.15}
Set $X^{[0]}=\la a,b\ra$. Suppose we have chosen  $x^1,\dots,x^i \in I_{a,b}$ and defined 
\begin{equation}\label{eq:4.12}
  X^{[i]}:=\la a,b,x^1,\dots,x^i\ra \quad  \text{ for some }~ i\geq 0.
\end{equation}
We choose $x^{i+1}\in I_{a,b}$ such that $x^{i+1}\not\in X^{[i]}$ and 
$(x^{i+1}\mymid x^j)=2^{-5}$ for all $1\leq j\leq i$ and define the subalgebra $X^{[i+1]}$ by  
\[
  X^{[i+1]}:=\langle X^{[i]},x^{i+1}\rangle=\la a,b,x^1,\dots,x^i, x^{i+1}\ra 
\]
as long as possible.
\end{df}

We will show that the structure of the Griess algebra of $X^{[n]}$ 
does not depend on the choice of $x^i$, $1\leq i\leq n$, and is uniquely determined up to 
isomorphism.
First, we observe from \cite{Ma2} that the Griess algebra of 
$\la x^1,\dots,x^n\ra$ is uniquely determined by our choice. 

\begin{prop}[\cite{Ma2}]\label{prop:4.15}
  Suppose $x^1,\cds,x^n$ are Ising vectors such that 
  $(x^i\mymid x^j)=2^{-5}$ and $x^k\not \in \la x^1,\dots,x^{k-1}\ra$
  for any  $1\leq i<j\leq n$ and $1<k\leq n$.
  Then the Griess algebra of $\la x^1,\dots,x^n\ra$ has a unique structure
  with a basis $\{ x^i,~x^j\circ x^k \mid 1\leq i\leq n,~1\leq j<k\leq n\}$.
\end{prop}

\pf
Let $i$, $j$, $k$ be distinct.
By Lemma \ref{lem:2.10},  we have $(x^i\circ x^j| x^k)=0$ and 
$\la x^i\circ x^j ,x^k\ra$ is the 2B-algebra.
Therefore $x^k_{(1)}\l( x^i\circ x^j\r)=0$ and $\sigma_{x^k}$ fixes $x^i\circ x^j$.
Then by conjugating $\{ x^i,x^j,x^i\circ x^j\}$ by $\sigma_{x^k}$ we see that 
$\{ x^i\circ x^j,x^i\circ x^k,x^j\circ x^k\}$ forms a 2A-triple.
Moreover, if $l$ is distinct from $i$, $j$ and $k$, then 
\[
  (x^i\circ x^j \mid x^k\circ x^l)
  = \l( x^i +x^j-4x^i_{(1)}x^j \,\big|\, x^k\circ x^l\r)
  = -4\l( x^i\,\Big|\,x^j_{(1)}(x^k\circ x^l)\r)=0.
\]
Therefore, $\la x^i\circ x^j, x^k\circ x^l\ra$ is a 2B-algebra.
Hence, the Griess subalgebra generated by $x^i$, $1\leq i\leq n$, is uniquely determined 
and linearly spanned by Ising vectors $x^i$, $1\leq i\leq n$ and $x^j\circ x^k$, 
$1\leq j<k\leq n$.
Such Ising vectors can be realized inside $V_{\sqrt{2}A_n}^+$ and known to be 
linearly independent (cf.~\cite{DLMN,LSY,Ma1}).
This completes the proof.
\qed

\begin{rem}\label{rem:4.16}
The Griess algebra of $\la x^1,\dots,x^n\ra$ coincides with that of 
$M_{A_n}$ in \cite{LSY} (cf.~Proposition \ref{prop:a.3}).
Set $\tilde{x}^1=x^1$ and $\tilde{x}^i=x^i\circ x^{i-1}$ for $2\leq i\leq n$.
It is clear that $\la x^1,\dots,x^n\ra=\la \tilde{x}^1,\dots,\tilde{x}^n\ra$. 
We have $(\tilde{x}^i \mymid \tilde{x}^j)=0$ if $\abs{i-j}>1$ and 
$(\tilde{x}^k \mymid \tilde{x}^{k+1})=2^{-5}$ for $1\leq k<n$.
Then the associated involutions $\sigma_{\tilde{x}^1},\dots,\sigma_{\tilde{x}^n}$ 
acting on the subalgebra $\la \tilde{x}^1,\dots,\tilde{x}^n\ra$ satisfy 
the Coxeter relation of type $A_n$ and hence 
$\la \sigma_{\tilde{x}^1},\dots,\sigma_{\tilde{x}^n}\ra$ is isomorphic to the 
symmetric group $\mathrm{S}_{n+1}$.
If we identify this group as a permutation group of the $(n+1)$-set 
$\{ 0,1,2,\dots,n\}$, then the involutions $\sigma_{x^i}$ and 
$\sigma_{x^j\circ x^k}=\sigma_{x^j}\sigma_{x^k}\sigma_{x^j}$ 
correspond to the transpositions $(0~i)$ and $(j~k)$, respectively.
\end{rem}

\begin{thm}\label{thm:4.17}
  With reference to the above, the Griess algebra generated by 
  Ising vectors $a,b,x^1,\dots,x^n$ is uniquely determined
  and has the following basis:
  \[
    u=u_{a,b},~
    a,~b,~c=\tau_a b,~
    x^i,~a\circ x^i,~ b\circ x^i,~c\circ x^i,~ 
    x^j\circ x^k,~~ 1\leq i\leq n,~~ 1\leq j<k\leq n.
  \]
  The conformal vector of $X^{[n]}=\la a,b,x^1,\dots,x^n\ra$ is given by 
  \begin{equation}\label{eq:4.13}
  \begin{array}{l}
    \w^n = \dfr{3(3-n)}{2(n+7)}u + \dfr{16}{3(n+7)}\l( a+b+c
    +\dsum_{i=1}^n (a\circ x^i+b\circ x^i+c\circ x^i)\r)
    \medskip\\
    ~~~~+\dfr{4}{n+7} \l(\dsum_{i=1}^n x^i +\dsum_{1\leq j<k\leq n} x^j\circ x^k\r) 
  \end{array}
  \end{equation}
  and its central charge is equal to $(n+2)(5n+29)/5(n+7)$.
  The commutant $\com_{X^{[n]}} X^{[n-1]}$ has the $c=c_{n+4}$ conformal vector 
  \begin{equation}\label{eq:4.14}
    f^n=\w^n-\w^{n-1}.
  \end{equation}
  Therefore, $X^{[n]}$ has a unitary Virasoro frame 
  $L(c_3,0)\tensor L(c_4,0)\tensor \cds \tensor L(c_{n+4},0)$.
\end{thm}

\pf
Set $c=\tau_a b=\tau_b a$ and 
\begin{equation}\label{eq:4.15}
  A=\{ u_{a,b},a,~ b,~ c,~ x^i,~ a\circ x^i,~ b\circ x^i,~ c\circ x^i,~ 
  x^j\circ x^k~ \mid 1\leq i\leq n,~~ 1\leq j<k\leq n\}.
\end{equation}
It follows from Lemma \ref{lem:3.2} and (1) of Proposition \ref{prop:4.3} that 
both $\la a,b\circ x^i\ra$ and $\la a\circ x^i,b\circ x^j\ra$ are isomorphic to 
the 6A-algebra.
Therefore, we have the following subalgebras in the Griess algebra:
\[
\begin{array}{l}
  \la a,b\ra_2 =\Span\{ u_{a,b},a,b,c\},
  \medskip\\
  \la a,a\circ x^i\ra_2= \la a,x^i\ra_2=\Span\{ a,a\circ x^i,x^i\}, 
  \medskip\\
  \la a,b\circ x^i\ra_2 
  =\Span\{ u_{a,b},x^i,a,b,c,a\circ x^i,b\circ x^i, c\circ x^i\}, 
  \medskip\\
  \la a\circ x^i,b\circ x^i\ra_2 =\sigma_{x^i}\la a,b\ra_2
  = \Span\{ u_{a,b},a\circ x^i,b\circ x^i,c\circ x^i\},
  \medskip\\
  \la a\circ x^i, b\circ x^j\ra_2
  =\Span\{ u_{a,b},x^i\circ x^j,a\circ x^i,b\circ x^i,c\circ x^i,a\circ x^j,
    b\circ x^j,c\circ x^j\} .
\end{array}
\]
In particular, the following orthogonality relation holds.
\[
  (u_{a,b} \mymid x^i)=(u_{a,b}\mymid x^i\circ x^j)=0.
\]
Applying Proposition \ref{prop:4.15} to the subalgebra $\la a,x^1,\dots,x^n\ra$, we also 
have the following orthogonality.
\[
  (a\mymid x^j\circ x^k) =(a\circ x^i \mymid x^j)=(x^i\mymid x^j\circ x^k)
  =(a\circ x^i \mymid x^j\circ x^k)=0.
\]
On $X^{[n]}$, all $\tau_{x^i}$ and $\tau_{x^i\circ x^j}$ are trivial and 
we have $\tau_{e\circ x^i}=\tau_e$ for $e\in \{ a,b,c\}$ by (5) of Theorem \ref{thm:2.9}.
Then $\la \tau_a,\tau_b\ra\cong \mathrm{S}_3$  acts faithfully on each of the 3-sets 
$\{ a\circ x^i,b\circ x^i,c\circ x^i\}, 1\leq i\leq n$.
The group generated by $\sigma_{x^i}$,  $1\leq i\leq n$, is isomorphic to 
$\mathrm{S}_{n+1}$ and this group  acts transitively on the set of Ising vectors 
$\{ x^i, x^j\circ x^k \mid 1\leq i\leq n, 1\leq j<k\leq n\}$ (cf.~Remark \ref{rem:4.16}).
Set $G=\la \tau_a,\tau_b,\sigma_{x^1},\dots,\sigma_{x^n}\ra \subset \aut(X^{[n]})$.
Then $A$ is $G$-invariant and the products and inner products of vectors in $A$ are 
uniquely determined in the linear span of $A$.
Therefore, the linear span of $A$ forms a subalgebra in the Griess algebra.
That $A$ is linearly independent will be shown in Appendix \ref{sec:a.3} using an 
explicit construction.

Next we prove that $\w^n$ is the conformal vector of $X^{[n]}$.
It is clear that $\w^n$ is fixed by $G$.
By using $G$-invariance, one can directly verify that 
$\w^n$ satisfies $(\w^n \mymid t)=(t \mymid t)$ for all $t\in A$. 
Therefore, $\w^n$ is the conformal vector of $X^{[n]}$ by Lemma \ref{lem:4.4}.
The central charge of $\w^n$ is $2(\w^n|\w^n)=(n+2)(5n+29)/5(n+7)$.
Since $X^{[n-1]}$ is a subalgebra of $X^{[n]}$, both $\w^{n-1}$ and 
$f^n=\w^{n}-\w^{n-1}$ are mutually orthogonal Virasoro vectors of $X^{[n]}$ 
and $f^n$ is the conformal vector of 
$\com_{X^{[n]}} X^{[n-1]}$ by Theorem 5.1 of \cite{FZ}.
The central charge of $f^n$ is 
\[
  2(\w^n\mymid \w^n)-2(\w^{n-1}\mymid \w^{n-1})
  = \dfr{n^2+13n+36}{(n+6)(n+7)}=c_{n+4}.
\]
Set 
\begin{equation}\label{eq:4.16}
  v_{a,b}=-\dfr{5}{14}u_{a,b}+\dfr{16}{21}(a+b+c).
\end{equation}
Then $v_{a,b}$ is a $c=c_4=6/7$ Virasoro vector and the conformal vector 
$\w^0$ of $X^{[0]}=\la a,b\ra$ is an orthogonal sum $\w^0=u_{a,b}\dotplus v_{a,b}$ 
by (6) of Theorem \ref{thm:2.11}.
Therefore, we have the following orthogonal decompositions:
\begin{equation}\label{eq:4.17}
  \w^n = \w^0\dotplus f^1 \dotplus \cds \dotplus f^n
  = u_{a,b} \dotplus v_{a,b} \dotplus f^1 \dotplus \cds \dotplus f^n.
\end{equation}
This shows $X^{[n]}$ has a full subVOA isomorphic to 
$L(c_3,0)\tensor L(c_4,0)\tensor \cds \tensor L(c_{n+4},0)$.
This completes the proof.
\qed
\medskip

Set $D^{[0]}=\{ \tau_x \in \aut(V) \mid x\in I_{a,b}\}$ and 
inductively we define 
\begin{equation}\label{eq:4.18}
  D^{[i]}
  :=\{ \tau_y \in D^{[i-1]} \mid \tau_y\tau_{x^i}=\tau_{x^i}\tau_y \} 
\end{equation}

It is clear that $\la D^{[i]}\ra$ is a subgroup of the centralizer of 
$\la \tau_a,\tau_b,\tau_{x^1},\cdots,\tau_{x^n}\ra$ in $\aut(V)$.
For inductive arguments in the next section, we will consider the action of 
$\la D^{[i]}\ra$ restricted to a smaller subalgebra.

\begin{lem}\label{lem:4.18}
  Each involution in $D^{[n]}$  acts trivially on $X^{[n]}$.
\end{lem}

\pf
We prove this by induction on $n$.
By definition, both $\la a,x\ra$ and $\la b,x\ra$ are  2A-algebras for each 
$x\in D^{[0]}=I_{a,b}$, and it follows that $\tau_x$ fixes both $a$ and $b$.
Thus $D^{[0]}$ acts trivially on $X^{[0]}=\la a,b\ra$.
Suppose each involution of $D^{[i]}$ acts on $X^{[i]}=\la a,b,x^1,\cds,x^i\ra$ trivially.
Let $\tau_y\in D^{[i+1]}$.
Then $\tau_y \in D^{[i]}$ and $\tau_y$ acts trivially on $X^{[i]}$.
Since $\tau_y$ and $\tau_{x^{i+1}}$ commute, it follows from Theorem \ref{thm:2.6} 
that $\la y,x^{i+1}\ra$ is a dihedral algebra of type 1A, 2A or 2B.
In each case $\tau_{y}$ fixes $x^{i+1}$.
Therefore, $\tau_y$ acts trivially on $X^{[i+1]}=\la X^{[i]},x^{i+1}\ra$.
\qed
\medskip

By this lemma, we have a group homomorphism by restriction.
\begin{equation}\label{eq:4.19}
\begin{array}{cccc}
  \varphi^{[i]}: &\la D^{[i]}\ra &\longrightarrow &\aut(\com_VX^{[i]})
  \medskip\\
  & g & \longmapsto & g|_{\com_VX^{[i]}}
\end{array}
\end{equation}
Set $G^{[i]}:=\mathrm{Im}\, \varphi^{[i]}\subset \aut(\com_VX^{[i]})$.
It follows from Corollary \ref{cor:3.7} that $\l( G^{[i]},\varphi^{[i]}(D^{[i]})\r)$ 
is a 3-transposition group.
Since $X^{[n-1]}\subset X^{[n]}$ implies $\com_V X^{[n-1]}\supset \com_V X^{[n]}$, 
we obtain the following inductive structure:
\begin{equation}\label{eq:4.20}
\begin{array}{ccccccccc}
  ~~~X^{[0]} &\subset& ~~~X^{[1]} &\subset& ~~~X^{[2]} &\subset& ~~~X^{[3]} 
  &\subset& \cds
  \vspace{-10pt}\\
  \mbox{\large \rotatebox{270}{$\leadsto$}} 
  && \mbox{\large \rotatebox{270}{$\leadsto$}} 
  && \mbox{\large \rotatebox{270}{$\leadsto$}} 
  && \mbox{\large \rotatebox{270}{$\leadsto$}} 
  && 
  \\
  \com_V X^{[0]} &\supset& \com_V X^{[1]} &\supset& \com_V X^{[2]} &\supset& 
  \com_V X^{[3]} &\supset& \cds
  \\ 
  \circlearrowleft && 
  \circlearrowleft && 
  \circlearrowleft && 
  \circlearrowleft && 
  \smallskip\\
  G^{[0]} &&  G^{[1]} && G^{[2]} && G^{[3]} && \cds
  \\
  \mbox{\large \rotatebox{90}{$\twoheadrightarrow$}} 
  && \mbox{\large \rotatebox{90}{$\twoheadrightarrow$}} 
  && \mbox{\large \rotatebox{90}{$\twoheadrightarrow$}}
  && \mbox{\large \rotatebox{90}{$\twoheadrightarrow$}}
  && 
  \\ 
  \la D^{[0]}\ra & \supset & \la D^{[1]}\ra & \supset & \la D^{[2]}\ra & \supset &
  \la D^{[3]}\ra & \supset & \cds 
\end{array}
\end{equation}
Note that $\varphi^{[n-1]}(\tau_{x^n})\in G^{[n-1]}$ is an external automorphism of 
$\com_V X^{[n-1]}$ in the sense that $x^n \not\in \com_V X^{[n-1]}$.
By Theorem 5.2 of \cite{FZ} and \eqref{eq:4.17} we have 
\[
  \com_V X^{[n]}=\com_V\la u_{a,b},v_{a,b},f^1,\dots,f^n\ra .
\]
Since the decomposition in \eqref{eq:4.17} is orthogonal, we have 
$f^n \in \com_V X^{[n-1]}$.
We prove that $\varphi^{[n-1]}(\tau_{x^n})$ can be described as an internal 
automorphism defined by $f^n$.

\begin{thm}\label{thm:4.19}
  Let $X^{[n]}=\la a,b,x^1,\dots,x^n\ra$ be the subalgebra of a VOA $V$ defined 
  by \eqref{eq:4.12} and let 
  $\w^n=u_{a,b}\dotplus v_{a,b}\dotplus f^1\dotplus \cds \dotplus f^n$ 
  be the orthogonal decomposition of the conformal vector of $X^{[n]}$ as 
  in \eqref{eq:4.17}.
  \\
  (1) If $n$ is odd, then $\tau_{f^n}=\tau_{x^1}\tau_{x^2}\cds \tau_{x^n}$ 
  in $\aut(V)$ and $\tau_{f^n}=\varphi^{[n-1]}(\tau_{x^n})$ as an automorphism 
  of $\com_V X^{[n-1]}$.
  \\
  (2) If $n>0$ is even, then 
  $\tau_{f^n}=\tau_{f^{n-1}}=\tau_{x^1}\tau_{x^2}\cds \tau_{x^{n-1}}$ in $\aut(V)$ and 
  $\varphi^{[n-1]}(\tau_{f^n})$ is trivial in $\aut(\com_V X^{[n-1]})$.
  In this case, $f^n$ is of $\sigma$-type on $\com_V X^{[n-1]}$ and satisfies 
  $\sigma_{f^n}=\varphi^{[n-1]}(\tau_{x^n})$ as an automorphism of $\com_V X^{[n-1]}$.
\end{thm}

\pf
By direct computations, we have the following relations in $X^{[n]}$.
\begin{equation}\label{eq:4.21}
\begin{array}{l}
  a_{(1)}u_{a,b}=\dfr{2}{3}a+\dfr{5}{24}u_{a,b}-\dfr{7}{24}v_{a,b},
  ~~~ a_{(1)}v_{a,b}=\dfr{4}{3}a-\dfr{5}{24}u_{a,b}+\dfr{7}{24}v_{a,b},
  \medskip\\
  x^1_{(1)}v_{a,b}=\dfr{5}{7}x^1+\dfr{3}{14}v_{a,b}-\dfr{2}{7}f^1,
  ~~~ x^1_{(1)}f^1=\dfr{9}{7}x^1-\dfr{3}{14}v_{a,b}+\dfr{2}{7}f^1,
  \medskip\\
  (x^i\circ x^{i+1})_{(1)}f^{i}
  = \dfr{i+5}{i+7}\, x^i\circ x^{i+1}+\dfr{i+6}{4(i+7)}f^i-\dfr{i+8}{4(i+7)}f^{i+1},
  \medskip\\
  (x^i\circ x^{i+1})_{(1)}f^{i+1}
  = \dfr{i+9}{i+7}\, x^i\circ x^{i+1}-\dfr{i+6}{4(i+7)}f^i+\dfr{i+8}{4(i+7)}f^{i+1},
\end{array}
\end{equation}
where $1\leq i\leq n-1$.
By these relations, it follows from Proposition 4.2 of \cite{LY3} that 
the subalgebras $\la a,u_{a,b},v_{a,b}\ra$, $\la x^1,v_{a,b},f^1\ra$ 
and $\la x^i\circ x^{i+1},f^i,f^{i+1}\ra$ are uniquely determined 
and isomorphic to $A(\shf,c_3^1)$, $A(\shf,c_4^1)$ and $A(\shf,c_{i+4}^1)$ 
in (loc.~cit.), respectively.
Then by Theorem 4.6 of \cite{LY3}, we have the following relations 
in $\aut(V)$.
\[
\begin{array}{l}
  \tau_{u_{a,b}}=\tau_{v_{a,b}},~~~
  \tau_{v_{a,b}}\tau_{f^1}=\tau_{x^1},~~~
  \tau_{f^{i}}\tau_{f^{i+1}}
  =\begin{cases}
    \tau_{x^i\circ x^{i+1}}=\tau_{x^i}\tau_{x^{i+1}} & \mbox{if $i$ is even},
    \\
    1 & \mbox{if $i$ is odd}.
  \end{cases}
\end{array}
\]
It is also shown in Remark 2.12 of \cite{HLY2} that 
$\tau_{u_{a,b}}=\tau_{v_{a,b}}$ is trivial on $V$.
Therefore, $\tau_{f^1}=\tau_{x^1}$ and by induction we have
\[
  \tau_{f^n}=\begin{cases}
    \tau_{x^1}\tau_{x^2}\cds \tau_{x^n} & \mbox{if $n$ is odd,}
    \\
    \tau_{x^1}\tau_{x^2}\cds \tau_{x^{n-1}} & \mbox{if $n>0$ is even.}
    \end{cases}
\]
It is clear that $\tau_{x^i}$ are trivial on the commutant $\com_V X^{[n-1]}$ 
for $1\leq i\leq n-1$ since $x^1,\dots,x^{n-1} \in X^{[n-1]}$.
Thus as an automorphism of $\com_V X^{[n-1]}$ we have 
\[
  \tau_{f^n}
  =\begin{cases}
   \varphi^{[n-1]}(\tau_{x^1}\cds \tau_{x^n})
   =\varphi^{[n-1]}(\tau_{x^n}) & \mbox{if $n$ is odd,}
   \\
   \varphi^{[n-1]}(\tau_{x^1}\cds \tau_{x^{n-1}})=1 & \mbox{if $n>0$ is even.}
  \end{cases}
\]
If $n>0$ is even, it is shown in Theorem 4.8 of \cite{LY3} that 
$f^n$ is of $\sigma$-type on the commutant $\com_V \la f^{n-1}\ra$ 
and both $\sigma_{f^n}$ and $\tau_{x^{n-1}\circ x^n}$ define the same automorphism 
on $\com_V \la f^{n-1}\ra$.
Since $\com_V X^{[n-1]}\subset \com_V\la f^{n-1}\ra$, we have
\[
  \sigma_{f^n}
  =\varphi^{[n-1]}(\tau_{x^{n-1}\circ x^n})
  =\varphi^{[n-1]}(\tau_{x^{n-1}}\tau_{x^n})
  =\varphi^{[n-1]}(\tau_{x^n})
\]
as an automorphism of $\com_V X^{[n-1]}$.
This completes the proof.
\qed

\section{Conway-Miyamoto correspondence}

In this section, we apply our theorems to the moonshine vertex operator algebra 
and establish the Conway-Miyamoto correspondence for the Fischer 
3-transposition groups.
\medskip

Throughout this section, we identify the Monster simple group $\M$ with 
the automorphism group of the moonshine VOA $V^\natural$ (cf.~\cite{FLM}).
The moonshine VOA $V^\natural$ can be defined over the real numbers and has 
a compact real form $V^\natural_\R$ as in \cite{FLM,M4}.
It is shown in Proposition 5.2 of \cite{HLY1} that all Ising vectors in $V^\natural$ 
are contained in the compact real form $V^\natural_\R$.
Therefore, $V^\natural$ has a canonical real form defined as the real subalgebra 
generated by Ising vectors, which indeed coincides with $V^\natural_\R$. 
In the following discussion, we will consider subalgebras generated by 
Ising vectors of $V^\natural$ so that all subalgebras have compact real forms.

\begin{df}\label{df:5.1}
  Let $V$ be a VOA and $G$ a subgroup of $\aut(V)$.
  Let $I$ be a conjugacy class of involutions of $G$.
  We define the \emph{Conway-Miyamoto correspondence between involutions $I$ of $G$ and   
  $c=c_n$ Virasoro vectors of $V$} by means of the following conditions:
  \vspace{-5pt}
  \begin{enumerate}
  \setlength{\parskip}{0pt}
  \setlength{\itemsep}{2pt}
  \setlength{\itemindent}{-5pt}
  \item 
  For each $t\in I$, there exists a unique $c=c_n$ Virasoro vector $e_t \in V^{C_G(t)}$.
  \item 
  If the unique Virasoro vector $e_t$ is not of $\sigma$-type on $V$, then $\tau_{e_t}=t$ on $V$. 
  Otherwise, $\sigma_{e_t}=t$ on $V$.
  \end{enumerate}  
  \vspace{-7pt}
  The unique $c=c_n$ Virasoro vector $e_t$ of $V^{C_G(t)}$ is called the 
  \emph{axial vector} associated to $t$.
  We say that the Conway-Miyamoto correspondence between $G$ and $V$ is 
  \emph{bijective} if the axial vector $e_t$ is the unique $c=c_n$ Virasoro 
  vector $a$ of $V$ satisfying $\tau_a=t$ (or $\sigma_a=t$ in the case $e_t$ is 
  of $\sigma$-type on $V$).
\end{df}

Here are known examples of Conway-Miyamoto correspondences.

\begin{thm}[\cite{ATLAS,FLM,C,M1,Ma1,Ho,HLY1,HLY2}]\label{thm:5.2}~\\
  There are Conway-Miyamoto correspondences between the following groups and VOAs.
  \\
  (1) 2A-elements of the Monster and Ising vectors of $V^\natural$. 
  The correspondence is bijective.
  \\
  (2) 2A-elements of the Baby Monster and $c=c_2$ Virasoro vectors of 
  $\sigma$-type of $V\!\!B^\natural=\com_{V^\natural}\la e\ra$,
  where $e$ is an Ising vector of $V^\natural$. 
  The correspondence is bijective.
  \\
  (3) 2C-elements of the largest Fischer 3-transposition group and 
  $c=c_4$ Virasoro vectors of $\sigma$-type of 
  $V\!\!F^\natural=\com_{V^\natural}\la u_{a,b}\ra$, where $a$, $b$ are Ising vectors 
  of $V^\natural$ such that $(a \mymid b)=13\cd 2^{-10}$.
\end{thm}

By the bijective correspondence in (1) of Theorem \ref{thm:5.2} and 
Proposition \ref{prop:2.5}, we have

\begin{cor}\label{cor:5.3}
  Let $e^1,\dots,e^k$ be Ising vectors of $V^\natural$.
  Then the centralizer $C_\M(\tau_{e^1},\cds,\tau_{e^k})$ coincides with 
  the pointwise stabilizer of $\la e^1,\dots,e^k\ra$ in $\M$.
\end{cor}

Now we fix a pair of Ising vectors $a$, $b$ in $V^\natural$ such that
$\la a,b\ra$ is isomorphic to the 3A-algebra.
Such a pair corresponds to a pair of 2A-involutions $\tau_a$ and $\tau_b$ 
such that $\tau_a\tau_b$ is a 3A-element by (1) of Theorem \ref{thm:5.2}, 
and is known to be unique up to conjugation by $\M$ (cf.~\cite{ATLAS}).

Let $E_{V^\natural}$ be the set of Ising vectors of $V^\natural$ and 
$I_{a,b}=\{ x\in E_{V^\natural} \mid (a\mymid x)=(b\mymid x)=2^{-5}\}$ as before.
For $n=0,1,2$,  we consider the subalgebra $X^{[n]}=\la a,b,x^1,\cds,x^n\ra$ of 
$V^\natural$ defined as in Sec.~\ref{sec:4.3} and study the automorphism group 
$G^{[n]}=\la \varphi^{[i]}(D^{[i]})\ra$ of $\com_{V^\natural} X^{[n]}$ defined as in 
\eqref{eq:4.18} and \eqref{eq:4.19}.

\begin{prop}\label{prop:5.4}
  Let $a$, $b$, $x^1$, $x^2\in V^\natural$ be as above.
  \\
  (1) $\la D^{[0]}\ra =C_\M(\tau_a,\tau_b)\cong \Fi_{23}$ in $\aut(V^\natural)$ 
  and $G^{[0]}\cong \Fi_{23}$ in $\aut(\com_{V^\natural} X^{[0]})$.
  Moreover, for any $\tau_y\in D^{[0]}$, $\varphi^{[0]}(\tau_y)$ defines 
  a 2A-element of $\Fi_{23}$ on $\aut(\com_{V^\natural} X^{[0]})$.
  \\
  (2) $\la D^{[1]}\ra =C_\M(\tau_a,\tau_b,\tau_{x^1}) \cong 2.\Fi_{22}$ in 
  $\aut(V^\natural)$ and $G^{[1]}\cong \Fi_{22}$ in $\aut(\com_{V^\natural} X^{[1]})$.
  Moreover, for any $\tau_y\in D^{[1]}$, 
  $\varphi^{[1]}(\tau_y)$ defines a 2A-element of $\Fi_{22}$ 
  on $\aut(\com_{V^\natural} X^{[1]})$.
  \\
  (3) $\la D^{[2]}\ra =C_\M(\tau_a,\tau_b,\tau_{x^1},\tau_{x^2})\cong 
  2^2.\mathrm{PSU}_6(2)$ and $G^{[2]}\cong \mathrm{PSU}_6(2)$ in 
  $\aut(\com_{V^\natural} X^{[2]})$.
  Moreover, for any $\tau_y\in D^{[2]}$, 
  $\varphi^{[2]}(\tau_y)$ defines a 2A-element of $\mathrm{PSU}_6(2)$ 
  on $\aut(\com_{V^\natural} X^{[2]})$.
\end{prop}

\pf
First, we observe that $\la D^{[n]}\ra$ is a normal subgroup of 
$C_\M(\tau_a,\tau_b,\tau_{x^1},\dots,\tau_{x^n})$ by Proposition \ref{prop:2.5}.
\\
(1):~
Since $\M$ acts transitively on pairs of Ising vectors of $V^\natural$ 
generating 6A-subalgebras \cite{ATLAS}, 
there exists a 6A-subalgebra of $V^\natural$ containing $\la a,b\ra$ \cite{LM} 
and its central Ising vector belongs to $I_{a,b}$ so that $D^{[0]}$ is non-trivial.
Since $\la D^{[0]}\ra$ is a normal subgroup of $C_\M(\tau_a,\tau_b)$
and $C_\M(\tau_a,\tau_b)\cong \Fi_{23}$ is a simple group by \cite{ATLAS},
we have $\la D^{[0]}\ra =C_\M(\tau_a,\tau_b)\cong \Fi_{23}$.
By definition, $G^{[0]}$ is a homomorphic image of $\la D^{[0]}\ra$ and we have 
$G^{[0]}\cong \Fi_{23}$ in $\aut(\com_{V^\natural} X^{[0]})$.
By Corollary \ref{cor:3.7}, $D^{[0]}$ is the set of 3-transpositions, and 
only the 2A-involutions satisfy the 3-transposition property in $\Fi_{23}$ 
(cf.~\cite{ATLAS}, \cite[Page 45]{I}).
Therefore, each involution in $D^{[0]}$ defines a 2A-element of $\Fi_{23}$. 
\\
(2):~By (1) and Corollary \ref{cor:5.3} we have 
$C_\M(\tau_a,\tau_b,\tau_{x^1})\cong C_{\Fi_{23}}(\mathrm{2A})\cong 2.\Fi_{22}$ \cite{ATLAS}. 
Since $\la D^{[1]}\ra$ is normal in $C_\M(\tau_a,\tau_b,\tau_{x^1})$, 
we have $\la D^{[1]}\ra =C_\M(\tau_a,\tau_b,\tau_{x^1}) \cong 2.\Fi_{22}$ 
in $\aut(V^\natural)$ where the center is generated by $\tau_{x^1}$. 
Then $G^{[1]} \cong \Fi_{22}$ since the center $\tau_{x^1}$ acts trivially 
on $\aut(\com_{V^\natural} X^{[1]})$ by Lemma \ref{lem:4.18}.
Again it follows from Corollary \ref{cor:3.7} that $\la \varphi^{[1]}(D^{[1]})\ra$ 
is the set of 3-transpositions in $G^{[1]}\cong \Fi_{22}$ and only 2A-involutions  
satisfy the 3-transposition property in $\Fi_{22}$ (cf.~\cite{ATLAS}, \cite[Page 45]{I}).
Therefore $\varphi^{[1]}(\tau_y)$ defines a 2A-element of $\Fi_{22}$ on 
$\com_{V^\natural} X^{[1]}$ for any $\tau_y\in D^{[1]}$.
\\
(3):~ By (2) and Corollary \ref{cor:5.3},  we have 
$\la D^{[2]}\ra=C_\M(\tau_a,\tau_b,\tau_{x^1},\tau_{x^2})\cong 2^2.\mathrm{PSU}_6(2)$ \cite{ATLAS}.
Since the center $\la \tau_{x^1},\tau_{x^2}\ra$ acts trivially on 
$\com_{V^\natural} X^{[2]}$, we have $G^{[2]}\cong \mathrm{PSU}_6(2)$ in 
$\aut(\com_{V^\natural} X^{[2]})$.
Since only 2A-involutions satisfy the 3-transposition property in 
$\mathrm{PSU}_6(2)$ (cf.~\cite{ATLAS}, \cite[Page 45]{I}), 
$\varphi^{[2]}(\tau_y)$ is a 2A-element for any $\tau_y\in D^{[2]}$.
This completes the proof.
\qed

\begin{lem}\label{lem:5.5}
  Let $\w^0$, $\w^1$ and $\w^2$ be the conformal vectors of the subalgebras 
  $X^{[0]}=\la a,b\ra$, $X^{[1]}=\la a,b,x^1\ra$ and $X^{[2]}=\la a,b,x^1,x^2\ra$ 
  of $V^\natural$, respectively, and let $f^1=\w^1-\w^0$ and $f^2=\w^2-\w^1$ 
  be the conformal vectors of $\com_{X^{[1]}}X^{[0]}$ and $\com_{X^{[2]}} X^{[1]}$,  
  respectively.
  Then $f^1$ is fixed by $C_{G^{[0]}}(\varphi^{[0]}(\tau_{x^1}))$ and 
  $f^2$ is fixed by $C_{G^{[1]}}(\varphi^{[1]}(\tau_{x^2}))$.
\end{lem}

\pf
It follows from Proposition \ref{prop:5.4} that 
$C_{G^{[0]}}(\varphi^{[0]}(\tau_{x^1}))=\varphi^{[0]}(C_\M(\tau_a,\tau_b,\tau_{x^1}))$ and 
$C_{G^{[1]}}(\varphi^{[1]}(\tau_{x^2}))
=\varphi^{[1]}(C_\M(\tau_a,\tau_b,\tau_{x^1},\tau_{x^2}))$.
Then $C_\M(\tau_a,\tau_b,\tau_{x^1})$ fixes $f^1=\w^0-\w^1\in \la a,b,x^1\ra$ and 
$C_\M(\tau_a,\tau_b,\tau_{x^1},\tau_{x^2})$ fixes 
$f^2=\w^1-\w^2\in \la a,b,x^1,x^2\ra$ by Corollary \ref{cor:5.3}.
Thus the claim follows.
\qed

\subsection{Transpositions of $\Fi_{23}$ and $c=25/28$ Virasoro vectors}

Let $a$, $b$ be Ising vectors of $V^\natural$ such that $\la a,b\ra$ is 
the 3A-algebra.
In this subsection, we consider $X^{[0]}=\la a,b\ra\subset V^\natural$ and 
its commutant subalgebra $\VF_{23}= \com_{V^\natural}\la a,b\ra$.
Let   
\[
  D^{[0]}=\{ \tau_x \in \aut(V^\natural) \mid x\in I_{a,b}\}, ~~~
  G^{[0]}=\la \varphi^{[0]}(D^{[0]})\ra \subset \aut(\com_{V^\natural}\la a,b\ra)
\]
be defined as in \eqref{eq:4.18} and \eqref{eq:4.19}.
We have shown in Proposition \ref{prop:5.4} that $G^{[0]}\cong \Fi_{23}$.
In the following, we will identify $\Fi_{23}$ with 
$G^{[0]}\subset \aut(\com_{V^\natural}\la a,b\ra)$.
Let $x^1\in D^{[0]}$ and $X^{[1]}=\la a,b,x^1\ra$.
Let $f^1$ be the conformal vector of $\com_{X^{[1]}} X^{[0]}$.
Then it follows from Theorem \ref{thm:4.19} and Proposition \ref{prop:5.4} that 
$\tau_{f^1}=\varphi^{[0]}(\tau_{x^1})$ defines a 2A-element of $\Fi_{23}$.
Moreover, it is shown in Lemma \ref{lem:5.5} that $f^1$ is fixed by 
the centralizer $C_{\Fi_{23}}(\tau_{f^1})$.
We will show that $f^1$ is the unique $C_{\Fi_{23}}(\tau_{f^1})$-invariant 
$c=c_5$ Virasoro vector of $\com_{V^\natural}\la a,b\ra$.

Recall that the 3A-algebra $U_{\mathrm{3A}}=\la a,b\ra$ has 6 irreducible representations 
$U(h)$ which are distinguished by the top weights 
$h\in \{ 0,1/7,5/7,2/5,19/35,4/35\}$ (cf.~\cite{SY,LLY}).
We denote by $\irr\, U_{\mathrm{3A}}$ the set of equivalent classes of irreducible 
$U_{\mathrm{3A}}$-modules.
For $U(h)\in \irr\, U_{\mathrm{3A}}$,  we set
\[
  U(h)^c:=\hom_{\la a,b\ra}(U(h),V^\natural).
\]
Then we have the isotypical decomposition
\begin{equation}\label{eq:5.1}
  V^\natural=\bigoplus_{U(h)\in \, \irr\,\la a,b\ra} U(h)\tensor U(h)^c.
\end{equation}
The commutant $\com_{V^\natural}\la a,b\ra$  acts naturally on $U(h)^c$ and
one can consider \eqref{eq:5.1} as the decomposition  as a
$\la a,b\ra\tensor \com_{V^\natural}\la a,b\ra$-module.
We denote the top levels of $U(h)$ and $U(h)^c$ by $\mathrm{Top}~U(h)$ and 
$\mathrm{Top}~U(h)^c$, respectively.
The top weight and the dimension of the top level of $U(h)^c$ are obtained 
in Lemma 5.6 of \cite{HLY2}.

\renewcommand{\arraystretch}{1.5}

\begin{lem}\label{lem:5.6}
  The top levels of $U(h)^c$ are irreducible as $G^{[0]}\cong \Fi_{23}$-modules and 
  their characters are as in the following table.
  \[
  \begin{array}{|c||c|c|c|c|c|c|}\hline
    U(h) & U(0)  & U(\sfr{5}{7}) & U(\sfr{19}{35}) & U(\sfr{2}{5})
    & U(\sfr{1}{7}) & U(\sfr{4}{35})
    \\ \hline
    \dim \mathrm{Top}~U(h) & 1 & 3 & 3 & 1 & 1 & 2
    \\ \hline
    \mbox{\rm Top weight of}~U(h)^c & 0  & 9/7 & 51/35 & 8/5 & 13/7 & 66/35
    \\ \hline
    \dim \mathrm{Top}~U(h)^c & 1 & 782 & 3588 & 5083 & 25806 & 60996
    \\ \hline
    \Fi_{23}\mbox{\rm -character} & \chi_1  & \chi_2 & \chi_3 & \chi_4 & \chi_5 & \chi_7
    \\ \hline
  \end{array}
  \]
  Here $\chi_i$ denotes the irreducible characters of $Fi_{23}$ labeled as  
  in \cite[page 178]{ATLAS}.
  By abuse of notations, we also use $\chi_i$ to denote the irreducible module 
  affording the character $\chi_i$.
  Then we have the decomposition $\l(\com_{V^\natural}\la a,b\ra\r)_2=\chi_1+\chi_6\,
  (={\bf 1}+{\bf 30888})$ as a $\Fi_{23}$-module.
\end{lem}

\renewcommand{\arraystretch}{1}

The centralizer $C_{\Fi_{23}}(\tau_{f^1})$ is isomorphic to $2.\Fi_{22}$ 
by Proposition \ref{prop:5.4} and we have the following decomposition 
(cf.~\cite{ATLAS}).

\renewcommand{\arraystretch}{0.6}
\begin{prop}\label{prop:5.7}
  Consider the action of the centralizer 
  $C_{G^{[0]}}(\varphi^{[0]}(\tau_{x^1}))\cong 2.\Fi_{22}$ on the Griess algebra of 
  $\com_{V^\natural}\la a,b\ra$ and the top levels $\mathrm{Top}\,U(h)^c$ of $U(h)^c$ 
  in \eqref{eq:5.1}.
  As $2.\Fi_{22}$-modules,  we have the following decompositions.
  \[
  \begin{array}{l}
    \begin{array}{cccccccccccc}
      \l(\com_{V^\natural}\la a,b\ra\r)_2 =& \mathbf{1} &+& \mathbf{1} &+& \mathbf{429} 
      &+& \mathbf{3080} &+& \mathbf{13650} &+& \mathbf{13728},
      \\
      & \chi_1 && \chi_1 && \chi_3 && \chi_7 && \chi_9 && \chi_{73}
    \end{array}
    \bigskip \\
    \begin{array}{cccccc}
      \mathrm{Top}~U(\sfr{5}{7})^c=& \mathbf{1} &+& \mathbf{429} &+& \mathbf{352},
      \\
      & \chi_1 && \chi_3 && \chi_{66}
    \end{array}
    \bigskip \\
    \begin{array}{cccccc}
      \mathrm{Top}~U(\sfr{19}{35})^c=& \mathbf{78} &+& \mathbf{1430} &+& \mathbf{2080},
      \\
      & \chi_2 && \chi_5 && \chi_{68}
    \end{array}
    \bigskip \\
    \begin{array}{cccc}
      \mathrm{Top}~U(\sfr{2}{5})^c=& \mathbf{3003} &+& \mathbf{2080},
      \\
      & \chi_6 && \chi_{67}
    \end{array}
    \bigskip \\
    \begin{array}{cccccccc}
      \mathrm{Top}~U(\sfr{1}{7})^c=& \mathbf{1001} &+& \mathbf{10725} &+& \mathbf{352} 
      &+& \mathbf{13728},
      \\
      & \chi_4 && \chi_8 && \chi_{66} && \chi_{74}
    \end{array}
    \bigskip \\
    \begin{array}{cccccccc}
      \mathrm{Top}~U(\sfr{4}{35})^c=& \mathbf{1430} &+& \mathbf{30030} &+& \mathbf{2080} 
      &+& \mathbf{27456},
      \\
      & \chi_5 && \chi_{10} && \chi_{67} && \chi_{75}
    \end{array}
  \end{array}
  \]
  where $\chi_i$ are the irreducible characters of $\Fi_{22}$ and $2.\Fi_{22}$ 
  labeled as in \cite[pages 156--157]{ATLAS}.
\end{prop}
\renewcommand{\arraystretch}{1}

By Lemma \ref{lem:5.5} and Proposition \ref{prop:5.7}, the Griess algebra of 
the $C_{\Fi_{23}}(\tau_{f^1})$-invariants of $\com_{V^\natural}\la a,b\ra$
is 2-dimensional spanned by the $c=782/35$ conformal vector and 
the $c=25/28$ Virasoro vector $f^1\in \com_{X^{[1]}}X^{[0]}$.
Therefore, $f^1$ is the axial vector of the 2A-element $\tau_{f^1}$ of $\Fi_{23}$ 
and we have established the Conway-Miyamoto correspondence for $\Fi_{23}$.

\begin{thm}\label{thm:5.8}
  Let $E_{V^\natural}$ be the set of Ising vectors of $V^\natural$ and let $a$, 
  $b\in E_{V^\natural}$ be a pair such that  $\la a,b \ra\cong U_{3A}$, i.e.,  $(a\mymid b)=13\cd 2^{-10}$.
  Then there exists a Conway-Miyamoto correspondence between 2A-elements of 
  $\Fi_{23}$ and $c=25/28$ Virasoro vectors of $\com_{V^\natural}\la a,b\ra$.
  More precisely, set 
  $I_{a,b}=\{ x\in E_{V^\natural} \mid \la a,x\ra\cong \la b,x\ra \cong U_{2A} \}$, 
  i.e., $(a\mymid x)=(b\mymid x)=2^{-5}$ for any $x\in I_{a,b}$.   
  Then the following hold.
  \\
  (1)~ $C_\M(\tau_a,\tau_b)=\la \tau_x \mid x\in I_{a,b}\ra\cong \Fi_{23}$ 
  acts faithfully on $\com_{V^\natural}\la a,b\ra$.
  \\
  (2)~ There is a one to one correspondence between $I_{a,b}$ and 2A-involutions of 
  $\Fi_{23}$ via the Miyamoto involution 
  $x\longmapsto \tau_x|_{\com_{V^\natural}\la a,b\ra}$. 
  \\
  (3)~ Given a 2A-involution $t$ of $\Fi_{23}$, there exists a unique 
  $C_{\Fi_{23}}(t)$-invariant $c=25/28$ Virasoro vector $f$ of 
  $\com_{V^\natural}\la a,b\ra$ such that $\tau_f=t$.
  Moreover, $\tau_f$ on $V^\natural$ defines a $2A$-involution of $\M$.
\end{thm}

For the inductive argument from $G^{[0]}=\Fi_{23}$ to $G^{[1]}=\Fi_{22}$,  we determine 
the decomposition of the Griess algebra of $\com_{V^\natural} X^{[1]}$ as a $\Fi_{22}$-module.

\begin{lem}\label{lem:5.9}
  Consider the subalgebras $X^{[0]}=\la a,b\ra \subset X^{[1]}=\la a,b,x^1\ra
  \subset V^\natural$ as above and let $f^1$ be the conformal vector of 
  $\com_{X^{[1]}}X^{[0]}$.
  Then the zero-mode $\o(f^1)$ acts on the Griess algebra of $V^\natural$ semisimply 
  with possible eigenvalues
  \[
     0,~ 2,~ \dfr{9}{7},~ \dfr{3}{4},~ \dfr{5}{14},~ \dfr{1}{28},~ 
     \dfr{3}{28},~ \dfr{15}{28},~ 
     \dfr{5}{32},~ \dfr{3}{224},~ \dfr{15}{224},~ \dfr{99}{224},~ 
     \dfr{143}{224}.
  \]
  More precisely, $\mathrm{o}(f^1)$ acts semisimply on the Griess algebra of 
  $\com_{V^\natural}\la a,b\ra$ and the top levels $\mathrm{Top}\,U(h)^c$ 
  with possible eigenvalues as follows.
  \[
  \begin{array}{ll}
    (\com_{V^\natural}\la a,b\ra)_2 : 0,\sfr{3}{4},\sfr{5}{32} ,
    &\mathrm{Top}\, U(\sfr{5}{7})^c :\sfr{9}{7},\sfr{1}{28},\sfr{15}{224},
    \medskip\\
    \mathrm{Top}\, U(\sfr{2}{5})^c : 0,\sfr{5}{32},
    &\mathrm{Top}\, U(\sfr{1}{7})^c : \sfr{5}{14},\sfr{3}{28},\sfr{143}{224},
    \sfr{3}{224},
    \medskip\\
    \mathrm{Top}\, U(\sfr{19}{35})^c : \sfr{5}{14},\sfr{3}{28},\sfr{3}{224},~~~
    &\mathrm{Top}\, U(\sfr{4}{35})^c: \sfr{15}{28},\sfr{1}{28},\sfr{99}{224},
    \sfr{15}{224}.
  \end{array}
  \]
\end{lem}

\pf
Since the axial vector $f^1$ is fixed by $C_{\Fi_{23}}(\tau_{f^1})$ by 
Theorem \ref{thm:5.8}, its zero-mode $\mathrm{o}(f^1)$ acts as a scalar on 
each $C_{\Fi_{23}}(\tau_{f^1})$-irreducible component.
By \eqref{eq:5.1}, we have the following decomposition of the Griess algebra 
of $V^\natural$.
\begin{equation}\label{eq:5.2}
  V_2^\natural=\la a,b\ra_2 \oplus \l(\com_{V^\natural}\la a,b\ra\r)_2 \oplus
  \bigoplus_{U(h)\in \irr\,\la a,b\ra \atop h>0} \mathrm{Top}~U(h)\tensor 
  \mathrm{Top}~U(h)^c.
\end{equation}
The 6A-algebra $\la a,b\circ x^1\ra = \la a,b,x^1\ra$ contains a 
rational full subVOA $\la a,b,f^1\ra =\la a,b\ra \tensor \la f^1\ra$
by (8) of Theorem \ref{thm:2.13}.
As explained in the beginning of this section, all Ising vectors of $V^\natural$ 
are contained in the compact real form $V^\natural_\R$ of $V^\natural$ so that 
$V^\natural$ is semisimple as a $\la a,b,x^1\ra$-module.
The uniqueness of a simple VOA structure of the 6A-algebra 
as well as the classification of irreducible modules over the 6A-algebra 
are established in \cite{DJY}.
By (loc.~cit.), an irreducible $U_{\mathrm{6A}}$-module is isomorphic to one of the following as 
a $U_{\mathrm{3A}}\tensor L(\sfr{25}{28},0)$-module.
\begin{equation}\label{eq:5.3}
\begin{array}{ll}
  {}[0, 0]
  \oplus [\sfr{1}{7}, \sfr{34}{7}]
  \oplus [\sfr{5}{7}, \sfr{9}{7}],
  & [0, \sfr{5}{32}]
  \oplus [\sfr{1}{7}, \sfr{675}{224}]
  \oplus [\sfr{5}{7}\tensor [\sfr{99}{224}],
  \medskip\\
  {}[\sfr{1}{7}, \sfr{3}{224}]
  \oplus [0, \sfr{165}{32}]
  \oplus [\sfr{5}{7}, \sfr{323}{224}],
  & [\sfr{19}{35}, \sfr{3}{224}]
  \oplus [\sfr{2}{5}, \sfr{165}{32}]
  \oplus [\sfr{4}{35}, \sfr{323}{224}],
  \medskip\\
  {}[\sfr{1}{7}, \sfr{5}{14}]
  \oplus [0, \sfr{15}{2}]
  \oplus [\sfr{5}{7}, \sfr{39}{14}],
  & [\sfr{19}{35}, \sfr{5}{14}]
  \oplus [\sfr{2}{5}, \sfr{15}{2}]
  \oplus [\sfr{4}{35}, \sfr{39}{14}],
  \medskip\\
  {}[\sfr{2}{5}, 0]
  \oplus [\sfr{19}{35}, \sfr{34}{7}]
  \oplus [\sfr{4}{35}, \sfr{9}{7}],
  & [\sfr{4}{35}, \sfr{15}{224}]
  \oplus [\sfr{2}{5}, \sfr{57}{32}]
  \oplus [\sfr{19}{35}, \sfr{143}{224}],
  \medskip\\
  {}[\sfr{4}{35}, \sfr{1}{28}]
  \oplus [\sfr{2}{5}, \sfr{3}{4}]
  \oplus [\sfr{19}{35}, \sfr{45}{28}],
  & [\sfr{1}{7}, \sfr{3}{28}]
  \oplus [0, \sfr{13}{4}]
  \oplus [\sfr{5}{7},  \sfr{15}{28}],
  \medskip\\
  {}[0, \sfr{3}{4}]
  \oplus [\sfr{5}{7}, \sfr{1}{28}]
  \oplus [\sfr{1}{7}, \sfr{45}{28}],
  & [\sfr{1}{7}, \sfr{143}{224}]
  \oplus [\sfr{5}{7}, \sfr{15}{224}]
  \oplus [0, \sfr{57}{32}],
  \medskip\\
  {}[\sfr{19}{35}, \sfr{3}{28}]
  \oplus [\sfr{4}{35}, \sfr{15}{28}]
  \oplus [\sfr{2}{5}, \sfr{13}{4}],
  & [\sfr{2}{5}, \sfr{5}{32}]
  \oplus [\sfr{4}{35}, \sfr{99}{224}]
  \oplus [\sfr{19}{35}, \sfr{675}{224}],
\end{array}
\end{equation}
where $[h,k]$ denotes $U(h)\tensor L(\sfr{25}{28},k)$.
By considering the decomposition in \eqref{eq:5.2} together with 
the possible shapes in \eqref{eq:5.3}, we obtain the list of possible 
eigenvalues of $\o(f^1)$ as in the assertion.
\qed

In order to determine the $\o(f^1)$-spectrum of the Griess algebra of $V^\natural$, 
we need traces $\tr_{V^\natural_2}\mathrm{o}(f)^i$ for $1\leq i\leq 5$, which can 
be computed by the Matsuo-Norton trace formulae in \cite{Ma1}.

\begin{lem}[\cite{Ma1}]\label{lem:5.11}
  The traces of $\mathrm{o}(f^1)^i$, $1\leq i\leq 5$, on the Griess algebra 
  are as follows.
  \[
  \begin{array}{c}
    \tr_{V^\natural_2} \o(f^1)=\dfr{410175}{28}_,~~~~~
    \tr_{V^\natural_2} \o(f^1)^2=\dfr{2411375}{784}_,~~~~~
    \tr_{V^\natural_2} \o(f^1)^3=\dfr{27230625}{21952}_,
    \medskip\\
    \tr_{V^\natural_2} \o(f^1)^4=\dfr{793401325}{1229312}_,~~~~~~~~
    \tr_{V^\natural_2} \o(f^1)^5=\dfr{15221783625}{39337984}_.
  \end{array}
  \]
\end{lem}

By Lemma \ref{lem:5.11} and possible eigenvalues in Lemma \ref{lem:5.9},  
we can determine the $\o(f^1)$-spectrum compatible with the decompositions 
in Proposition \ref{prop:5.7}.
The result is as follows.

\begin{prop}\label{prop:5.12}
  The zero-mode $\mathrm{o}(f^1)$ acts on the decomposition in \eqref{eq:5.2} 
  as follows.
  \[
  \begin{array}{l}
    \begin{array}{ccccccccccccc}
      \l(\com_{V^\natural}\la a,b\ra\r)_2 &=& \mathbf{1} &+& \mathbf{1} &+& \mathbf{429}
      &+& \mathbf{3080} &+& \mathbf{13650} &+& \mathbf{13728}
      \\
      \mathrm{o}(f^1) &:& 2 && 0 && \sfr{3}{4} && 0 && 0 && \sfr{5}{32}
      \\
      \tau_{f^1} &:& 1 && 1 && 1 && 1 && 1 && -1
    \end{array}
    \bigskip \\
    \begin{array}{ccccccc}
      \mathrm{Top}~U(\sfr{5}{7})^c &=& \mathbf{1} &+& \mathbf{429} &+& \mathbf{352}
      \\
      \mathrm{o}(f^1) &:& \sfr{9}{7} && \sfr{1}{28} && \sfr{15}{224}
      \\
      \tau_{f^1} &:& 1 && 1 && -1
    \end{array}
    \bigskip \\
    \begin{array}{ccccccc}
      \mathrm{Top}~U(\sfr{19}{35})^c &=& \mathbf{78} &+& \mathbf{1430} &+& \mathbf{2080}
      \\
      \mathrm{o}(f^1) &:& \sfr{5}{14} && \sfr{3}{28} && \sfr{3}{224}
      \\
      \tau_{f^1} &:& 1 && 1 && -1
    \end{array}
    \\
    \begin{array}{ccccc}
      \mathrm{Top}~U(\sfr{2}{5})^c &=& \mathbf{3003} &+& \mathbf{2080}
      \\
      \mathrm{o}(f^1) &:& 0 && \sfr{5}{32}
      \\
      \tau_{f^1} &:& 1 && -1
    \end{array}
    \bigskip \\
    \begin{array}{ccccccccc}
      \mathrm{Top}~U(\sfr{1}{7})^c &=& \mathbf{1001} &+& \mathbf{10725} &+& \mathbf{352} 
      &+& \mathbf{13728}
      \\
      \mathrm{o}(f^1) &:& \sfr{5}{14} && \sfr{3}{28} && \sfr{143}{224} && \sfr{3}{224}
      \\
      \tau_{f^1} &:& 1 && 1 && -1 && -1
    \end{array}
   \bigskip \\
    \begin{array}{ccccccccc}
      \mathrm{Top}~U(\sfr{4}{35})^c &=& \mathbf{1430} &+& \mathbf{30030} &+& \mathbf{2080} 
      &+& \mathbf{27456}
      \\
      \mathrm{o}(f^1) &:& \sfr{15}{28} && \sfr{1}{28} && \sfr{99}{224} && \sfr{15}{224}
      \\
      \tau_{f^1} &:& 1 && 1 && -1 && -1
    \end{array}
  \end{array}
  \]
\end{prop}

\pf
Since $\tau_{f^1}$ is the central element in $C_{\Fi_{23}}(\tau_{f^1})\cong 2.\Fi_{22}$, 
it acts as 1 on the irreducible $2.\Fi_{22}$-components corresponding to 
the characters $\chi_i$ with $1\leq i\leq 10$ and as $-1$ on those to $\chi_i$ with 
$i\geq 66$ (cf.~\cite{ATLAS}).
If $\tau_{f^1}$ is trivial then the possible eigenvalues of $\o(f^1)$ are
0, 2, 9/7, 3/4, 5/14, 1/28, 3/28, 15/28 and if $\tau_{f^1}$ is non-trivial 
then the possible eigenvalues are 5/32, 3/224, 15/224, 99/224, 143/224 
by Lemma \ref{lem:5.9}.
This information together with Lemma \ref{lem:5.9} leads to the unique 
assignment of eigenvalues compatible with the set of traces in 
Lemma \ref{lem:5.11} and we obtain the eigenspace decompositions as in the 
assertion.
\qed

\subsection{Transpositions of $\Fi_{22}$ and $c=11/12$ Virasoro vectors}

Let $X^{[0]}=\la a,b\ra \subset X^{[1]}=\la a,b,x^1\ra \subset 
X^{[2]}=\la a, b, x^1, x^2\ra \subset V^\natural$ be the subalgebras of 
$V^\natural$ defined as in \eqref{eq:4.12}, and let $f^1$ and $f^2$ be the conformal 
vectors of $\com_{X^{[1]}}X^{[0]}$ and $\com_{X^{[2]}}X^{[1]}$, respectively.
Let 
\[
  D^{[1]}=\{ \tau_y \in D^{[0]} \mid \tau_y\tau_{x^1}=\tau_{x^1}\tau_y\},~~~
  G^{[1]}=\la \varphi^{[1]}(D^{[1]})\ra \subset \aut(\com_{V^\natural}\la a,b,x^1\ra ) 
\]
be defined as in \eqref{eq:4.18} and \eqref{eq:4.19}.
We have shown in Proposition \ref{prop:5.4} that $G^{[1]}\cong \Fi_{22}$.
In the following we will identify $\Fi_{22}$ with $G^{[1]}$.
Since $\la a,b,f^1\ra\cong \la a,b\ra\tensor \la f^1\ra$ is a full subalgebra of 
$X^{[1]}=\la a,b,x^1\ra$, the commutant $\com_{V^\natural}\la a,b,x^1\ra$ coincides with 
the 0-eigenspace of $\o(f^1)$ in $\com_{V^\natural}\la a,b\ra$.
Therefore, by Proposition \ref{prop:5.12}, 
we have the following decomposition of the Griess algebra of 
$\VF_{22}=\com_{V^\natural}\la a,b,x^1\ra$ as a $\Fi_{22}$-module.

\begin{prop}\label{prop:5.13}
As a $G^{[1]}\simeq \Fi_{22}$-module, the Griess algebra of 
$\com_{V^\natural}\la a,b,x^1\ra$ decomposes as follows.
\[
\renewcommand{\arraystretch}{0.8}
\begin{array}{ccccccc}
  \l(\com_{V^\natural} \la a,b,x^1\ra \r)_2 &=& \mathbf{1} &+& \mathbf{3080} 
  &+& \mathbf{13650},
  \\
  && \chi_1 && \chi_7 && \chi_9
\end{array}
\renewcommand{\arraystretch}{1}
\]
where $\chi_1$, $\chi_7$ and $\chi_9$ are irreducible characters of 
$\Fi_{22}$ labeled as in \cite[page 156]{ATLAS}.
\end{prop}

It is shown in Proposition \ref{prop:5.4} that $\varphi^{[1]}(\tau_{x^2})$ defines 
a 2A-element of $\Fi_{22}$ and its centralizer is 
$C_{\Fi_{22}}(\varphi^{[1]}(\tau_{x^2}))\cong 2.\mathrm{PSU}_6(2)$.
We have the following decompositions.

\begin{lem}\label{lem:5.14}
As a $C_{\Fi_{22}}(\varphi^{[1]}(\tau_{x^2}))\cong 2.\mathrm{PSU}_6(2)$-module, we have 
\[
\renewcommand{\arraystretch}{0.8}
\begin{array}{ccccccccccc}
  \mathbf{3080} &=& \mathbf{1} &+& \mathbf{252} &+& \mathbf{440} &+& \mathbf{1155} 
  &+& \mathbf{1232},
  \\
  \chi_7 && \mu_1 && \mu_4 && \mu_6 && \mu_{11} && \mu_{51}
  \bigskip\\
  \mathbf{13650} &=& \mathbf{252} &+& \mathbf{1155} &+& \mathbf{1155} &+& 
  \mathbf{4928}&+&\mathbf{6160},
  \\
  \chi_9 && \mu_4 && \mu_{12} && \mu_{13} && \mu_{20} && \mu_{57}
\end{array}
\renewcommand{\arraystretch}{1}
\]
where $\chi_i$ are the irreducible characters of $\Fi_{22}$ as in the previous 
proposition, and $\mu_i$ are the irreducible characters of $\mathrm{PSU}_6(2)$ and 
$2.\mathrm{PSU}_6(2)$ labeled as in \cite[page 116]{ATLAS}.
\end{lem}

It follows from Lemma \ref{lem:5.5}, Proposition \ref{prop:5.13} and 
Lemma \ref{lem:5.14} that the Griess algebra of the 
$C_{\Fi_{22}}(\varphi^{[1]}(\tau_{x^2}))$-invariants of 
$\com_{V^\natural}\la a,b,x^1\ra$ is 2-dimensional spanned 
by the $c=429/20$ conformal vector and the $c=11/12$ Virasoro vector $f^2$  
in $\com_{X^{[2]}}X^{[1]}$.
It also follows from (2) of Theorem \ref{thm:4.19} that $f^2$ is of $\sigma$-type 
and satisfies 
$\sigma_{f^2}=\varphi^{[1]}(\tau_{x^2})$ on $\com_{V^\natural} \la a,b,x^1\ra$. 
Therefore, we have obtained the Conway-Miyamoto correspondence between 
2A-involutions of $\Fi_{22}$ and $c=11/12$ Virasoro vectors of $\sigma$-type 
of $\com_{V^\natural}\la a,b,x^1\ra$.

\begin{thm}\label{thm:5.15}
  Let $E_{V^\natural}$ be the set of Ising vectors of $V^\natural$ and let $a$, 
  $b\in E_{V^\natural}$ be a pair such that  
  $\la a,b\ra\cong U_{3A}$, i.e., $(a\mymid b)=13\cd 2^{-10}$.   
  Set $I_{a,b}=\{ x\in E_{V^\natural} \mid \la a,x\ra\cong \la b,x\ra \cong U_{2A}\}$, 
  i.e., $(a\mymid x)=(b\mymid x)=2^{-5}$ for any $x\in I_{a,b}$, and take $x^1\in I_{a,b}$.
  Then there exists a Conway-Miyamoto correspondence between 2A-elements of 
  $\Fi_{22}$ and $c=11/12$ Virasoro vectors of $\sigma$-type of 
  $\com_{V^\natural}\la a,b,x^1\ra$.
  More precisely, the following hold.
  \\
  (1)~ Set 
  $D^{[1]}=\{ \tau_y \mid y\in I_{a,b},~ \tau_{x^1}\tau_y=\tau_y\tau_{x^1}\}$.
  Then $\la D^{[1]}\ra = C_\M(\tau_a,\tau_b,\tau_{x^1}) \cong 2.\Fi_{22}$  
  acts on $\com_{V^\natural}\la a,b,x^1\ra$ with the kernel $\la \tau_{x^1}\ra$.
  \\
  (2)~ Let $\varphi^{[1]} : \la D^{[1]}\ra \rightarrow 
  \aut(\com_{V^\natural}\la a,b,x^1\ra)$ be the homomorphism given in (1) and 
  $G^{[1]}$ the image of $\varphi^{[1]}$.
  Then $G^{[1]}\cong \Fi_{22}$ and $\varphi^{[1]}(D^{[1]}))$ is the set of 
  2A-involutions of $\Fi_{22}$.
  \\
  (3)~ Given a 2A-involution $t$ of $\Fi_{22}$, there exists a unique 
  $C_{\Fi_{22}}(t)$-invariant $c=11/12$ Virasoro vector $f$ of $\sigma$-type 
  of $\com_{V^\natural}\la a,b,x^1\ra$ such that $\sigma_f=t$.
\end{thm}

\begin{rem}\label{rem:5.16}
  In principle, we can continue to consider the next case 
  $X^{[2]}\subset X^{[3]}=\la a,b,x^1,x^2,x^3\ra$ and 
  study $\com_{V^\natural} X^{[2]}$ as a $G^{[2]}\cong \mathrm{PSU}_6(2)$-module.
  Let $f^3$ be the $c=c_7=14/15$ Virasoro vector of $\com_{X^{[3]}} X^{[2]}$.
  We know that $\tau_{f^3}=\varphi^{[2]}(\tau_{x^3})$ defines a 2A-element of $\mathrm{PSU}_2(6)$ 
  on $\com_{V^\natural} X^{[2]}$ by (1) of Theorem \ref{thm:4.19} and Proposition \ref{prop:5.4}.
  In order to verify the Conway-Miyamoto correspondence for $\mathrm{Fi}_{21}=\mathrm{PSU}_6(2)$, 
  we have to compute the $C_{\Fi_{21}}(\tau_{f^3})$-invariants of the Griess 
  algebra of $\com_{V^\natural} X^{[2]}$. 
  However, the Griess algebra of $\com_{V^\natural} X^{[1]}$ splits into many pieces 
  as in Lemma \ref{lem:5.14} and its 0-eigenspace of $\o(f^2)$ is technically difficult 
  to determine so that we cannot obtain the decomposition of the Griess algebra of 
  $\com_{V^\natural} X^{[2]}$ at this moment.
  It is likely that the $C_{\mathrm{Fi}_{21}}(\tau_{f^3})$-invariants of the Griess algebra 
  of of $\com_{V^\natural} X^{[2]}$ has a dimension more than 2, and the Conway-Miyamoto 
  correspondence for $\Fi_{21}$ seems to fail.
  However, even the $C_{\Fi_{21}}(\tau_{f^3})$-invariants is larger, we still have 
  a chance that it contains a unique $c=c_7$ Virasoro vector.
\end{rem}

\begin{appendices}
\section{Explicit constructions}

We will give explicit constructions of the subalgebras 
discussed in Section \ref{sec:4}. 

\subsection{Compact real form of a lattice VOA}\label{sec:a.1}

Let $L$ be an even positive definite lattice.
Let $V_L$ be the associated lattice VOA defined over $\C$ and $V_{L,\R}$ 
the one defined over $\R$ (cf.~\cite{FLM}).
Let $\theta$ be a lift of the $(-1)$-isometry of $L$ on $V_{L,\R}$ and 
$V_{L,\R}=V_{L,\R}^+\oplus V_{L,\R}^-$ the eigenspace decomposition such that 
$\theta$ acts by $\pm 1$ on $V_{L,\R}^\pm$.
Then the invariant bilinear form on $V_{L,\R}$ is positive definite on 
$V_{L,\R}^+$ and negative definite on $V_{L,\R}^-$.
Therefore, the real subspace 
\begin{equation}\label{eq:a.1}
  (V_L)_\R:=V_{L,\R}^+ \oplus \sqrt{-1} V_{L,\R}^-
\end{equation}
forms a compact real subalgebra of $V_L$.
In the following We will use the compact form $(V_L)_\R$ in \eqref{eq:a.1} of $V_L$.

\begin{rem}\label{rem:a.1}
  In the construction of a lattice vertex operator algebra $V_L$ one needs to 
  implement a 2-cocycle  $\varepsilon \in Z^2(L,\{\pm 1\})$ such that 
  $\varepsilon(a,b)\varepsilon(b,a)=(-1)^{(a|b)}$ for $a$, $b\in L$.
   If $L$ is of rank one or doubly even, then we can take $\varepsilon$ to be trivial.
\end{rem}

\subsection{Virasoro vectors associated to root systems}\label{sec:a.2}
We recall the construction of certain Virasoro vectors from \cite{DLMN}.
Let $R$ be a root lattice and $\Phi(R)$ its root system. 
Fix a system of simple roots and let $\Phi^+(R)$ and $\Phi^-(R)$ be the set of 
positive and negative roots, respectively.
We use $\sqrt{2}R$ to denote the $\Z$-submodule $\sqrt{2}\tensor_\Z R$ 
of $\R\tensor_\Z R$ which forms a doubly even lattice.
By \cite{DLMN}, the conformal vector $\omega_R$ of the lattice VOA $V_{\sqrt{2}R}$ is 
given by
\[
  \omega_R = \dfrac{1}{4h}\sum_{\al\in\Phi(R)}{\al_{(-1)}}^2 \vac ,
\]
where $h$ is the Coxeter number of $R$.
For $\alpha\in \sqrt{2}\,R$ with $(\alpha|\alpha)=4$ we set 
\begin{equation}\label{eq:a.2}
  w^\pm (\alpha) := \dfr{1}{16}{\alpha_{(-1)}}^2\vac \pm \dfr{1}{4}(e^\alpha + e^{-\alpha}).
\end{equation}
Then $w^\pm(\alpha)$ are Ising vectors (cf.~\cite{DMZ}). 
The following is straight forward.

\begin{lem}\label{lem:a.2}
  Let $\alpha$, $\beta\in \sqrt{2}R$ be squared norm 4 vectors.
  Then
  \[
    \big( w^\varepsilon(\alpha)\,\big|\,w^{\varepsilon'}(\beta)\big) =
    \begin{cases}
    2^{-2} &  \mbox{if $\alpha=\pm \beta$ and $\varepsilon=\varepsilon'$},
    \medskip\\
    0 & \mbox{if $(\alpha | \beta)=0$ or $\alpha=\pm \beta$ with $\varepsilon=-\varepsilon'$},
    \medskip\\
    2^{-5} & \mbox{if $(\alpha|\beta)=\pm 2$},
    \medskip\\
    \end{cases}
  \]
  where $w^\pm(\alpha)$ and $w^\pm(\beta)$ are defined as in \eqref{eq:a.2}.
  Moreover, we have
  \[
    w^\varepsilon(\alpha)_{(1)} w^{\epsilon'}(\beta)
    = \dfr{1}{4} \l( w^\varepsilon(\alpha)+w^{\varepsilon'}(\beta)
    -w^{-\varepsilon\varepsilon'}(\alpha-\beta)\r)
  \]
  if $(\alpha|\beta)=\pm 2$.
\end{lem}

We also define
\begin{equation}\label{cv1}
\begin{array}{ll}
  s_R
  & := \dfrac{4}{(h+2)}\dsum_{\alpha \in \sqrt{2}\,\Phi^+(R)} w^-(\alpha),
  \medskip\\
  t_R&: = \w_R - s_R = \dfrac{2}{h+2}\,\omega_R+\dfrac{1}{h+2}
  \dsum_{\alpha\in\sqrt{2}\,\Phi(R)}e^{\alpha}.
\end{array}
\end{equation}
Then it is shown in~\cite{DLMN} that $s_R$ and $t_R$ are mutually orthogonal 
Virasoro vectors.
The central charge of $t_R$ is
$2n/(n+3)$ if $R$ is of type $A_n$, $1$ if $R$ is of type $D_n$, and
$6/7$, $7/10$ and $1/2$ if $R$ is of type $E_6$, $E_7$ and $E_8$, respectively.
By Lemma \ref{lem:a.2}, the linear span of $w^-(\alpha)$ with 
$\alpha\in \sqrt{2}\,\Phi(R)$ forms a subalgebra of the Griess algebra of $V_{\sqrt{2}R}^+$.
Indeed, this is the Griess algebra of the subVOA generated by the Ising vectors 
$w^-(\alpha)$, $\alpha \in \sqrt{2}\,\Phi(R)$.
The Virasoro vector $s_R$ provides its conformal vector and by the orthogonality 
we have the following (cf.~Proposition 5.1 and Lemma 5.7 of \cite{LSY}).

\begin{prop}[\cite{LSY}]\label{prop:a.3}
  $\la w^-(\alpha) \mymid \alpha \in \sqrt{2}\,\Phi(R)\ra 
  = \com_{V_{\sqrt{2}R}^+} \la t_R\ra$.
\end{prop}

The subVOA $\la w^-(\alpha) \mymid \alpha \in \sqrt{2}\,\Phi(R)\ra$ is denoted 
by $M_R$ and its automorphism group is determined in \cite{LSY}.

\begin{rem}\label{EE8Ising}
If $R=E_8$, then
\begin{equation}\label{eE}
  t_{E_8}=\dfrac{1}{16} \omega_{E_8} +\dfrac{1}{32}\sum_{\al \in \sqrt{2}\,\Phi(E_8)} e^{\alpha}
\end{equation}
is an Ising vector. 
We say that $t_{E_8}$ is the \emph{standard Ising vector} of $V_{\sqrt{2}E_8}^+$.
Recall that the dual lattice $(\sqrt{2}E_8)^*=  \dfrac{1}{\sqrt{2}} E_8$.  
For $x\in  (\sqrt{2}E_8)^*$, define a $\Z$-linear map
\[
\begin{array}{ccccc}
  (x| \,\cdot\,) &:& \sqrt{2}E_8 &\longrightarrow & \Z_2 
  \\
  && y & \longmapsto& (x|y)  \mod 2.
\end{array}
\]
Clearly the map
\[
\begin{array}{ccccc}
  \varphi &:& (\sqrt{2}E_8)^* &\longrightarrow& \mathrm{Hom}_\Z(\sqrt{2}E_8, \Z_2)
  \\
  && x &\longmapsto& (x| \,\cdot\, )
\end{array}
\]
is a group homomorphism and $\ker \varphi= \sqrt{2}E_8$. 
The map  $(x|\,\cdot\,)$ thus induces an automorphism $\varphi_x$ of $V_{\sqrt{2}E_8}$ 
given by
\begin{equation}\label{phi_x}
  \varphi_x(u\otimes e^\al) 
  =(-1)^{(x|\al)} u\otimes e^\al \quad
  \text{ for } u\in M(1)\text{ and } \al \in \sqrt{2}E_8.
\end{equation}
In this case,
\[
  \varphi_x \,t_{E_8}
  = \dfrac{1}{16}\omega_{E_8} + \dfrac{1}{32}\sum_{\al \in \sqrt{2}\,\Phi(E_8)} 
  (-1)^{(x|\al)} e^{\al}
\]
is also an Ising vector.  
Since $\varphi_x$ commutes $\theta$ (the lift of the $-1$ isometry), 
the Ising vector $\varphi_x t_{E_8}$ is contained in $V_{\sqrt{2}E_8}^+$.   
The Ising vectors in $V_{\sqrt{2}E_8}^+$ have been classified in \cite{GrO,LSY,LS}. 
There are $240+256=496$ Ising vectors and they can be divided into two different types:
\medskip\\
{\bf $~~~A_1$-type:} $w^\pm (\al)$, $\al \in \sqrt{2}E_8$, $(\al|\al)=4$~~
($\abs{\Phi(E_8)}=240$ Ising  vectors),
\medskip\\
{\bf $~~~E_8$-type:} $\varphi_x\, t_{E_8}$ for some $x \in \dfrac{1}{\sqrt{2}}E_8$~~
($2^8=256$ Ising vectors).
\end{rem}

\begin{rem}\label{yj}
  Let $\al_1, \dots, \al_n$ be the simple roots of a root lattice of type $A_n$, 
  that is, $(\al_i|\al_i)=2$,  $(\al_i|\al_j)= -1$ if $|i-j|= 1$ and $(\al_i|\al_j)= 0$ 
  otherwise.  
  Then 
  \[
    \Span_\Z\{ \al_1, \dots, \al_k\}\cong A_k
  \] 
  for any $1 \leq k \leq n$ and we have  a sequence of sublattices
  \[
    A_1\subset A_2\subset \cdots \subset A_n.
  \]
  For $1 \leq k \leq n$, denote $s^k=s_{A_k}$. 
  We also set
  \begin{equation}\label{eq:a.6}
    \eta^1 = s^1\quad \text{ and } \quad \eta^k= s^k -s^{k-1}
  \end{equation}
  for $2\leq k\leq n$. 
  It is shown in \cite{DLMN} (see also \cite{GKO}) that $\eta^k$ is a simple 
  Virasoro vector of central charge
  \[
    c_k=1- \dfrac{6}{(k+2)(k+3)},~~~ 1\leq k\leq n.
  \]  
  Moreover, we have the following orthogonal decomposition of the conformal vector 
  of $V_{\sqrt{2}A_n}$.
  \begin{equation}\label{eq:a.7}
    \w_{A_n}
    =\eta^1 \dotplus \eta^2 \dotplus \dots \dotplus \eta^n \dotplus t_{A_n}.
  \end{equation}
\end{rem}

\subsection{Construction of $X^{[n]}$}\label{sec:a.3}



Let $n\geq 1$ and $\epsilon_1, \dots  \epsilon_{n+5}$ an orthogonal basis of $\R^n$ 
such that $(\epsilon_i|\epsilon_j)=2\delta_{i,j}$.
Set 
\[
  K^{[k]}=\Z \epsilon_1\oplus \Z \epsilon_2 \oplus \cds \oplus \Z \epsilon_k,
  ~~ 1\leq k\leq n+5.
\]
Then $K^{[k]}\cong  A_1^{\oplus k}$.
Set $\gamma = \fr{1}{2}(\epsilon_1+\epsilon_2+\epsilon_3+\epsilon_4)$ 
and $L=K^{[n+5]}\sqcup (K^{[n+5]}+\gamma)$.
Then $L\cong D_4\oplus A_1^{\oplus (n+5)}$. 
We consider lattice VOAs 
\[
  V_{K^{[1]}}\subset V_{K^{[2]}} \subset \cds \subset V_{K^{[n+5]}} \subset V_L.
\] 
Note that $V_{K^{[n+5]}}$ is the full subVOA of $V_L$.
We denote the conformal vector of $V_L$ by $\w_L$.

\begin{rem}\label{rem:a.6}
  Since $(\alpha|\beta)\in 2\Z$ for any $\alpha, \beta\in K^{[n+5]}$, we can define $V_{K^{[n+5]}}$ with a 
  trivial 2-cocycle in $Z^2(K^{[n+5]},\{\pm 1\})$.
  In the construction of $V_L$, we choose a 2-cocycle in $Z^2(L,\{\pm 1\})$ such 
  that its restriction on $K^{[n+5]}$ is trivial.
  This is possible since we can form $V_L=V_{K^{[n+5]}}\oplus V_{K^{[n+5]}+\gamma}$ 
  as a $\Z_2$-graded simple current extension of $V_{K^{[n+5]}}$.
\end{rem}

For $1\leq k\leq n+5$, we set
\begin{equation}\label{eq:a.8}
\begin{array}{l}
  H^k={{\epsilon_1}_{(-1)}}^2\vac +{{\epsilon_2}_{(-1)}}^2\vac + \cds 
  + {{\epsilon_k}_{(-1)}}^2\vac,
  \medskip\\
  E^k=e^{\epsilon_1}+e^{\epsilon_2}+\cds +e^{\epsilon_k},
  \medskip\\
  F^k=e^{-\epsilon_1}+ e^{-\epsilon_2}+\cds +e^{-\epsilon_k}.
\end{array}
\end{equation}
Then $H^k$, $E^k$, $F^k$ form an $\mathfrak{sl}_2$-triple in the weight one subspace
of $V_{K^{[k]}}$ and generate a subVOA isomorphic to level $k$ affine VOA 
$L_{\hat{\mathfrak{sl}}_2}(k,0)$ associated to $\hat{\mathfrak{sl}}_2$.  Let
\begin{equation}\label{eq:a.9}
  \Omega^k=\dfr{1}{4(k+2)}\l( H^k_{(-1)}H^k+2E^k_{(-1)}F^k+2F^k_{(-1)}E^k\r)
\end{equation}
be the Sugawara element.
Then $\Omega^k$ is the conformal vector of $L_{\hat{\mathfrak{sl}}_2}(k,0)$ (cf.~\cite{FZ}).
Consider the sublattice 
\begin{equation}\label{eq:a.10}
  M^{[k]}=\{ \al \in K^{[k+1]}\mid (\al\mid \epsilon_1+\cdots + \epsilon_{n+5})=0\},~~~
\end{equation}
Then $M^{[k]}\cong \sqrt{2}A_k$ and we can define the Virasoro vectors 
\begin{equation}\label{eq:a.11}
  s^k=\dfr{4}{k+3}\dsum_{1\leq i<j\leq k+1}w^-(\epsilon_i-\epsilon_j),~~~
  \eta^1=w^-(\epsilon_1-\epsilon_2),~~~
  \eta^k=s^k-s^{k-1}
\end{equation}
for $2\leq k\leq n+4$ as in \eqref{cv1} and \eqref{eq:a.6}. 
It is shown in \cite{GKO,LYd} that the Virasoro vectors $\eta^k$ can be written as 
\[
  \eta^k=\Omega^k+\dfr{1}{4}{{\epsilon_{k+1}}_{(-1)}}^2\vac-\Omega^{k+1}.
\]
Therefore, we have the following orthogonal decompositions. 
\begin{equation}\label{eq:a.12}
  \w_L=s^{n+4} \dotplus \Omega^{n+5}
  = \eta^1\dotplus \eta^2 \dotplus \cds  \dotplus \eta^{n+4} \dotplus \Omega^{n+5}.
\end{equation}
As a consequence, we have the following  full subVOA of $V_L$.
\begin{equation}\label{eq:a.13}
  L(c_1,0)\tensor L(c_2,0)\tensor \cds \tensor L(c_{n+4},0)\tensor 
  L_{\hat{\mathfrak{sl}}_2}(n+5,0).
\end{equation}
We note that all $\eta^k$, $1\leq k\leq n+4$, and $\Omega^{n+5}$ are defined 
inside the compact form $V_{L,\R}^+$.

\begin{lem}\label{lem:a.7}
  The weight one subspace of $\com_{V_L}\la\eta^1,\eta^2,\Omega^{n+5}\ra$ is 
  trivial.
\end{lem}

\pf
It suffices to show $\ker_{V_L} {\eta^1}_{(1)} \cap \ker_{V_L} {\eta^2}_{(1)} 
\cap\ker_{V_L} {\Omega^{n+5}}_{(1)} \cap (V_L)_1 =0$.
Since $V_L=V_{K^{[n+5]}}\oplus V_{K^{[n+5]}+\gamma}$, we have a decomposition 
$(V_L)_1=(V_{K^{[n+5]}})_1\oplus (V_{K^{[n+5]}+\gamma})_1$.
The weight one subspace of $V_{K^{[n+5]}}$ is $3(n+5)$-dimensional spanned by 
${\epsilon_i}_{(-1)}\vac$ and $e^{\pm\epsilon_i}$ for $1\leq i\leq n+5$.
By a direct computation one has
\[
\begin{array}{l}
  {\Omega^{n+5}}_{(1)} {\epsilon_i}_{(-1)}\vac
  =\dfr{2}{n+7} {\epsilon_i}_{(-1)}\vac +\dfr{1}{n+7} H^{n+5}, 
  \medskip\\
  {\Omega^{n+5}}_{(1)} e^{\epsilon_i}
  =\dfr{2}{n+7}e^{\epsilon_i} +\dfr{1}{n+7} E^{n+5}, 
  \medskip\\
  {\Omega^{n+5}}_{(1)} e^{-\epsilon_i}
  =\dfr{2}{n+7}e^{-\epsilon_i} +\dfr{1}{n+7}F^{n+5}.
\end{array}
\]
Therefore, the characteristic polynomial of ${\Omega^{n+5}}_{(1)}$ on $(V_{K^{[n+5]}})_1$ is 
\[
  (x-1)^3\l( x-\dfr{2}{n+7}\r)^{3n+12}
\] 
and hence $\ker_{V_L} {\Omega^{n+5}}_{(1)} \cap (V_{K^{[n+5]}})_1=0$.
The weight one subspace of $V_{K^{[n+5]}+\gamma}$ is spanned by the vectors 
$e^{\alpha}$ with $\alpha=\fr{1}{2}(\pm \epsilon_1\pm \epsilon_2\pm \epsilon_3\pm \epsilon_4)$. 
Set $\gamma_1=\fr{1}{2}(\epsilon_1+\epsilon_2-\epsilon_3-\epsilon_4)$, 
$\gamma_2=\fr{1}{2}(\epsilon_1-\epsilon_2-\epsilon_3+\epsilon_4)$,
$\gamma_3=\fr{1}{2}(\epsilon_1-\epsilon_2+\epsilon_3-\epsilon_4)$ and 
$\tilde{e}(\gamma_i)=e^{\gamma_i}+e^{-\gamma_i}$ for $1\leq i\leq 3$.
Then by a direct computation one has
\[
\begin{array}{l}
  \ker_{V_L} {\Omega^{n+5}}_{(1)} \cap (V_{K^{[n+5]}+\gamma})_1
  =\Span\{\, \tilde{e}(\gamma_1)-\tilde{e}(\gamma_2),\,
   \tilde{e}(\gamma_2)-\tilde{e}(\gamma_3)\,\},
  \medskip\\
  \ker_{V_L} {\eta^1}_{(1)}\cap \ker_{V_L} {\Omega^{n+5}}_{(1)} \cap (V_{K^{[n+5]}+\gamma})_1
  = \Span\{\, 2\tilde{e}(\gamma_1)-\tilde{e}(\gamma_2)-\tilde{e}(\gamma_3)\,\} ,
  \medskip \\
  {\eta^2}_{(1)}\l( 2\tilde{e}(\gamma_1)-\tilde{e}(\gamma_2)-\tilde{e}(\gamma_3)\r)
  =\dfr{3}{5}\l( 2\tilde{e}(\gamma_1)-\tilde{e}(\gamma_2)-\tilde{e}(\gamma_3)\r) .
\end{array}
\]
Therefore, $\ker_{V_L} {\eta^1}_{(1)} \cap \ker_{V_L} {\eta^2}_{(1)} 
\cap\ker_{V_L} {\Omega^{n+5}}_{(1)} \cap (V_{K^{[n+5]}+\gamma})_1 =0$ 
and thus the weight one subspace of $\com_{V_L}\la\eta^1,\eta^2,\Omega^{n+5}\ra$ 
is trivial.
\qed

\begin{cor}\label{cor:a.8}
  The subalgebra $\com_{V_{L,\R}^+}\la \eta^1,\eta^2,\Omega^{n+5}\ra$ is a compact 
  VOA of OZ-type.
\end{cor}

We will realize the (2A,3A)-generated subalgebra 
$X^{[n]}=\la a,b,x^1,\dots,x^n\ra$ in Section \ref{sec:4.3} as a subalgebra of 
$\com_{V_{L,\R}^+}\la\eta^1,\eta^2,\Omega^{n+5}\ra$.
Set
\begin{equation}\label{eq:a.14}
\begin{array}{ll}
  q
  &:= \dsum_{k=1}^3 (2\gamma+\epsilon_k-\epsilon_4-4\epsilon_5)_{(-1)} 
  \l( e^{\gamma -\epsilon_k-\epsilon_4}-e^{-\gamma +\epsilon_k+\epsilon_4}\r)
  \\
  &~~~~
   -4 \dsum_{k=1}^3 \l( e^{\gamma -\epsilon_k-\epsilon_5}
  +e^{-\gamma +\epsilon_k+\epsilon_5}\r)
 +12\l( e^{\gamma-\epsilon_4-\epsilon_5}+e^{-\gamma+\epsilon_4+\epsilon_5}\r) \in V_{L,\R}^+.
\end{array}
\end{equation}
The following is shown in \cite{SY}.

\begin{lem}[\cite{SY}]\label{lem:a.9}
  The vector $q$ is a highest weight vectors for $\eta^1$, $\eta^2$, 
  $\eta^3$, $\eta^4$, $\Omega^{n+5}$ and $w^-(\epsilon_4-\epsilon_5)$ 
  with highest weights $0$, $0$, $2/3$, $4/3$, $0$ and $1/16$, respectively. 
  In particular, $q\in \com_{V_{L,\R}^+}\la\eta^1,\eta^2,\Omega^{n+5}\ra$.
\end{lem}


We set
\begin{equation}\label{eq:a.15}
  a:=w^-(\epsilon_4-\epsilon_5),~~
  b:=-\dfr{1}{2}a+\dfr{15}{64}\eta^3+\dfr{21}{32}\eta^4+\dfr{\sqrt{3}}{2^7}q,~~
  c:=\tau_a b.
\end{equation}
Then it is shown in \cite{SY} that $a$ and $b$ are Ising vectors in $V_L$ 
such that $\la a,b\ra$ is isomorphic to the 3A-algebra with the characteristic 
Virasoro frame $\eta^3\dotplus \eta^4$.
For $1\leq i\leq n$, we set 
\begin{equation}\label{eq:a.16}
  x^i:=w^-(\epsilon_5-\epsilon_{i+5}).
\end{equation}
Using Lemma \ref{lem:a.2} we can verify the following.

\begin{lem}\label{lem:a.10}
  $(a\mymid x^i)=(b \mymid x^i)=(x^j \mymid x^k)=2^{-5}$ 
  for $1\leq i \leq n$ and $1\leq j<k\leq n$.
\end{lem}

\begin{lem}\label{lem:a.11}
  $\la a,b,x^1,\dots,x^n\ra$ is a subalgebra of 
  $\com_{V_{L,\R}^+}\la \eta^1,\eta^2,\Omega^{n+5}\ra$.
\end{lem}

\pf
It follows from Proposition \ref{prop:a.3} and Eq.~\eqref{eq:a.12} that 
$a=w^-(\epsilon_4-\epsilon_5)$ and $x^i=w^-(\epsilon_5-\epsilon_{i+5})$, $1\leq i\leq n$, 
belong to the commutant of $\la \Omega^{n+5}\ra$.
It also follows from Proposition \ref{prop:a.3} that $s^2=\eta^1\dotplus \eta^2$ is 
the conformal vector of $\la w^-(\epsilon_1-\epsilon_2),w^-(\epsilon_2-\epsilon_3)\ra$. 
By Lemma \ref{lem:a.2}, $a$ and $x^i$, $1\leq i\leq n$, are orthogonal 
to both of $w^-(\epsilon_1-\epsilon_2)$ and $w^-(\epsilon_2-\epsilon_3)$. 
Therefore $\la a,x^1,\dots,x^n\ra \subset \com_{V_{L,\R}^+}\la \eta^1,\eta^2,\Omega^{n+5}\ra$.
That $b\in \com_{V_{L,\R}^+}\la \eta^1,\eta^2,\Omega^{n+5}\ra$ follows from 
Lemma \ref{lem:a.9}.
Hence, $\la a,b,x^1,\dots,x^n\ra$ is a subalgebra of 
$\com_{V_{L,\R}^+}\la \eta^1,\eta^2,\Omega^{n+5}\ra$.
\qed

\begin{rem}\label{rem:a.12}
  The Ising vectors in \eqref{eq:a.15} and \eqref{eq:a.16} are defined 
  inside $V_{L,\R}^+$ so that $\la a,b,x^1,\dots,x^n\ra$ has a compact real form 
  as a subalgebra of $\com_{V_{L,\R}^+}\la \eta^1,\eta^2,\Omega^{n+5}\ra$ which is of 
  OZ-type by Corollary \ref{cor:a.8}.
  Since $\la a,b\ra$ is the 3A-type and all $\la a,x^i\ra$, $\la b,x^i\ra$ 
  and $\la x^j,x^k\ra$ are the 2A-type by Lemma \ref{lem:a.10}, 
  we have obtained a realization of the (2A,3A)-generated subalgebra 
  $X^{[n]}=\la a,b,x^1,\dots,x^n\ra$ discussed in Section \ref{sec:4.3} 
  inside $V_{L,\R}^+$.
\end{rem}

\begin{prop}\label{prop:a.13}
  $\la a,b,x^1,\dots,x^n\ra$ has the Virasoro frame 
  $\eta^3\dotplus \eta^4\dotplus \cds \dotplus \eta^{n+4}$.
  In particular, $\la a,b,x^1,\dots,x^n\ra$ is a full subalgebra of 
  $\com_{V_{L,\R}^+}\la\eta^1,\eta^2,\Omega^{n+5}\ra$.
\end{prop}

\pf
Since $\la a,b\ra$ is isomorphic to the 3A-algebra with Virasoro frame 
$\eta^3\dotplus \eta^4$, we have 
\[
  a+b+c=\dfr{15}{32}\eta^3+\dfr{21}{16}\eta^4.
\]
By \eqref{eq:a.11} we have
\begin{equation}\label{eq:a.17}
\begin{array}{ll}
  a+b+c&=\dfr{15}{32}\eta^3+\dfr{21}{16}\eta^4
  = \dfr{15}{32}(s^3-s^2)+\dfr{21}{16}(s^4-s^3)
  \medskip\\
  &= -\dfr{15}{32}s^2+\dfr{3}{16}\dsum_{1\leq i<j\leq 4}w^-(\epsilon_i-\epsilon_j)
  +\dfr{3}{4}\dsum_{1\leq i\leq 4}w^-(\epsilon_i-\epsilon_5).
  \end{array}
\end{equation}
It is clear from expressions that $s^2$ and $w^-(\epsilon_i-\epsilon_j)$ 
are orthogonal to $w^-(\epsilon_5-\epsilon_{k+5})$ for $1\leq i<j\leq 4$ and 
$1\leq k\leq n$, whereas 
$(w^-(\epsilon_i-\epsilon_5)\mymid w^-(\epsilon_5-\epsilon_{k+5}))=2^{-5}$ 
and one has 
$w^-(\epsilon_i-\epsilon_5)\circ w^-(\epsilon_5-\epsilon_{k+5})
=w^-(\epsilon_i-\epsilon_{k+5})$ by Lemma \ref{lem:a.2}.
Therefore we have
\begin{equation}\label{eq:a.18}
\begin{array}{l}
  (a+b+c)\circ x^k
  = -\dfr{15}{32}s^2+\dfr{3}{16}\dsum_{1\leq i<j\leq 4}w^-(\epsilon_i-\epsilon_j)
  +\dfr{3}{4}\dsum_{1\leq i\leq 4}w^-(\epsilon_i-\epsilon_{k+5}).
\end{array}
\end{equation}
So we have 
\[
\begin{array}{l}
  a+b+c+\dsum_{k=1}^n (a+b+c)\circ x^k
  \medskip\\
  =-\dfr{15(n+1)}{32}s^2-\dfr{3(n+1)}{16}\dsum_{1\leq i<j\leq 4}w^-(\epsilon_i-\epsilon_j)
  +\dfr{3}{4}\dsum_{1\leq i\leq 4\atop 5\leq k\leq n+5}w^-(\epsilon_i-\epsilon_k).
\end{array}
\]
Now we consider the conformal vector $\w^n$ of $\la a,b,x^1,\dots,x^n\ra$ 
in \eqref{eq:4.13}. 
\[
\begin{array}{l}
  \w^n
  =\dfr{3(3-n)}{2(n+7)}\eta^3 + \dfr{16}{3(n+7)}\l( a+b+c
  +\dsum_{i=1}^n (a+b+c)\circ x^i\r)
 \medskip\\
  ~~~~~~~
  +\dfr{4}{n+7} \l(\dsum_{i=1}^n x^i +\dsum_{1\leq j<k\leq n} x^j\circ x^k\r) 
  \medskip\\
  =\dfr{3(3-n)}{2(n+7)}\l( \dfr{2}{3}\dsum_{1\leq i\leq 4} w^-(\epsilon_i-\epsilon_j)
  -s^2\r) 
  \medskip\\
  ~~+ \dfr{16}{3(n+7)}\l( -\dfr{15(n+1)}{32}s^2
  +\dfr{3(n+1)}{16}\dsum_{1\leq i<j\leq 4}w^-(\epsilon_i-\epsilon_j)
  +\dfr{3}{4}\dsum_{1\leq i\leq 4\atop 5\leq k\leq n+5}w^-(\epsilon_i-\epsilon_k)\r)
  \medskip\\
  ~~+\dfr{4}{n+7} \l(\dsum_{i=1}^n w^-(\epsilon_5-\epsilon_{i+5}) 
  +\dsum_{1\leq j<k\leq n} w^-(\epsilon_{j+5}-\epsilon_{k+5})\r) 
  \medskip\\
  = -s^2+\dfr{4}{n+7}\dsum_{1\leq i<j\leq n+5} w^-(\epsilon_i-\epsilon_j)
  = s^{n+4}-s^2.
\end{array}
\]
Thus, we have 
\[
  \w^n=s^{n+4}-s^2= (s^3-s^2)+(s^4-s^3)+\cds +(s^{n+5}-s^{n+4})
  =\eta^3\dotplus \eta^4\dotplus \cds \dotplus \eta^{n+4}.
\]
Since $\w^k$ is the conformal vector of a subalgebra 
$\la a,b,x^1,\dots,x^k\ra$ of $\la a,b,x^1,\dots,x^n\ra$ for $1\leq k\leq n$, 
the Virasoro vector $\eta^{k+4}=\w^k-\w^{k-1}$ belongs to $\la a,b,x^1,\dots,x^n\ra$.
Therefore, $\la a,b,x^1,\dots,x^n\ra$ has a Virasoro frame 
$\w^n=\eta^3\dotplus \eta^4 \dotplus \cds \dotplus \eta^n$.
\qed
\medskip

Denote $Y^{[0]}=\la a,b\ra\subset V_{L,\R}^+$ and 
$Y^{[k]}=\la a,b,x^1,\dots,x^k\ra\subset V_{L,\R}^+$ for $1\leq k\leq n$.
Then $Y^{[k]}$ provides a realization of $X^{[k]}$ defined in \eqref{eq:4.12}.
Using the Virasoro frame 
$\eta^3\dotplus \eta^4 \dotplus \cds \dotplus \eta^{n+4}$ we can prove that the commutant 
$\com_{Y^{[n]}} Y^{[n-1]}$ is generated by the $c=c_{n+4}$ Virasoro vector 
$\eta^{n+4}=\w^n-\w^{n-1}$ which is denoted by $f^n$ in Theorem \ref{thm:4.17}.
Consider the full subalgebra of $V_L$ displayed in \eqref{eq:a.13}.
The irreducible decomposition with respect to this subalgebra is known as 
the GKO construction.
By Eq.~(2.20) of \cite{GKO} (see also Lemma 3.1 of \cite{LLY}), 
we have the following decompositions:
\begin{equation}\label{eq:a.19}
\begin{array}{l}
  \ds 
  V_{K^{[n+5]}}=
  \!\!\!\!
  \bigoplus_{1\leq i_k\leq k+1 \atop i_k \equiv 1 \,\mathrm{mod}\, 2}\!\!\!\!
  L(c_1,h^{(1)}_{i_1,i_2})\tensor 
  \cds \tensor L(c_{n+4},h^{(n+4)}_{i_{n+4},i_{n+5}})\tensor L_{\hat{\mathfrak{sl}}_2}(n+5,i_{n+5}-1), 
  \medskip\\
  \ds
  V_{K^{[n+5]}+\gamma}=
  \!\!\! \!\!\!\! 
  \bigoplus_{1\leq j_k\leq k+1 \atop j_k \equiv 1+\delta_k \,\mathrm{mod}\, 2}\!\!\!\!\!\!\!
  L(c_1,h^{(1)}_{j_1,j_2})\tensor 
  \cds \tensor L(c_{n+4},h^{(n+4)}_{j_{n+4},j_{n+5}})\tensor L_{\hat{\mathfrak{sl}}_2}(n+5,j_{n+5}-1),
\end{array}
\end{equation}
where $\delta_1=1$, $\delta_2=0$, $\delta_3=1$ and $\delta_k=0$ for 
$4\leq k\leq n+5$.
By the decompositions above, we see that the commutant of 
$\la \eta^1,\eta^2,\dots,\eta^{n+3},\Omega^{n+5}\ra$ in $V_{L,\R}$ is exactly $\la \eta^{n+4}\ra$.
Therefore, $\com_{Y^{[n]}} Y^{[n-1]}$ is generated by the $c=c_{n+4}$ Virasoro vector 
$\eta^{n+4}=\w^n-\w^{n-1}$. 

\begin{conj}
  The subalgebra $X^{[n]}=\la a,b,x^1,\dots,x^n\ra$ defined in \eqref{eq:4.12} is unique 
  up to isomorphism and coincides with the subalgebra 
  $\com_{V_{L,\R}^+}\la \eta^1,\eta^2,\Omega^{n+5}\ra$ of $V_{L,\R}^+$.
  In particular, $\com_{X^{[n]}} X^{[n-1]}$ is generated by the $c=c_{n+4}$ Virasoro vector 
  $f^n$ in Theorem \ref{thm:4.17}.
\end{conj}

\begin{rem}
  We will construct $V^\natural$ as an extension of 
  $\com_{V_{L,\R}^+}\la \eta^1,\eta^2,\Omega^{n+5}\ra$ in \cite{CLY}.
\end{rem}

Finally, we prove that the Virasoro vectors in \eqref{eq:4.15} are linearly independent.
By \eqref{eq:a.16} and \eqref{eq:a.17}, we know that the following vectors are 
linearly independent in $V_{K^{[n+5]}}$.
\[
\begin{array}{l}
  u_{a,b}=\eta^3,~~
  x^i=w^-(\epsilon_5-\epsilon_{i+5}),~~
  a=w^-(\epsilon_4-\epsilon_5),~~
  a\circ x^i=w^-(\epsilon_4-\epsilon_{i+5}), ~~
  \medskip\\
  a+b+c,~~ 
  (a+b+c)\circ x^i,~~
  x^j\circ x^k=w^-(\epsilon_{j+5}-\epsilon_{k+5}), ~~
  1\leq i\leq n,~~
  1\leq j<k\leq n.
\end{array}
\]
On the other hand, by \eqref{eq:a.14} and \eqref{eq:a.15} we have 
\[
\begin{array}{ll}
  q=\dfr{2^6}{\sqrt{3}}(b-c)
  &\!\!\!\! = \dsum_{k=1}^3 (2\gamma+\epsilon_k-\epsilon_4-4\epsilon_5)_{(-1)} 
  \l( e^{\gamma -\epsilon_k-\epsilon_4}-e^{-\gamma +\epsilon_k+\epsilon_4}\r)
  \\
  &~~~~
   -4 \dsum_{k=1}^3 \l( e^{\gamma -\epsilon_k-\epsilon_5}
  +e^{-\gamma +\epsilon_k+\epsilon_5}\r)
 +12\l( e^{\gamma-\epsilon_4-\epsilon_5}+e^{-\gamma+\epsilon_4+\epsilon_5}\r) .
\end{array}
\]
Since $\sigma_{x^i}$ acts on $V_L$ by a reflection
associated to $\epsilon_5-\epsilon_{i+5}$, we have 
\[
\begin{array}{ll}
  \dfr{2^6}{\sqrt{3}}(b-c)\circ x^i
  &\!\!\!\! = \dsum_{k=1}^3 (2\gamma+\epsilon_k-\epsilon_4-4\epsilon_{i+5})_{(-1)} 
  \l( e^{\gamma -\epsilon_k-\epsilon_4}-e^{-\gamma +\epsilon_k+\epsilon_4}\r)
  \\
  &~~~~ 
   -4 \dsum_{k=1}^3 \l( e^{\gamma -\epsilon_k-\epsilon_{i+5}}
  +e^{-\gamma +\epsilon_k+\epsilon_{i+5}}\r)
  +12\l( e^{\gamma-\epsilon_4-\epsilon_{i+5}}+e^{-\gamma+\epsilon_4+\epsilon_{i+5}}\r) .
\end{array}
\]
From the expressions above, we see that 
$b-c$ and $(b-c)\circ x^i$, $1\leq i\leq n$, are linearly independent 
in $V_{K^{[n+5]}+\gamma}$.
Thus, by considering the Gram matrix, all the Virasoro vectors in \eqref{eq:4.15} are 
also linearly independent.

\subsection{Construction of $\la a,b,x,y\ra$ with $\la x,y\ra \cong U_{\mathrm{3A}}$}\label{sec:a.4}
In this subsection, we will construct explicitly four Ising vectors $a,b,x,y$ in the 
VOA $V_\Lambda^+$ associated to the Leech lattice $\Lambda$ such that 
$\la a,b\ra \cong \la x,y\ra\cong U_{\mathrm{3A}}$ and 
$\la a,x\ra \cong \la b, x\ra\cong \la a,y\ra \cong \la b, y\ra\cong U_{\mathrm{2A}}$.
Therefore, the Griess algebra discussed in Section \ref{sec:4.2} does exist.

\medskip

For explicit calculations,  we will use the notion of \textit{hexacode balance\/} (or
MOG) to denote the codewords of the Golay code $\mathcal{G}$ and the vectors in the Leech
lattice $\Lambda$~\cite{cs,Gr}. Namely,  we arrange the set $\Omega=\{1,2, \ldots,24\}$
into a $4\times 6$ array such that the six columns forms a sextet.
Recall that the Conway group $\Co_0$ is the orthogonal group $\mathrm{O}(\Lambda)$ of 
$\Lambda$ and contains 4 conjugacy classes of involutions (see \cite{cs}). 
If $i$ is an involution with trace $8$ on $\Lambda$, then the $(-1)$-eigenlattice is 
isomorphic to $\sqrt{2}E_8$.

Consider the following 4 involutions in $\mathrm{O}(\Lambda)$.

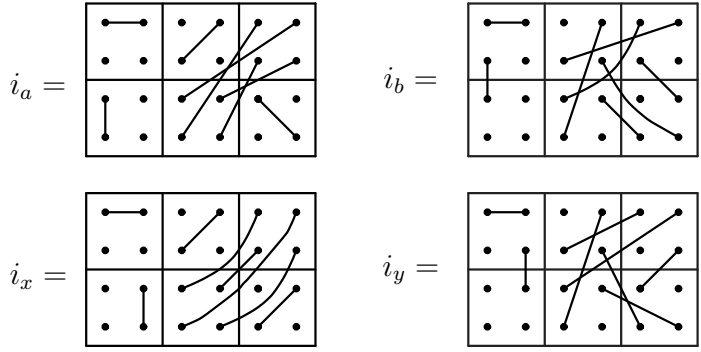
\begin{center}
\input{fig1.tex}
\end{center}

\begin{center}
\input{fig2.tex}
\end{center}
The automorphisms above should be viewed as permutations on $24$ coordinates and
they act on the Leech lattice $\Lambda$ from the left.

Note that
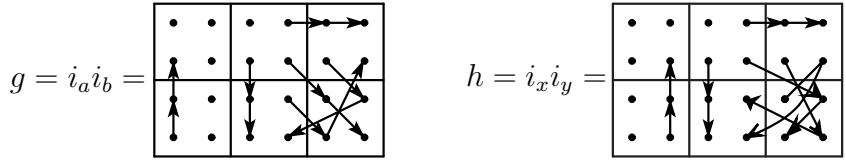
\begin{center}
\input{fig3.tex}
\end{center}
where a composition of above permutations is executed from right to left.

We have $\la i_a, i_b\ra \cong \la i_x, i_y\ra \cong \mathrm{S}_3$. 
Note that $i_x$ and $i_y$ commute with $\la i_a, i_b\ra$. 
Thus, we also have  $\la i_a, i_b, i_x\ra \cong \la i_a,i_b, i_y\ra \cong 
\mathrm{Dih}_{12}$, a dihedral group of order 12.

\begin{nota}
Let $A$, $B$, $X$, $Y$ be the $(-1)$-eigenlattices of $i_a$, $i_b$, $i_x$ and $i_y$, 
respectively. 
Then  $A\cong B\cong X\cong Y\cong \sqrt{2}E_8$. 
Moreover, $A+B\cong X+Y \cong \dih{6}{14}$ and $A+B+X\cong A+B+Y \cong \dih{12}{16}$ 
by the analysis in \cite{GL}.
\end{nota}

\begin{nota}[cf.~Remark \ref{EE8Ising}]\label{eM}
For any $\sqrt{2}E_8$-sublattice $M$, let
\[
  t_M=\dfrac{1}{16} \omega_M +\dfrac{1}{32}\sum_{\al\in M\atop (\al|\al)=4} e^\al
\]
be the standard Ising vector in $V_M^+$.
\end{nota}

Since $A$ and $B$ are doubly even, it is possible to choose a 2-cocycle 
$\varepsilon\in \Z^2(\Lambda,\{\pm 1\})$ such that $\varepsilon$ is trivial on 
$A$ and $B$ (see \cite[Notation 5.38]{LY2}). 
Let $a:=t_A$ and $b:=t_B$ be the standard Ising vector in $V_A^+$ and $V_B^+$, respectively.

\begin{lem}\label{abx}
Let $M$, $N$, $E$ be $\sqrt{2}E_8$-sublattices of the Leech lattice. 
Suppose $M+N\cong \dih{6}{14}$ and $M+E\cong N+E\cong \dih{4}{12}$.
Then there exists an Ising vector $e\in V_E^+$ such that $(t_M \mymid e)=(t_N \mymid e)=2^{-5}$.
\end{lem}

\pf
Let $i_L$ be the SSD involution associated to a $\sqrt{2}E_8$-sublattice $L$.
Since $i_Ei_N$ has order 2 and trace $8$, the $(-1)$-eigenlattice of $i_Ei_N$ is 
isometric to $\sqrt{2}E_8$. 
Let $N'$  be the $(-1)$-eigenlattice of $i_Ei_N$. 
Then $\la i_M, i_Ni_E\ra \cong \mathrm{Dih}_{12}$ and hence 
$M+N'\cong M+N+E \cong \dih{12}{16}$ as defined in \cite{GL}.
Let $e'$ be an Ising vector of $E_8$-type in $V_{N'}^+$ 
(see Remark \ref{EE8Ising} for definition). 
Then $\la t_M, e'\ra \cong U_{\mathrm{6A}}$ since $\tau_{t_M}\tau_{e'}$ has order 6.

Let $\tilde{e}$ be the central Ising vector in $\la t_M, e'\ra \cong U_{\mathrm{6A}}$ 
(cf.~Eq.~\eqref{eq:2.9}). 
Then $\tilde{e}\in V_E^+$. 
Let $\tilde{t}_{N} = e'\circ t_N$. 
Then $\tilde{t}_N \in V_N^+$ and
\begin{equation}\label{te}
  (t_M \mymid \tilde{e})
  = (\tilde{t}_N \mymid \tilde{e})
  = \dfrac{1}{32}\quad \text{ and } \quad (t_M  \mymid \tilde{t}_N)=\dfrac{13}{2^{10}}.
\end{equation}
By the classification of Ising vectors in \cite{LSY} and \cite{LS} 
(cf.~Remark \ref{EE8Ising}), 
$\tilde{t}_N = \varphi_x (t_N)$ for some $x\in N^*$.
Moreover, by Proposition 3.4 of \cite{GL3}, there exists an $\al\in M$ such that 
$\varphi_\al|_{V_N^+}= \varphi_x$.
Hence, we have
\[
  \varphi_\al (t_M)=t_M\quad 
  \text{and} \quad 
  \varphi_\al( \tilde{t}_N) 
  = \varphi_\al \varphi_x (t_N) =t_N.
\]
Now let $e=\varphi_\al(\tilde{e})$. 
Then we have
\[
  (t_M \mymid e)=(t_N \mymid e)=\dfrac{1}{32}
\]
by Equation \eqref{te}.
\qed

Finally, by Lemma \ref{abx}, there exist Ising vectors $x\in V_X^+$ and 
$y\in V_Y^+$ such that 
$\la a,x\ra \cong \la b, x\ra\cong \la a,y\ra \cong \la b, y\ra\cong U_{\mathrm{2A}}$. 
Since $X+Y\cong \dih{6}{14}$, we have $(x \mymid y)= 13\cd 2^{-10}$ or $5\cd 2^{-10}$. 
By Theorem \ref{thm:3.5}, the case $(x \mymid y)= 5\cd 2^{-10}$ is impossible.
Hence we have $(x \mymid y)= 13\cd 2^{-10}$ and $\la x,y\ra \cong U_{\mathrm{3A}}$ as desired.

\end{appendices}

\small

\end{document}

%% file: fig1.tex
\unitlength 0.1in
\begin{picture}( 40.0000,  8.0000)(  2.0000,-10.0000)
%

%

%

%
\special{pn 13}%
\special{pa 2600 200}%
\special{pa 2600 1000}%
\special{fp}%
\special{pa 2600 1000}%
\special{pa 3800 1000}%
\special{fp}%
\special{pa 3800 1000}%
\special{pa 3800 200}%
\special{fp}%
\special{pa 3800 200}%
\special{pa 2600 200}%
\special{fp}%
\special{pa 2600 600}%
\special{pa 3800 600}%
\special{fp}%
\special{pa 3000 1000}%
\special{pa 3000 200}%
\special{fp}%
\special{pa 3400 200}%
\special{pa 3400 1000}%
\special{fp}%
%
\special{pn 20}%
\special{sh 1}%
\special{ar 2700 300 10 10 0  6.28318530717959E+0000}%
\special{sh 1}%
\special{ar 2900 300 10 10 0  6.28318530717959E+0000}%
\special{sh 1}%
\special{ar 3100 300 10 10 0  6.28318530717959E+0000}%
\special{sh 1}%
\special{ar 3300 300 10 10 0  6.28318530717959E+0000}%
\special{sh 1}%
\special{ar 3500 300 10 10 0  6.28318530717959E+0000}%
\special{sh 1}%
\special{ar 3700 300 10 10 0  6.28318530717959E+0000}%
\special{sh 1}%
\special{ar 3700 500 10 10 0  6.28318530717959E+0000}%
\special{sh 1}%
\special{ar 3500 500 10 10 0  6.28318530717959E+0000}%
\special{sh 1}%
\special{ar 3300 500 10 10 0  6.28318530717959E+0000}%
\special{sh 1}%
\special{ar 3100 500 10 10 0  6.28318530717959E+0000}%
\special{sh 1}%
\special{ar 2900 500 10 10 0  6.28318530717959E+0000}%
\special{sh 1}%
\special{ar 2700 500 10 10 0  6.28318530717959E+0000}%
\special{sh 1}%
\special{ar 2700 700 10 10 0  6.28318530717959E+0000}%
\special{sh 1}%
\special{ar 2900 700 10 10 0  6.28318530717959E+0000}%
\special{sh 1}%
\special{ar 2900 900 10 10 0  6.28318530717959E+0000}%
\special{sh 1}%
\special{ar 2700 900 10 10 0  6.28318530717959E+0000}%
\special{sh 1}%
\special{ar 3100 900 10 10 0  6.28318530717959E+0000}%
\special{sh 1}%
\special{ar 3100 700 10 10 0  6.28318530717959E+0000}%
\special{sh 1}%
\special{ar 3300 700 10 10 0  6.28318530717959E+0000}%
\special{sh 1}%
\special{ar 3300 900 10 10 0  6.28318530717959E+0000}%
\special{sh 1}%
\special{ar 3500 900 10 10 0  6.28318530717959E+0000}%
\special{sh 1}%
\special{ar 3700 900 10 10 0  6.28318530717959E+0000}%
\special{sh 1}%
\special{ar 3700 700 10 10 0  6.28318530717959E+0000}%
\special{sh 1}%
\special{ar 3500 700 10 10 0  6.28318530717959E+0000}%
\special{sh 1}%
\special{ar 3500 700 10 10 0  6.28318530717959E+0000}%
%
\special{pn 13}%
\special{pa 2700 300}%
\special{pa 2900 300}%
\special{fp}%
\special{pa 2700 500}%
\special{pa 2700 700}%
\special{fp}%
\special{pa 3100 900}%
\special{pa 3300 300}%
\special{fp}%
\special{pa 3100 500}%
\special{pa 3700 300}%
\special{fp}%
\special{pa 3300 700}%
\special{pa 3500 900}%
\special{fp}%
\special{pa 3500 500}%
\special{pa 3700 700}%
\special{fp}%
%
\special{pn 13}%
\special{pa 3100 700}%
\special{pa 3130 688}%
\special{pa 3162 674}%
\special{pa 3192 660}%
\special{pa 3220 646}%
\special{pa 3248 630}%
\special{pa 3276 614}%
\special{pa 3302 596}%
\special{pa 3328 576}%
\special{pa 3350 556}%
\special{pa 3372 534}%
\special{pa 3390 510}%
\special{pa 3408 484}%
\special{pa 3426 456}%
\special{pa 3442 428}%
\special{pa 3456 400}%
\special{pa 3470 370}%
\special{pa 3484 340}%
\special{pa 3496 308}%
\special{pa 3500 300}%
\special{sp}%
%
\special{pn 13}%
\special{pa 3300 500}%
\special{pa 3316 530}%
\special{pa 3332 558}%
\special{pa 3348 586}%
\special{pa 3364 614}%
\special{pa 3382 642}%
\special{pa 3400 668}%
\special{pa 3418 694}%
\special{pa 3438 718}%
\special{pa 3460 740}%
\special{pa 3484 762}%
\special{pa 3508 782}%
\special{pa 3532 802}%
\special{pa 3558 820}%
\special{pa 3586 836}%
\special{pa 3614 852}%
\special{pa 3642 868}%
\special{pa 3670 884}%
\special{pa 3700 900}%
\special{pa 3700 900}%
\special{sp}%
%
\special{pn 13}%
\special{pa 600 200}%
\special{pa 600 1000}%
\special{fp}%
\special{pa 600 1000}%
\special{pa 1800 1000}%
\special{fp}%
\special{pa 1800 1000}%
\special{pa 1800 200}%
\special{fp}%
\special{pa 1800 200}%
\special{pa 600 200}%
\special{fp}%
\special{pa 600 600}%
\special{pa 1800 600}%
\special{fp}%
\special{pa 1000 1000}%
\special{pa 1000 200}%
\special{fp}%
\special{pa 1400 200}%
\special{pa 1400 1000}%
\special{fp}%
%
\special{pn 20}%
\special{sh 1}%
\special{ar 700 300 10 10 0  6.28318530717959E+0000}%
\special{sh 1}%
\special{ar 900 300 10 10 0  6.28318530717959E+0000}%
\special{sh 1}%
\special{ar 1100 300 10 10 0  6.28318530717959E+0000}%
\special{sh 1}%
\special{ar 1300 300 10 10 0  6.28318530717959E+0000}%
\special{sh 1}%
\special{ar 1500 300 10 10 0  6.28318530717959E+0000}%
\special{sh 1}%
\special{ar 1700 300 10 10 0  6.28318530717959E+0000}%
\special{sh 1}%
\special{ar 1700 500 10 10 0  6.28318530717959E+0000}%
\special{sh 1}%
\special{ar 1500 500 10 10 0  6.28318530717959E+0000}%
\special{sh 1}%
\special{ar 1300 500 10 10 0  6.28318530717959E+0000}%
\special{sh 1}%
\special{ar 1100 500 10 10 0  6.28318530717959E+0000}%
\special{sh 1}%
\special{ar 900 500 10 10 0  6.28318530717959E+0000}%
\special{sh 1}%
\special{ar 700 500 10 10 0  6.28318530717959E+0000}%
\special{sh 1}%
\special{ar 700 700 10 10 0  6.28318530717959E+0000}%
\special{sh 1}%
\special{ar 900 700 10 10 0  6.28318530717959E+0000}%
\special{sh 1}%
\special{ar 900 900 10 10 0  6.28318530717959E+0000}%
\special{sh 1}%
\special{ar 700 900 10 10 0  6.28318530717959E+0000}%
\special{sh 1}%
\special{ar 1100 900 10 10 0  6.28318530717959E+0000}%
\special{sh 1}%
\special{ar 1100 700 10 10 0  6.28318530717959E+0000}%
\special{sh 1}%
\special{ar 1300 700 10 10 0  6.28318530717959E+0000}%
\special{sh 1}%
\special{ar 1300 900 10 10 0  6.28318530717959E+0000}%
\special{sh 1}%
\special{ar 1500 900 10 10 0  6.28318530717959E+0000}%
\special{sh 1}%
\special{ar 1700 900 10 10 0  6.28318530717959E+0000}%
\special{sh 1}%
\special{ar 1700 700 10 10 0  6.28318530717959E+0000}%
\special{sh 1}%
\special{ar 1500 700 10 10 0  6.28318530717959E+0000}%
\special{sh 1}%
\special{ar 1500 700 10 10 0  6.28318530717959E+0000}%
%
\special{pn 13}%
\special{pa 700 300}%
\special{pa 700 300}%
\special{fp}%
\special{pa 700 300}%
\special{pa 900 300}%
\special{fp}%
\special{pa 700 700}%
\special{pa 700 900}%
\special{fp}%
\special{pa 1100 500}%
\special{pa 1300 300}%
\special{fp}%
\special{pa 1500 700}%
\special{pa 1700 900}%
\special{fp}%
\special{pa 1300 700}%
\special{pa 1700 500}%
\special{fp}%
\special{pa 1500 500}%
\special{pa 1300 900}%
\special{fp}%
\special{pa 1100 700}%
\special{pa 1700 300}%
\special{fp}%
\special{pa 1500 300}%
\special{pa 1100 900}%
\special{fp}%
\put(2.0000,-7.0000){\makebox(0,0)[lb]{$i_a=$}}%
\put(23.0000,-6.0000){\makebox(0,0){$i_b=$}}%
\end{picture}%

%% file: fig2.tex
\unitlength 0.1in
\begin{picture}( 40.0000,  8.0000)(  2.0000,-10.0000)
%

%

%

%
\special{pn 13}%
\special{pa 2600 200}%
\special{pa 2600 1000}%
\special{fp}%
\special{pa 2600 1000}%
\special{pa 3800 1000}%
\special{fp}%
\special{pa 3800 1000}%
\special{pa 3800 200}%
\special{fp}%
\special{pa 3800 200}%
\special{pa 2600 200}%
\special{fp}%
\special{pa 2600 600}%
\special{pa 3800 600}%
\special{fp}%
\special{pa 3000 1000}%
\special{pa 3000 200}%
\special{fp}%
\special{pa 3400 200}%
\special{pa 3400 1000}%
\special{fp}%
%
\special{pn 20}%
\special{sh 1}%
\special{ar 2700 300 10 10 0  6.28318530717959E+0000}%
\special{sh 1}%
\special{ar 2900 300 10 10 0  6.28318530717959E+0000}%
\special{sh 1}%
\special{ar 3100 300 10 10 0  6.28318530717959E+0000}%
\special{sh 1}%
\special{ar 3300 300 10 10 0  6.28318530717959E+0000}%
\special{sh 1}%
\special{ar 3500 300 10 10 0  6.28318530717959E+0000}%
\special{sh 1}%
\special{ar 3700 300 10 10 0  6.28318530717959E+0000}%
\special{sh 1}%
\special{ar 3700 500 10 10 0  6.28318530717959E+0000}%
\special{sh 1}%
\special{ar 3500 500 10 10 0  6.28318530717959E+0000}%
\special{sh 1}%
\special{ar 3300 500 10 10 0  6.28318530717959E+0000}%
\special{sh 1}%
\special{ar 3100 500 10 10 0  6.28318530717959E+0000}%
\special{sh 1}%
\special{ar 2900 500 10 10 0  6.28318530717959E+0000}%
\special{sh 1}%
\special{ar 2700 500 10 10 0  6.28318530717959E+0000}%
\special{sh 1}%
\special{ar 2700 700 10 10 0  6.28318530717959E+0000}%
\special{sh 1}%
\special{ar 2900 700 10 10 0  6.28318530717959E+0000}%
\special{sh 1}%
\special{ar 2900 900 10 10 0  6.28318530717959E+0000}%
\special{sh 1}%
\special{ar 2700 900 10 10 0  6.28318530717959E+0000}%
\special{sh 1}%
\special{ar 3100 900 10 10 0  6.28318530717959E+0000}%
\special{sh 1}%
\special{ar 3100 700 10 10 0  6.28318530717959E+0000}%
\special{sh 1}%
\special{ar 3300 700 10 10 0  6.28318530717959E+0000}%
\special{sh 1}%
\special{ar 3300 900 10 10 0  6.28318530717959E+0000}%
\special{sh 1}%
\special{ar 3500 900 10 10 0  6.28318530717959E+0000}%
\special{sh 1}%
\special{ar 3700 900 10 10 0  6.28318530717959E+0000}%
\special{sh 1}%
\special{ar 3700 700 10 10 0  6.28318530717959E+0000}%
\special{sh 1}%
\special{ar 3500 700 10 10 0  6.28318530717959E+0000}%
\special{sh 1}%
\special{ar 3500 700 10 10 0  6.28318530717959E+0000}%
%
\special{pn 13}%
\special{pa 2700 300}%
\special{pa 2900 300}%
\special{fp}%
\special{pa 2900 500}%
\special{pa 2900 700}%
\special{fp}%
\special{pa 3100 900}%
\special{pa 3300 300}%
\special{fp}%
\special{pa 3100 500}%
\special{pa 3500 300}%
\special{fp}%
\special{pa 3300 500}%
\special{pa 3500 900}%
\special{fp}%
\special{pa 3300 700}%
\special{pa 3700 900}%
\special{fp}%
%
\special{pn 13}%
\special{pa 3100 700}%
\special{pa 3700 300}%
\special{fp}%
%
\special{pn 13}%
\special{pa 3500 700}%
\special{pa 3700 500}%
\special{fp}%
%
\special{pn 13}%
\special{pa 600 200}%
\special{pa 600 1000}%
\special{fp}%
\special{pa 600 1000}%
\special{pa 1800 1000}%
\special{fp}%
\special{pa 1800 1000}%
\special{pa 1800 200}%
\special{fp}%
\special{pa 1800 200}%
\special{pa 600 200}%
\special{fp}%
\special{pa 600 600}%
\special{pa 1800 600}%
\special{fp}%
\special{pa 1000 1000}%
\special{pa 1000 200}%
\special{fp}%
\special{pa 1400 200}%
\special{pa 1400 1000}%
\special{fp}%
%
\special{pn 20}%
\special{sh 1}%
\special{ar 700 300 10 10 0  6.28318530717959E+0000}%
\special{sh 1}%
\special{ar 900 300 10 10 0  6.28318530717959E+0000}%
\special{sh 1}%
\special{ar 1100 300 10 10 0  6.28318530717959E+0000}%
\special{sh 1}%
\special{ar 1300 300 10 10 0  6.28318530717959E+0000}%
\special{sh 1}%
\special{ar 1500 300 10 10 0  6.28318530717959E+0000}%
\special{sh 1}%
\special{ar 1700 300 10 10 0  6.28318530717959E+0000}%
\special{sh 1}%
\special{ar 1700 500 10 10 0  6.28318530717959E+0000}%
\special{sh 1}%
\special{ar 1500 500 10 10 0  6.28318530717959E+0000}%
\special{sh 1}%
\special{ar 1300 500 10 10 0  6.28318530717959E+0000}%
\special{sh 1}%
\special{ar 1100 500 10 10 0  6.28318530717959E+0000}%
\special{sh 1}%
\special{ar 900 500 10 10 0  6.28318530717959E+0000}%
\special{sh 1}%
\special{ar 700 500 10 10 0  6.28318530717959E+0000}%
\special{sh 1}%
\special{ar 700 700 10 10 0  6.28318530717959E+0000}%
\special{sh 1}%
\special{ar 900 700 10 10 0  6.28318530717959E+0000}%
\special{sh 1}%
\special{ar 900 900 10 10 0  6.28318530717959E+0000}%
\special{sh 1}%
\special{ar 700 900 10 10 0  6.28318530717959E+0000}%
\special{sh 1}%
\special{ar 1100 900 10 10 0  6.28318530717959E+0000}%
\special{sh 1}%
\special{ar 1100 700 10 10 0  6.28318530717959E+0000}%
\special{sh 1}%
\special{ar 1300 700 10 10 0  6.28318530717959E+0000}%
\special{sh 1}%
\special{ar 1300 900 10 10 0  6.28318530717959E+0000}%
\special{sh 1}%
\special{ar 1500 900 10 10 0  6.28318530717959E+0000}%
\special{sh 1}%
\special{ar 1700 900 10 10 0  6.28318530717959E+0000}%
\special{sh 1}%
\special{ar 1700 700 10 10 0  6.28318530717959E+0000}%
\special{sh 1}%
\special{ar 1500 700 10 10 0  6.28318530717959E+0000}%
\special{sh 1}%
\special{ar 1500 700 10 10 0  6.28318530717959E+0000}%
%
\special{pn 13}%
\special{pa 700 300}%
\special{pa 700 300}%
\special{fp}%
\special{pa 700 300}%
\special{pa 900 300}%
\special{fp}%
\special{pa 900 700}%
\special{pa 900 900}%
\special{fp}%
\special{pa 1100 500}%
\special{pa 1300 300}%
\special{fp}%

\special{pn 13}%
\special{pa 1100 700}%
\special{pa 1130 688}%
\special{pa 1162 674}%
\special{pa 1192 660}%
\special{pa 1220 646}%
\special{pa 1248 630}%
\special{pa 1276 614}%
\special{pa 1302 596}%
\special{pa 1328 576}%
\special{pa 1350 556}%
\special{pa 1372 534}%
\special{pa 1390 510}%
\special{pa 1408 484}%
\special{pa 1426 456}%
\special{pa 1442 428}%
\special{pa 1456 400}%
\special{pa 1470 370}%
\special{pa 1484 340}%
\special{pa 1496 308}%
\special{pa 1500 300}%
\special{sp}%

\special{pn 13}%
\special{pa 1300 900}%
\special{pa 1330 888}%
\special{pa 1362 874}%
\special{pa 1392 860}%
\special{pa 1420 846}%
\special{pa 1448 830}%
\special{pa 1476 814}%
\special{pa 1502 796}%
\special{pa 1528 776}%
\special{pa 1550 756}%
\special{pa 1572 734}%
\special{pa 1590 710}%
\special{pa 1608 684}%
\special{pa 1626 656}%
\special{pa 1642 628}%
\special{pa 1656 600}%
\special{pa 1670 570}%
\special{pa 1684 540}%
\special{pa 1696 508}%
\special{pa 1700 500}%
\special{sp}%

\special{pn 13}%
\special{pa 1500 900}%
\special{pa 1700 700}%
\special{fp}%
\special{fp}%
\special{pa 1500 500}%
\special{pa 1300 700}%
\special{fp}%

\special{pn 13}%
\special{pa 1100 900}%
\special{pa 1130 888}%
\special{pa 1162 874}%
\special{pa 1192 860}%
\special{pa 1220 842}%
\special{pa 1400 705}%
\special{pa 1420 675}%
\special{pa 1448 650}%
\special{pa 1642 420}%
\special{pa 1656 400}%
\special{pa 1670 370}%
\special{pa 1684 340}%
\special{pa 1696 308}%
\special{pa 1700 300}%
\special{sp}%
\put(2.0000,-7.0000){\makebox(0,0)[lb]{$i_x=$}}%
\put(23.0000,-6.0000){\makebox(0,0){$i_y=$}}%
\end{picture}%

%% file: fig3.tex
\unitlength 0.1in
\begin{picture}( 44.9500,  8.0000)(  1.0500,-10.0000)
%
\special{pn 13}%
\special{pa 3400 200}%
\special{pa 3400 1000}%
\special{fp}%
\special{pa 3400 1000}%
\special{pa 4600 1000}%
\special{fp}%
\special{pa 4600 1000}%
\special{pa 4600 200}%
\special{fp}%
\special{pa 4600 200}%
\special{pa 3400 200}%
\special{fp}%
\special{pa 3400 600}%
\special{pa 4600 600}%
\special{fp}%
\special{pa 3800 1000}%
\special{pa 3800 200}%
\special{fp}%
\special{pa 4200 200}%
\special{pa 4200 1000}%
\special{fp}%
%
\special{pn 20}%
\special{sh 1}%
\special{ar 3700 300 10 10 0  6.28318530717959E+0000}%
\special{sh 1}%
\special{ar 3500 300 10 10 0  6.28318530717959E+0000}%
\special{sh 1}%
\special{ar 3900 300 10 10 0  6.28318530717959E+0000}%
\special{sh 1}%
\special{ar 4100 300 10 10 0  6.28318530717959E+0000}%
\special{sh 1}%
\special{ar 4300 300 10 10 0  6.28318530717959E+0000}%
\special{sh 1}%
\special{ar 4500 300 10 10 0  6.28318530717959E+0000}%
\special{sh 1}%
\special{ar 4500 500 10 10 0  6.28318530717959E+0000}%
\special{sh 1}%
\special{ar 4300 500 10 10 0  6.28318530717959E+0000}%
\special{sh 1}%
\special{ar 4100 500 10 10 0  6.28318530717959E+0000}%
\special{sh 1}%
\special{ar 3900 500 10 10 0  6.28318530717959E+0000}%
\special{sh 1}%
\special{ar 3700 500 10 10 0  6.28318530717959E+0000}%
\special{sh 1}%
\special{ar 3500 500 10 10 0  6.28318530717959E+0000}%
\special{sh 1}%
\special{ar 3500 700 10 10 0  6.28318530717959E+0000}%
\special{sh 1}%
\special{ar 3700 700 10 10 0  6.28318530717959E+0000}%
\special{sh 1}%
\special{ar 3700 900 10 10 0  6.28318530717959E+0000}%
\special{sh 1}%
\special{ar 3500 900 10 10 0  6.28318530717959E+0000}%
\special{sh 1}%
\special{ar 3900 900 10 10 0  6.28318530717959E+0000}%
\special{sh 1}%
\special{ar 3900 700 10 10 0  6.28318530717959E+0000}%
\special{sh 1}%
\special{ar 4100 700 10 10 0  6.28318530717959E+0000}%
\special{sh 1}%
\special{ar 4100 900 10 10 0  6.28318530717959E+0000}%
\special{sh 1}%
\special{ar 4300 900 10 10 0  6.28318530717959E+0000}%
\special{sh 1}%
\special{ar 4500 900 10 10 0  6.28318530717959E+0000}%
\special{sh 1}%
\special{ar 4500 700 10 10 0  6.28318530717959E+0000}%
\special{sh 1}%
\special{ar 4300 700 10 10 0  6.28318530717959E+0000}%
\special{sh 1}%
\special{ar 4300 700 10 10 0  6.28318530717959E+0000}%

\special{pn 13}%
\special{pa 3700 900}%
\special{pa 3700 900}%
\special{fp}%
%
\special{pn 13}%
\special{pa 3700 900}%
\special{pa 3700 716}%
\special{fp}%
\special{sh 1}%
\special{pa 3700 716}%
\special{pa 3680 782}%
\special{pa 3700 768}%
\special{pa 3720 782}%
\special{pa 3700 716}%
\special{fp}%
\special{pa 3700 700}%
\special{pa 3700 516}%
\special{fp}%
\special{sh 1}%
\special{pa 3700 516}%
\special{pa 3680 582}%
\special{pa 3700 568}%
\special{pa 3720 582}%
\special{pa 3700 516}%
\special{fp}%
%
\special{pn 13}%
\special{pa 3900 500}%
\special{pa 3900 686}%
\special{fp}%
\special{sh 1}%
\special{pa 3900 686}%
\special{pa 3920 618}%
\special{pa 3900 632}%
\special{pa 3880 618}%
\special{pa 3900 686}%
\special{fp}%
\special{pa 3900 700}%
\special{pa 3900 886}%
\special{fp}%
\special{sh 1}%
\special{pa 3900 886}%
\special{pa 3920 818}%
\special{pa 3900 832}%
\special{pa 3880 818}%
\special{pa 3900 886}%
\special{fp}%
%
\special{pn 13}%
\special{pa 4100 300}%
\special{pa 4286 300}%
\special{fp}%
\special{sh 1}%
\special{pa 4286 300}%
\special{pa 4218 280}%
\special{pa 4232 300}%
\special{pa 4218 320}%
\special{pa 4286 300}%
\special{fp}%
\special{pa 4300 300}%
\special{pa 4486 300}%
\special{fp}%
\special{sh 1}%
\special{pa 4486 300}%
\special{pa 4418 280}%
\special{pa 4432 300}%
\special{pa 4418 320}%
\special{pa 4486 300}%
\special{fp}%
%
\special{pn 13}%
\special{pa 4100 500}%
\special{pa 4490 690}%
\special{fp}%
\special{sh 1}%
\special{pa 4490 690}%
\special{pa 4458 620}%
\special{pa 4450 675}%
\special{pa 4410 680}%
\special{pa 4490 690}%
\special{fp}%
\special{pa 4500 700}%
\special{pa 4300 900}%
\special{fp}%
\special{sh 1}%
\special{pa 4310 890}%
\special{pa 4380 860}%
\special{pa 4358 830}%
\special{pa 4362 800}%
\special{pa 4310 890}%
\special{fp}%
%
\special{pn 13}%
\special{pa 4300 500}%
\special{pa 4500 900}%
\special{fp}%
\special{sh 1}%
\special{pa 4490 890}%
\special{pa 4480 820}%
\special{pa 4460 860}%
\special{pa 4430 845}%
\special{pa 4490 890}%
\special{fp}%
\special{pa 4500 900}%
\special{pa 4100 700}%
\special{fp}%
\special{sh 1}%
\special{pa 4100 710}%
\special{pa 4160 690}%
\special{pa 4130 720}%
\special{pa 4120 750}%
\special{pa 4100 710}%
\special{fp}%
%
\special{pn 13}%
\special{pa 4500 500}%
\special{pa 4300 700}%
\special{fp}%

\special{pn 13}%
\special{pa 4100 900}%
\special{pa 4130 888}%
\special{pa 4162 874}%
\special{pa 4192 860}%
\special{pa 4220 846}%
\special{pa 4248 830}%
\special{pa 4276 814}%
\special{pa 4302 796}%
\special{pa 4328 776}%
\special{pa 4350 756}%
\special{pa 4372 734}%
\special{pa 4390 710}%
\special{pa 4408 684}%
\special{pa 4426 656}%
\special{pa 4442 628}%
\special{pa 4456 600}%
\special{pa 4470 570}%
\special{pa 4484 540}%
\special{pa 4496 508}%
\special{pa 4500 500}%
\special{sp}%
\special{sh 1}%
\special{pa 4100 900}%
\special{pa 4160 830}%
\special{pa 4130 875}%
\special{pa 4180 890}%
\special{pa 4100 900}%
\special{fp}%

%
\special{pn 13}%
\special{pa 1000 200}%
\special{pa 1000 1000}%
\special{fp}%
\special{pa 1000 1000}%
\special{pa 2200 1000}%
\special{fp}%
\special{pa 2200 1000}%
\special{pa 2200 200}%
\special{fp}%
\special{pa 2200 200}%
\special{pa 1000 200}%
\special{fp}%
\special{pa 1000 600}%
\special{pa 2200 600}%
\special{fp}%
\special{pa 1400 1000}%
\special{pa 1400 200}%
\special{fp}%
\special{pa 1800 200}%
\special{pa 1800 1000}%
\special{fp}%
%
\special{pn 20}%
\special{sh 1}%
\special{ar 1100 300 10 10 0  6.28318530717959E+0000}%
\special{sh 1}%
\special{ar 1300 300 10 10 0  6.28318530717959E+0000}%
\special{sh 1}%
\special{ar 1500 300 10 10 0  6.28318530717959E+0000}%
\special{sh 1}%
\special{ar 1700 300 10 10 0  6.28318530717959E+0000}%
\special{sh 1}%
\special{ar 1900 300 10 10 0  6.28318530717959E+0000}%
\special{sh 1}%
\special{ar 2100 300 10 10 0  6.28318530717959E+0000}%
\special{sh 1}%
\special{ar 2100 500 10 10 0  6.28318530717959E+0000}%
\special{sh 1}%
\special{ar 1900 500 10 10 0  6.28318530717959E+0000}%
\special{sh 1}%
\special{ar 1700 500 10 10 0  6.28318530717959E+0000}%
\special{sh 1}%
\special{ar 1500 500 10 10 0  6.28318530717959E+0000}%
\special{sh 1}%
\special{ar 1300 500 10 10 0  6.28318530717959E+0000}%
\special{sh 1}%
\special{ar 1100 500 10 10 0  6.28318530717959E+0000}%
\special{sh 1}%
\special{ar 1100 700 10 10 0  6.28318530717959E+0000}%
\special{sh 1}%
\special{ar 1300 700 10 10 0  6.28318530717959E+0000}%
\special{sh 1}%
\special{ar 1300 900 10 10 0  6.28318530717959E+0000}%
\special{sh 1}%
\special{ar 1100 900 10 10 0  6.28318530717959E+0000}%
\special{sh 1}%
\special{ar 1500 900 10 10 0  6.28318530717959E+0000}%
\special{sh 1}%
\special{ar 1500 700 10 10 0  6.28318530717959E+0000}%
\special{sh 1}%
\special{ar 1700 700 10 10 0  6.28318530717959E+0000}%
\special{sh 1}%
\special{ar 1700 900 10 10 0  6.28318530717959E+0000}%
\special{sh 1}%
\special{ar 1900 900 10 10 0  6.28318530717959E+0000}%
\special{sh 1}%
\special{ar 2100 900 10 10 0  6.28318530717959E+0000}%
\special{sh 1}%
\special{ar 2100 700 10 10 0  6.28318530717959E+0000}%
\special{sh 1}%
\special{ar 1900 700 10 10 0  6.28318530717959E+0000}%
\special{sh 1}%
\special{ar 1900 700 10 10 0  6.28318530717959E+0000}%
%
\special{pn 13}%
\special{pa 1100 900}%
\special{pa 1100 900}%
\special{fp}%
%
\special{pn 13}%
\special{pa 1100 900}%
\special{pa 1100 716}%
\special{fp}%
\special{sh 1}%
\special{pa 1100 716}%
\special{pa 1080 782}%
\special{pa 1100 768}%
\special{pa 1120 782}%
\special{pa 1100 716}%
\special{fp}%
\special{pa 1100 700}%
\special{pa 1100 516}%
\special{fp}%
\special{sh 1}%
\special{pa 1100 516}%
\special{pa 1080 582}%
\special{pa 1100 568}%
\special{pa 1120 582}%
\special{pa 1100 516}%
\special{fp}%
%
\special{pn 13}%
\special{pa 1500 500}%
\special{pa 1500 686}%
\special{fp}%
\special{sh 1}%
\special{pa 1500 686}%
\special{pa 1520 618}%
\special{pa 1500 632}%
\special{pa 1480 618}%
\special{pa 1500 686}%
\special{fp}%
\special{pa 1500 700}%
\special{pa 1500 886}%
\special{fp}%
\special{sh 1}%
\special{pa 1500 886}%
\special{pa 1520 818}%
\special{pa 1500 832}%
\special{pa 1480 818}%
\special{pa 1500 886}%
\special{fp}%
%
\special{pn 13}%
\special{pa 1700 300}%
\special{pa 1886 300}%
\special{fp}%
\special{sh 1}%
\special{pa 1886 300}%
\special{pa 1818 280}%
\special{pa 1832 300}%
\special{pa 1818 320}%
\special{pa 1886 300}%
\special{fp}%
\special{pa 1900 300}%
\special{pa 2086 300}%
\special{fp}%
\special{sh 1}%
\special{pa 2086 300}%
\special{pa 2018 280}%
\special{pa 2032 300}%
\special{pa 2018 320}%
\special{pa 2086 300}%
\special{fp}%
%
\special{pn 13}%
\special{pa 1900 500}%
\special{pa 2090 690}%
\special{fp}%
\special{sh 1}%
\special{pa 2090 690}%
\special{pa 2058 630}%
\special{pa 2052 652}%
\special{pa 2030 658}%
\special{pa 2090 690}%
\special{fp}%
\special{pa 2100 700}%
\special{pa 1710 890}%
\special{fp}%
\special{sh 1}%
\special{pa 1710 890}%
\special{pa 1780 880}%
\special{pa 1758 868}%
\special{pa 1762 844}%
\special{pa 1710 890}%
\special{fp}%
%
\special{pn 13}%
\special{pa 1700 700}%
\special{pa 1890 890}%
\special{fp}%
\special{sh 1}%
\special{pa 1890 890}%
\special{pa 1858 830}%
\special{pa 1852 852}%
\special{pa 1830 858}%
\special{pa 1890 890}%
\special{fp}%
\special{pa 1900 900}%
\special{pa 2090 510}%
\special{fp}%
\special{sh 1}%
\special{pa 2090 510}%
\special{pa 2044 562}%
\special{pa 2068 558}%
\special{pa 2080 580}%
\special{pa 2090 510}%
\special{fp}%
%
\special{pn 13}%
\special{pa 1700 500}%
\special{pa 1890 690}%
\special{fp}%
\special{sh 1}%
\special{pa 1890 690}%
\special{pa 1858 630}%
\special{pa 1852 652}%
\special{pa 1830 658}%
\special{pa 1890 690}%
\special{fp}%
\special{pa 1900 700}%
\special{pa 2090 890}%
\special{fp}%
\special{sh 1}%
\special{pa 2090 890}%
\special{pa 2058 830}%
\special{pa 2052 852}%
\special{pa 2030 858}%
\special{pa 2090 890}%
\special{fp}%
\put(6.0000,-6.0000){\makebox(0,0){$g=i_ai_b=$}}%
\put(30.0000,-6.0000){\makebox(0,0){$h=i_xi_y=$}}%
\end{picture}%

%% file: fi23_R2_rev2.bbl
\begin{thebibliography}{AAAA99}

\bibitem[ATLAS]{ATLAS}
  J.H. Conway, R.T. Curtis, S.P. Norton, R.A. Parker and  R.A. Wilson,
  ATLAS of finite groups. 
  Clarendon Press, Oxford, 1985.

\bibitem[As97]{As}
  M. Aschbacher, 
  $3$-transposition groups.
  Cambridge Tracts in Mathematics, 124. Cambridge University Press, Cambridge, 1997.

\bibitem[B86]{B}
  R.E. Borcherds, 
  Vertex algebras, Kac-Moody algebras and the Monster.
  \textit{Proc.\ Nat.\ Acad.\ Sci. USA} \textbf{83} (1986), 3068--3071.

\bibitem[C85]{C}
  J.H. Conway, 
  A simple construction for the Fischer-Griess Monster group. 
  \textit{Invent. Math.} \textbf{79} (1985), 513--540.

\bibitem[CS99]{cs}
  J.\,H.\,Conway and N.\,J.\,A. Sloane, 
  Sphere packings, lattices and groups. 
  3rd Edition, Springer, New York, 1999.


\bibitem[CH95]{CH}
  H. Cuypers and J.I. Hall, 
  The $3$-transposition groups with trivial center.
  \textit{J. Algebra} \textbf{178}, (1995), 149--193.

\bibitem[CLY]{CLY}
  T. Creutzig, C.H. Lam and H. Yamauchi, An exceptional construction of the moonshine 
  vertex operator algebra, in preparation.







\bibitem[DJY19]{DJY}
  C. Dong, X. Jiao and N. Yu, 
  6A-Algebra and its representations.
  \href{https://arxiv.org/abs/1902.06951}{arXiv:1902.06951}.

\bibitem[DLMN98]{DLMN}
  C. Dong, H. Li, G. Mason and S.P. Norton, 
  Associative subalgebras of Griess algebra and related topics.
  Proc. of the Conference on the Monster and Lie algebra at the Ohio State
  University, May 1996, ed. by J. Ferrar and K. Harada, Walter de
  Gruyter, Berlin - New York, 1998.

\bibitem[DLY09]{DLY}
  C. Dong. C.H. Lam and H. Yamada,
  W-algebras related to parafermion algebras.
  \textit{J. Algebra} \textbf{322} (2009), no.7, 2366--2403.


\bibitem[DMZ94]{DMZ}
  C. Dong, G. Mason and Y. Zhu, Discrete series of the Virasoro algebra and
  the moonshine module.
  Proc. Symp. Pure. Math., American Math. Soc. \textbf{56} II (1994), 295--316.


\bibitem[DW16]{DW}
  C. Dong and Q. Wang,  
  Quantum dimensions and fusion rules for parafermion vertex operator algebras.
  \textit{Proc. Amer. Math. Soc.} \textbf{144} (2016), no. 4, 1483--1492. 





\bibitem[Fi71]{Fi}
  B. Fischer, 
  Finite groups generated by 3-transpositions. I.
  \textit{Invent. Math.} \textbf{13}, 232--246.

\bibitem[FLM88]{FLM}
  I.B. Frenkel, J. Lepowsky and A. Meurman, 
  Vertex Operator Algebras and the Monster. 
  Academic Press, New York, 1988.


\bibitem[FZ92]{FZ}
  I.B. Frenkel and Y. Zhu, 
  Vertex operator algebras associated to representation of affine and Virasoro algebras. 
  \textit{Duke Math. J.} \textbf{66} (1992), 123--168.

\bibitem[GKO86]{GKO}
  P. Goddard, A. Kent and D. Olive, 
  Unitary representations of the Virasoro and super-Virasoro algebras. 
  \textit{Comm. Math. Phys.} \textbf{103} (1986), 105--119.

\bibitem[G82]{G}
  R.L. Griess, 
  The friendly giant. 
  \textit{Invent. Math.} \textbf{69} (1982), 1--102.

\bibitem[G98a]{GrO}
  R.L. Griess, 
  A vertex operator algebra related to $E_8$ with automorphism group ${\rm O}\sp +(10,2)$.
  \textit{The Monster and Lie algebras} (Columbus, OH, 1996),  43--58, Ohio State
  Univ. Math. Res. Inst. Publ., 7, de Gruyter, Berlin, 1998.

\bibitem[G98b]{Gr}
  R.L. Griess, 
  Twelve Sporadic Groups. 
  Springer Verlag, 1998.



\bibitem[GL11]{GL}
  R.L. Griess and C. H. Lam, $EE_8$ lattices and dihedral groups.
  \textit{Pure and Applied Math Quarterly} (special issue for Jacques Tits),
  \textbf{7} (2011), no. 3, 621--743.

\bibitem[GL12a]{GLSDC}
  R.L. Griess and C.H. Lam, Diagonal lattices and rootless $EE_8$ pairs.
  \textit{J. Pure and Appl. Algebra} \textbf{216} (2012), no. 1, 154--169.

\bibitem[GL12b]{GL3}
  R.L. Griess and C.H. Lam,
  Moonshine paths for 3A  and 6A nodes of the extended $E_8$-diagram.
  \textit{J. Algebra} \textbf{379} (2013), 85--112.
  doi:10.1016/j.jalgebra.2012.12.019

\bibitem[H10]{Ho}
  G. H\"{o}hn, The group of symmetries of the shorter moonshine module.
  \textit{Abh. Math. Semin. Univ. Hambg.} \textbf{80} (2010), no. 2, 275--283. 

\bibitem[HLY12a]{HLY1}
  G. H\"{o}hn, C.H. Lam and H. Yamauchi,
  McKay's $E_7$ observation on the Baby Monster.
  \textit{Internat. Math. Res. Notices} (2012), doi:10.1093/imrn/rnr009.

\bibitem[HLY12b]{HLY2}
  G. H\"{o}hn, C.H. Lam and H. Yamauchi,
  McKay's $E_6$ observation on the largest Fischer group.
  \textit{Comm. Math. Physics} \textbf{310} Vol. 2 (2012), 329--365. 

\bibitem[HRS15]{HRS}
  J.I. Hall, F. Rehren and S. Shpectorov, 
  Primitive axial algebras of Jordan type.
  \textit{J. Algebra} \textbf{437} (2015), 79--115.

\bibitem[I76]{I}
  I. Martine Isaacs, 
  Character theory of finite groups. 
  Pure and Applied Mathematics, No. 69. Academic Press, New York-London, 1976.



\bibitem[KMY00]{KMY}
  M. Kitazume, M. Miyamoto and H. Yamada, 
  Ternary codes and vertex operator algebras. 
  \textit{J. Algebra} \textbf{223} (2000), 379--395.


\bibitem[LLY03]{LLY}
  C.H. Lam, N. Lam and H. Yamauchi, 
  Extension of unitary Virasoro vertex operator algebra by a simple module.
  \textit{Internat. Math. Res. Notices} \textbf{11} (2003), 577--611.

\bibitem[LM06]{LM}
  C.H. Lam and M. Miyamoto, 
  Niemeier lattices, Coxeter elements, and McKay's $E_8$-observation on the Monster simple group.
  \textit{Internat. Math. Res. Notices} Article ID 35967 (2006), 1--27.

\bibitem[LSY07]{LSY}
  C.H. Lam, S. Sakuma and H. Yamauchi, 
  Ising vectors and automorphism groups of commutant subalgebras related to root systems.
  \textit{Math. Z.} \textbf{255} vol.3 (2007), 597--626.

\bibitem[LS07]{LS}
  C.H. Lam and H. Shimakura,
  Ising vectors in the vertex operator algebra $V_\Lambda^+$ associated with 
  the Leech lattice $\Lambda$.
  \textit{Internat. Math. Res. Notices} \textbf{2007} Art. ID rnm132, 21 pp.

\bibitem[LY04]{LYd}
  C.H. Lam and H. Yamada, 
  Decomposition of the lattice vertex operator algebra $V_{\sqrt{2}A_l}$. 
  \textit{J. Algebra} \textbf{272} (2004), 614--624.

\bibitem[LYY07]{LYY1}
  C.H. Lam, H. Yamada and H. Yamauchi, 
  Vertex operator algebras, extended $E_8$-diagram, and McKay's observation on 
  the Monster simple group.
  \textit{Trans. Amer. Math. Soc.} \textbf{359} (2007), 4107--4123.

\bibitem[LYY05]{LYY2}
  C.H. Lam, H. Yamada and H. Yamauchi, 
  McKay's observation and vertex operator algebras generated by two conformal 
  vectors of central charge 1/2. 
  \textit{Internat. Math. Res. Papers} \textbf{3} (2005), 117--181.

\bibitem[LY07]{LY1}
  C.H. Lam and H. Yamauchi, 
  A characterization of the moonshine vertex operator algebra by means of Virasoro frames. 
  \textit{International Mathematics Research Notices} \textbf{2007} (2007), 
  article ID:rnm003, 9 pages.
  
\bibitem[LY14]{LY2}
  C.H. Lam and H. Yamauchi, 
  On 3-transposition groups generated by $\sigma$-involutions associated to 
  $c=4/5$ Virasoro vectors.
  \textit{J. Algebra}, \textbf{416} (2014), 84--121.

\bibitem[LY16]{LY3}
  C.H. Lam and H. Yamauchi, 
  3-dimensional Griess algebras and Miyamoto involutions. 
  Preprint.
 \href{http://arxiv.org/abs/1604.04470}{\tt arXiv:1604.04470}.

\bibitem[LSu]{LSo}
  C.\,H. Lam and C.\,S. Su, 
  Griess algebra generated by two $3A$-algebras  with a common axis.
  \textit{J. Math. Soc. Japan} \textbf{67} (2015), 453--476.

\bibitem[Ma01]{Ma1}
  A. Matsuo, 
  Norton's trace formulae for the Griess algebra of 
  a vertex operator algebra with larger symmetry.
  \textit{Commun. Math. Phys.} \textbf{224} (2001), 565--591.

\bibitem[Ma05]{Ma2}
  A. Matsuo, 
  $3$-transposition groups of symplectic type and vertex operator algebras. 
  {\it J. Math. Soc. Japan} \textbf{57} (2005), no. 3, 639--649. 
  \href{http://arxiv.org/abs/math/0311400}{\tt arXiv:math/0311400}.

\bibitem[MN93]{MN}
  W. Meyer and W. Neutsch, 
  Associative subalgebras of the Griess algebra. 
  \textit{J. Algebra} \textbf{158} (1993), 1--17.

\bibitem[Mi96]{M1}
  M. Miyamoto, 
  Griess algebras and conformal vectors in vertex operator algebras. 
  \textit{J. Algebra} \textbf{179} (1996), 528--548.

\bibitem[Mi01]{M2}
  M. Miyamoto, 
  $3$-state Potts model and automorphisms of vertex operator algebras of order $3$. 
  \textit{J. Algebra} \textbf{239} (2001), 56--76.

\bibitem[Mi03]{M3}
  M. Miyamoto, VOAs generated by two conformal vectors whose
  $\tau$-involutions generate $\mathrm{S}_3$. {\it J. Algebra} {\bf 268}
  (2003), 653--671.

\bibitem[Mi04]{M4}
  M. Miyamoto, 
  A new construction of the moonshine vertex operator algebra over the real number field. 
  \textit{Ann. Math.} \textbf{159} (2004), 535--596.

\bibitem[R15]{R}
  F. Rehren,  
  Linear idempotents in Matsuo algebras. 
  To appear in \textit{Indiana. Uni. Math. J.}  

\bibitem[S07]{S}
  S. Sakuma, 
  6-transposition property of $\tau$-involutions of vertex operator algebras.
  {\it Internat. Math. Res. Notices} {\bf 2007},  no. {\bf 9}, Art. ID rnm 030, 19 pp.

\bibitem[SY03]{SY}
  S. Sakuma and H. Yamauchi, 
  Vertex operator algebra with two Miyamoto involutions generating $S_3$. 
  \textit{J. Algebra} \textbf{267} (2003), 272--297.

\bibitem[W93]{W}
  W. Wang, Rationality of Virasoro vertex operator algebras.
  \textit{Internat. Math. Res. Notices} \textbf{71} (1993), 197--211.

\bibitem[Y14]{Y}
  H. Yamauchi, 
  Extended Griess algebras and Matsuo-Norton trace formulae. 
  \textit{Conformal Field Theory, Automorphic Forms and Related Topics}, Springer, 2014, 
  75--107.

\end{thebibliography}
